\documentclass[reqno]{amsart}
\usepackage{eurosym}
\usepackage{amsmath}
\usepackage{amsfonts}
\usepackage{amssymb}
\usepackage{amsfonts}
\usepackage{epsfig}
\usepackage{caption}
\usepackage{subcaption}
\usepackage{multirow}
\usepackage{csquotes}
\usepackage{stackengine}
\usepackage{makecell}
\stackMath

\newtheorem{theorem}{Theorem}

\newtheorem{corollary}[theorem]{Corollary}

\newtheorem{example}[theorem]{Example}

\newtheorem{lemma}[theorem]{Lemma}

\newtheorem{proposition}[theorem]{Proposition}
\newtheorem{remark}[theorem]{Remark}

\begin{document}

\title[Asymptotic condition numbers] {Asymptotic condition numbers for linear ordinary differential equations}
\author[S. Maset] {S.\ Maset \\
	Dipartimento di Matematica, Informatica e Geoscienze \\
	Universit\`{a} di Trieste \\
	maset@units.it}

\begin{abstract}
	{We are interested in the relative
		conditioning of the problem $y_0\mapsto \mathrm{e}^{tA}y_0$, i.e., the
		relative conditioning of the action of the matrix exponential $\mathrm{e}%
		^{tA}$ on a vector with respect to perturbations of this vector. The present paper is a study of the long-time behavior of this
		conditioning. In other words, we are interested in studying the propagation to the solution $y(t)$ of perturbations of the initial value for a linear ordinary differential equation $y^\prime(t)=Ay(t)$, by measuring these perturbations with relative errors. We introduce three condition numbers: the first considers a specific initial value and a specific direction of perturbation; the second considers a specific initial value and the worst case by varying the direction of perturbation; and the third considers the worst case by varying both  the initial value
		and the direction of perturbation. The long-time behaviors of
		these three condition numbers are studied.}
\end{abstract}

\maketitle



\section{Introduction}

 We are interested in understanding how a perturbation of the initial value $y_{0}$ of the linear $n-$dimensional Ordinary Differential Equation (ODE) 
\begin{equation}
	\left\{ 
	\begin{array}{l}
		y^{\prime }\left( t\right) =Ay\left( t\right) ,\ t\in\mathbb{R}, \\ 
		y\left( 0\right) =y_{0},%
	\end{array}%
	\right.  \label{ode}
\end{equation}%
where $A\in \mathbb{C}^{n\times n}$ and $y_0,y(t)\in\mathbb{C}^n$, is propagated to the solution $y(t)=\mathrm{e}^{tA}y_0$ of (\ref{ode}) over a long time interval. This perturbation, propagating along the solution, is measured by a \emph{relative error}. In other words, we study the \emph{relative conditioning} of the problem
\begin{equation}
	y_0\mapsto y(t)=\mathrm{e}^{tA}y_0  \label{due}
\end{equation}
for large time $t$.

The relative conditioning of the matrix exponential function, i.e., the relative conditioning of the problem $A\mapsto \mathrm{e}^{A}$, or the problem
\begin{equation}
	A\mapsto \mathrm{e}^{tA} \label{tre}
\end{equation}
involving the time $t$, has been extensively studied: see \cite{Levis1969}, \cite{van2006}, \cite{Kagstrom1977}, \cite{Van1977}, \cite{moler2003}, \cite{Mohy2008}, \cite{Zhu2008}, \cite{Mohy2011}, \cite{Ed} and \cite{ALMOHY2017}. In the study of time evolutions, an important aspect is to understand how the relative conditioning depends on $t$. Many of the papers cited above have examined this issue for the problem (\ref{tre}). For a normal matrix $A$, it is known that the relative condition number of (\ref{tre}) grows linearly with $t$. On the other hand, for a general matrix $A$, the exact order of growth with $t$ and its dependence on the matrix $A$ are not known: we have only polynomial lower bounds and exponential upper bounds in $t$ for the relative condition number.

It is also of interest to study the relative conditioning of the action of the matrix exponential on a vector. This is particularly important in the context of ODEs (\ref{ode}), where the solution is given by the action of the matrix exponential $\mathrm{e}^{tA}$ on the initial value $y_0$. In this context, we can consider the relative conditioning of the problems (\ref{due}) and 
\begin{equation}
	A \mapsto \mathrm{e}^{tA}y_0, \label{quattro}
\end{equation}
where, unlike the matrix exponential function problem (\ref{tre}), $y_0$ is involved.
Despite their importance, there has been little attention in the literature to the relative conditioning of these problems (\ref{due}) and (\ref{quattro}).

For the case of $A$ normal, an analysis of the relative conditioning, focused on time $t$, was carried out in \cite{FarooqMaset2021} for the problem (\ref{quattro}).

The relative conditioning of problem (\ref{due}) can be perceived, at first glance, as a trivial issue: (\ref{due}) is a linear problem and its condition number can be immediately determined and computed (see Subsection \ref{cn} below). This perception is especially reinforced when the relative conditioning of (\ref{due}) is compared with the relative conditioning of the non-linear problems (\ref{tre}) and (\ref{quattro}). However, the relative conditioning of (\ref{due}) ceases to be trivial once we want to understand how it depends on $t$: see points A and B at the beginning of Section \ref{AF} below.

Analyzing how the relative conditioning of (\ref{due}) depends on $t$ can fill a gap in our understanding of linear dynamics. In fact,
\emph{while it is well understood how the absolute conditioning of (\ref{due}) depends on $t$, i.e., how absolute errors due to perturbations of $y_0$ propagate to $y(t)$ (they are governed by the real part of the rightmost eigenvalues of $A$ for large $t$, and by the pseudospectra of $A$ for non-large $t$), how the relative conditioning of (\ref{due}) depends on $t$, i.e., how relative errors propagate, is far less understood, even for large $t$}. Such an analysis was carried out in \cite{Maset2018} for the simple case of $A$ normal. In the present paper, we carry out this analysis for a general ODE (\ref{ode}), by considering the relative conditioning for large $t$ and providing a detailed description of it.

\subsection{The condition numbers}\label{cn}

Suppose that the initial value $y_{0}\neq 0$ in (\ref{ode}) is perturbed to $%
\widetilde{y}_{0}$ and, as a consequence, the solution $y$ is perturbed to $\widetilde{y}$. Let $\left\Vert \ \cdot \ \right\Vert $ be an arbitrary vector norm on $\mathbb{C}^{n}$.  We introduce
the normwise relative error 
\begin{equation*}
\varepsilon :=\frac{\left\Vert \widetilde{y}_{0}-y_{0}\right\Vert }{%
\left\Vert y_{0}\right\Vert }
\end{equation*}
of $\widetilde{y}_{0}$ and the normwise relative error 
\begin{equation*}
\delta \left( t\right) :=\frac{\left\Vert \widetilde{y}\left( t\right)
-y\left( t\right) \right\Vert }{\left\Vert y\left( t\right) \right\Vert }  \label{deltat}
\end{equation*}%
of $\widetilde{y}(t)$. By writing
\begin{equation*}
\widetilde{y}_0=y_0+\varepsilon \Vert y_0\Vert \widehat{z}_0,
\end{equation*}
where $\widehat{z}_0\in\mathbb{C}^n$, with $\Vert \widehat{z%
}_0 \Vert=1$, is the \emph{direction of perturbation}, we obtain 
\begin{equation*}
\delta \left( t\right) =K\left( t,y_{0},\widehat{z}_{0}\right) \cdot  \label{magnification}
\varepsilon ,
\end{equation*}%
where
\begin{equation}
K\left( t,y_{0},\widehat{z}_{0}\right) :=\frac{\left\Vert \mathrm{e}^{tA}%
\widehat{z}_{0}\right\Vert }{\left\Vert \mathrm{e}^{tA}\widehat{y}%
_{0}\right\Vert } \label{Ktayz0},
\end{equation}%
with  $\widehat{y}_{0}:=\frac{y_{0}}{\left\Vert y_{0}\right\Vert }$ the \emph{normalized initial value}. The number $K\left( t,y_{0},\widehat{z}_{0}\right)$ is called the \emph{directional pointwise condition number} of the problem (\ref{due}): it is called ``directional'' because it depends on $\widehat{z}_0$, and ``pointwise'' because it depends on $y_0$.

Along with the condition number (\ref{Ktayz0}), we also introduce two other condition numbers:
\begin{itemize}
	\item the \emph{pointwise condition number} of the
problem (\ref{due}) given by 
\begin{equation}
K\left( t,y_{0}\right) :=\max\limits_{\substack{ \widehat{z}_0\in \mathbb{C%
}^{n}  \\ \Vert \widehat{z}_0\Vert =1}}K\left( t,y_{0},\widehat{z}%
_0\right) =\frac{\left\Vert \mathrm{e}^{tA}\right\Vert }{\left\Vert \mathrm{e%
}^{tA}\widehat{y}_{0}\right\Vert },  \label{KAy0}
\end{equation}%
where $\left\Vert \mathrm{e}^{tA}\right\Vert $ is the matrix norm of $%
\mathrm{e}^{tA}$ induced by the vector norm $\left\Vert \ \cdot \
\right\Vert $ (see \cite{Burgisser2013} for the definition of
condition number of a general problem, which corresponds to the pointwise condition number);
\item the \emph{global condition number} of the problem (\ref{due}) given by
\begin{equation}
K\left( t\right) :=\max\limits_{\substack{ y_0\in \mathbb{C}^{n}  \\ %
y_0\neq 0}}K\left( t,y_{0}\right) =\left\Vert \mathrm{e}^{tA}\right\Vert
\cdot \left\Vert \mathrm{e}^{-tA}\right\Vert=\kappa\left(\mathrm{e}%
^{tA}\right),  \label{KA0}
\end{equation}
which equals the standard condition number $\kappa\left(\mathrm{e}%
^{tA}\right)$ of the matrix $\mathrm{e}^{tA}$.
\end{itemize}

Observe that $K\left( t,y_{0}\right)$ is the worst $%
K(t,y_0,\widehat{z}_0)$ by varying $\widehat{%
z}_0$, and $K(t)$ is the worst $K(t,y_0)$ by varying $y_0$, i.e., the worst
$K(t,y_0,\widehat{z}_0)$ by varying both $y_0$ and $\widehat{z}_0$.

The paper \cite{Maset2018} studied the condition numbers (\ref{Ktayz0}%
), (\ref{KAy0}) and (\ref{KA0}) in the particular case of $A $ normal. The present paper studies the general case.

The aim of the present paper is to analyze the asymptotic (long-time) behavior of the three condition
numbers  $K(t,y_0,%
\widehat{z}_0)$, $K(t,y_0)$ and $%
K(t)$, i.e., their behavior as $t$ approaches infinity (becomes large).

\subsection{Plan of the paper}

Besides this introduction, the paper contains five sections and there are three appendices in the Supplementary Material.

Section 2 develops notions and notations for understanding the asymptotic forms of the three condition numbers. Section 3 analyzes the asymptotic behaviors of the condition numbers $K(t,y_0,%
\widehat{z}_0)$ and $K(t,y_0)$, by introducing the asymptotic condition numbers $K_\infty(t,y_0,%
\widehat{z}_0)$ and $K_\infty(t,y_0)$. Section 4 analyzes the asymptotic behavior of the condition number $K(t)$, by introducing the asymptotic condition numbers $K_\infty^+(t)$ and $K_\infty(t)$.  Section 5 introduces the Rightmost Last Generalized Eigenvector (RLGE) condition and summarizes the most important results. Sections 4 and 5 also contain numerical examples illustrating the results obtained. Conclusions are in Section 6.

The three appendices contain the more technical material. They should be consulted as needed while reading the paper, and read in full only by readers interested in the mathematical details.
 Appendix \ref{JCFsection} develops a suitable formula for the matrix exponential $\mathrm{e}^{tA}$ along with other fundamental material related to this formula. Appendix \ref{lastsection} investigates the properties of the key matrices $Q_{jl}(t)$ that determine the asymptotic behavior of the condition numbers.  Appendix \ref{sssQj0} analyzes the linear operators $Q_{j0}|_{U^e_j}(t)$, which are important for defining the asymptotic condition number $K_\infty(t)$ related to $K(t)$.

In what follows, we often refer to a \emph{generic case} for $y_0$ or $\widehat{z}_0$. By a generic case for an element $v$ of a finite-dimensional space $V$, we mean that $v$ satisfies a property which is not satisfied only on a subspace, or more generally a manifold,  $M$ of $V$ with $\mathrm{dim} M < \mathrm{dim} V$. Equivalently, if $v$ is drawn at random from $V$ (with respect to any distribution absolutely continuous with respect to the Lebesgue measure), the generic case holds with probability $1$.

\section{Asymptotic forms}\label{AF}

In the next Sections \ref{asymptoticbehavior} and \ref{asymptoticbehavior2}, we analyze the asymptotic behavior of the three condition
numbers $K(t,y_0,%
\widehat{z}_0)$, $K(t,y_0)$ and $%
K(t)$. The asymptotic forms of the condition numbers, called \emph{asymptotic condition numbers}, are obtained by inserting the asymptotic forms of $\mathrm{e}^{tA}%
\widehat{z}_{0}$ and $\mathrm{e}^{tA}%
\widehat{y}_{0}$ in (\ref{Ktayz0}), of $\mathrm{e}^{tA}$ and $\mathrm{e}^{tA}%
\widehat{y}_{0}$ in (\ref{KAy0}), and of $\mathrm{e}^{tA}$ and $\mathrm{e}^{-tA}$ in (\ref{KA0}).

	All this might suggest that the work is a straightforward development of well-established asymptotic results. However, it actually involves substantial mathematical difficulties. Specifically, the following two points should be remarked.
	
	\begin{itemize}
		\item [A] Although it is not difficult to identify the asymptotic forms of $\mathrm{e}^{tA}%
		\widehat{z}_{0}$, $\mathrm{e}^{tA}%
		\widehat{y}_{0}$, $\mathrm{e}^{tA}$ and $\mathrm{e}^{-tA}$, proving that they are asymptotic forms according to the definition (given below in Subsection \ref{firstsection})
		requires a certain mathematical effort: see Remark \ref{domterms} below.
		\item [B] If we are interested in defining asymptotic condition numbers for the problem (\ref{due}), the asymptotic condition numbers derived from (\ref{KAy0}) and (\ref{KA0}) as $t\rightarrow +\infty$ may not be appropriate, since they represent the \emph{asymptotic worst cases} of (\ref{Ktayz0}), by varying $\widehat{z}_0$ only and both $\widehat{z}_0$ and $y_0$, respectively. Instead, we may be concerned with the \emph{worst asymptotic cases} of (\ref{Ktayz0}).  In other words: do \enquote{asymptotic} and \enquote{worst} commute? Answering this interesting question also requires a certain mathematical effort: see Theorems \ref{wcy} and \ref{asworst2} below.
	\end{itemize}

In this section, we determine the asymptotic forms of $\mathrm{e}^{tA}$ and $\mathrm{e}^{tA}u$, where $u\in\mathbb{C}^n$. Some preliminary notations and notions need to be introduced.

\subsection{Notations $\sim$ and $\approx$} \label{firstsection}

In this subsection, we make precise what we mean by asymptotic form.

Let $f(t)$ and $g(t)$ be scalar, vector or matrix functions of $%
t\in \mathbb{R}$. For $g$ such that $g(t)\neq 0$ for $t$ in a neighborhood of $+\infty$, we write
\begin{equation}
f\left( t\right) \sim g\left( t\right),\ t\rightarrow +\infty,  \label{sim}
\end{equation}
when
\begin{equation*}
\lim\limits_{t\rightarrow +\infty }\frac{\Vert f(t)-g(t)\Vert }{\Vert
g(t)\Vert }=0.
\end{equation*}
In case of scalars, vectors and matrices, $\Vert\;\cdot\;\Vert$ denotes, respectively, the modulus, a vector norm and a matrix norm.

We interpret (\ref{sim}) as indicating that $g$ is an \emph{asymptotic form} of $f$.

Observe that (\ref{sim}) means that the relative error of $f$ with respect to its asymptotic form $g$ asymptotically vanishes. At a finite time, one may ask how dominant the asymptotic form is, i.e., how large the relative error of $f$ with respect to $g$ is. Therefore, we introduce the following notation. For $t\in\mathbb{R}$ such that $g(t)\neq 0$ and $\epsilon\geq 0$, we write 
\begin{equation}
f\left( t\right) \approx g\left( t\right)\textit{\ with precision\ }\epsilon  \label{approx}
\end{equation}
when
\begin{equation*}
\frac{\Vert f\left( t\right) -g\left( t\right) \Vert }{\Vert g\left(
t\right) \Vert }\leq \epsilon.
\end{equation*}
Observe that (\ref{approx}) means that the relative error of $f$ with respect to $g$ at the time $t$ is not larger than $\epsilon$. The notation (\ref{approx}) serves to quantify how dominant the asymptotic form $g$ of $f$ is at the finite time $t$. 

\begin{remark}
\label{Remark fg} Note that, for vector or matrix functions $f(t)$ and $%
g(t)$, 
$$
f\left( t\right)\sim g\left( t\right),\ t\rightarrow +\infty,
$$
implies 
$$
\Vert f\left( t\right) \Vert \sim \Vert g\left(
t\right) \Vert,\ t\rightarrow +\infty,
$$
and, for $\epsilon\geq 0$,
$$
f\left( t\right)\approx g\left( t\right)\text{\ with precision\ }\epsilon
$$
implies
$$
\Vert f\left( t\right) \Vert \approx \Vert g\left(
t\right) \Vert\text{\ with precision\ }\epsilon.
$$
This follows by
$$
\left\vert \left\Vert f(t)\right\Vert -\left\Vert g(t)\right\Vert\right\vert\leq \left\Vert f(t)-g(t)\right\Vert.
$$
\end{remark}

\subsection{Partition of the spectrum} \label{cc1}

In this subsection, we introduce a partition of the spectrum of $A$ suitable for the study of the asymptotic forms of $\mathrm{e}^{tA}$ and $\mathrm{e}^{tA}u$, and thus of the asymptotic forms of the condition numbers.

\begin{figure}[tbp]
\includegraphics[width=1\textwidth]{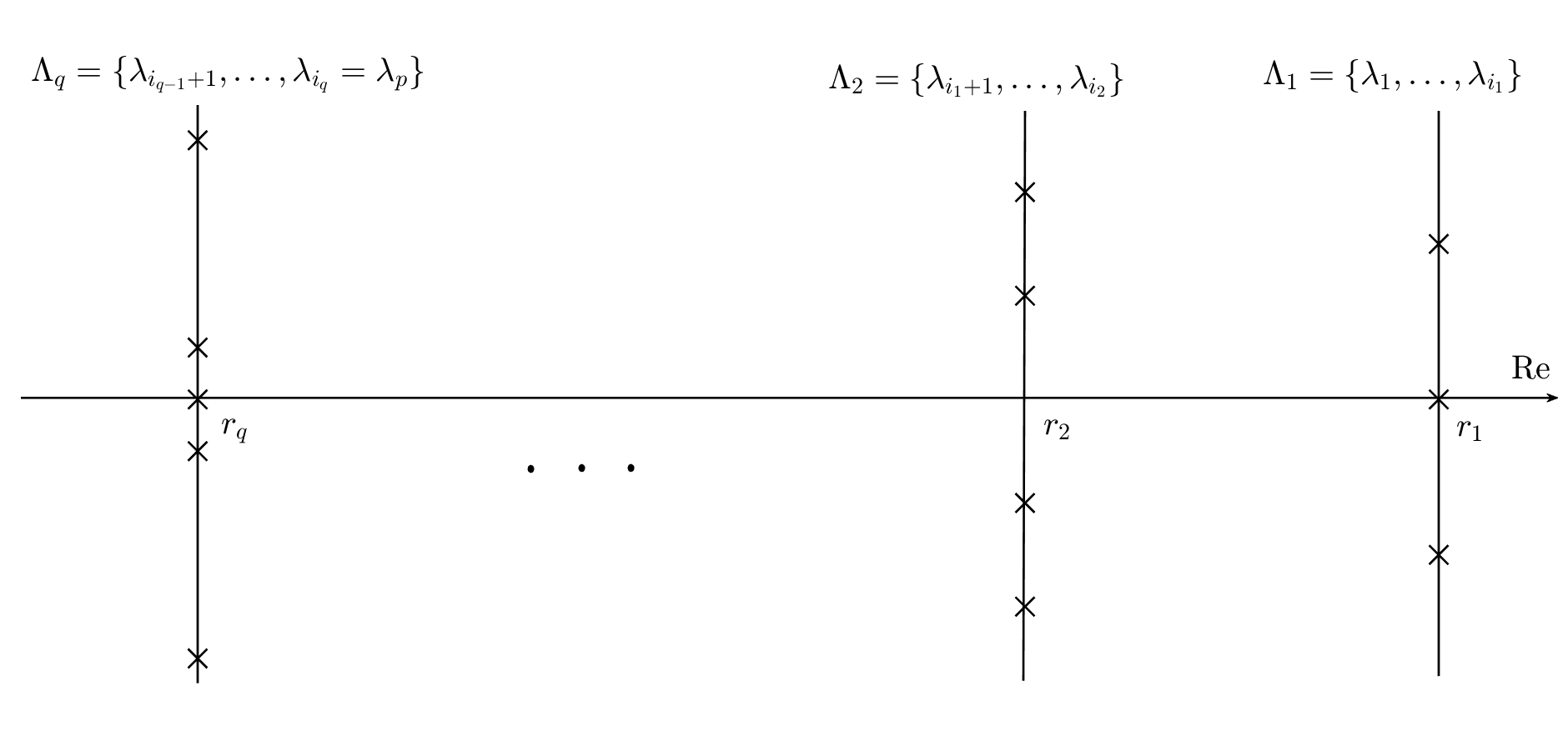}
\caption{Partition of the eigenvalues $\protect\lambda_1,\ldots,\protect%
\lambda_p$ of $A$ by decreasing real part. The eigenvalues are marked by ``x''.}
\label{Fig0} 
\end{figure}

The spectrum $\Lambda=\left\{\lambda _{1},\ldots ,\lambda _{p}\right\} $ of $A$,
where $\lambda _{1},\ldots ,\lambda _{p}$ are the distinct eigenvalues of $A$%
, is partitioned by decreasing real parts (see Figure \ref{Fig0}) in the
subsets $\Lambda_j$, $j\in\{1,\ldots,q\}$, given by 
\begin{eqnarray*}
&&\Lambda _{j}:=\{\lambda _{i_{j-1}+1},\lambda _{i_{j-1}+2},\ldots ,\lambda
_{i_{j}}\}  \label{Lambdaj} \\
&&\mathrm{Re}\left( \lambda _{i_{j-1}+1}\right) =\mathrm{Re}\left( \lambda
_{i_{j-1}+2}\right) =\cdots =\mathrm{Re}\left( \lambda _{i_j}\right) =r_{j}, 
\notag
\end{eqnarray*}
where the $q$ distinct real parts $r_j$, $j\in\{1,\ldots,q\}$, of the eigenvalues of $A$ satisfy $r_{1}>r_2>\cdots>r_{q}$. Observe that $\Lambda_1$ and $\Lambda_q$ are the sets of the rightmost and leftmost, respectively, eigenvalues of $A$. 

The rightmost real part $r_1$ is the \emph{spectral abscissa} of the matrix $A$. While the spectral abscissa plays a crucial role in the long-time behaviour of the matrix exponential $\mathrm{e}^{tA}$ and the absolute conditioning of the problem~(\ref{due}), it has no particular significance for the long-time behaviour of the relative conditioning. Instead, what matters here are (the asymptotic forms of) the matrix term $C(t)$ and the vector term $C(t,u)$ such that $\mathrm{e}^{tA}=\mathrm{e}^{r_1 t}C(t)$ and $\mathrm{e}^{tA}u=\mathrm{e}^{r_1 t}C(t,u)$. Therefore, we need to be more refined than simply stating that $\Vert \mathrm{e}^{tA}\Vert =O(\mathrm{e}^{r_1 t})$ as $t\to +\infty$.

\subsection{The matrix exponential} \label{cc2}
In this subsection, we introduce the tool for determining the asymptotic forms of $\mathrm{e}^{tA}$ and $\mathrm{e}^{tA}u$.  It is the formula (\ref{expexp}) below.

For the matrix exponential $\mathrm{e}^{tA}$, we have (see Appendix \ref{JCFsection})
\begin{equation}
\mathrm{e}^{tA}=\sum\limits_{i=1}^{p}\mathrm{e}^{\lambda_{i}t}\sum\limits_{l=0}^{m_{i}-1}\frac{t^{l}}{l!}P_{il}, \label{etAmp}
\end{equation}%
where $m_i$ is the ascent of the eigenvalue $\lambda _{i}$ and the matrices $P_{il}$ are given by
	\begin{equation}
	P_{il}=\frac{1}{2\pi\mathrm{i}}\int\limits_{\Gamma_i} (z-\lambda_i)^l (zI-A)^{-1} dz=(A-\lambda_iI)^lP_{i0}, \label{PilS11}
\end{equation}
with $\Gamma_i$ a positively oriented simple closed contour in the complex plane enclosing $\lambda_i$ and excluding the other eigenvalues of $A$. In particular, the matrix $P_{i0}$ is the projection onto the generalized eigenspace corresponding to $\lambda_i$. Here and in what follows, $\mathrm{i}$ denotes the imaginary unit (not to be confused with the index $i$).

The matrices $P_{il}$ are introduced, expressed and analyzed in Appendix \ref{JCFsection} using the Jordan Canonical Form (JCF) of $A$. The JCF is used for the theoretical analysis of the matrices $P_{il}$. For their computation, or for the computation of their action on a vector, formula (\ref{PilS11}) should be preferred (see \cite{DiNapoliPolizziSaad} for how this can be done) and a Schur decomposition may be used to control round-off errors.

By using in (\ref{etAmp}) the partition of the spectrum of $A$ given above, we obtain
\begin{equation}
\mathrm{e}^{tA}=\sum\limits_{j=1}^{q}\mathrm{e}^{r_{j}t}\sum\limits_{\lambda
_{i}\in \Lambda _{j}}\mathrm{e}^{\mathrm{i}\omega
_{i}t}\sum\limits_{l=0}^{m_{i}-1}\frac{t^{l}}{l!}P_{il},  \label{exptA0}
\end{equation}%
where $\omega _{i}$ denotes the imaginary part of the eigenvalue $\lambda _{i}$.

By exchanging the two inner sums, we can rewrite (\ref{exptA0})  as
\begin{equation}
\mathrm{e}^{tA}=\sum\limits_{j=1}^{q}\mathrm{e}^{r_{j}t}\sum%
\limits_{l=0}^{L_j}\frac{t^{l}}{l!}Q_{jl}(t),  \label{expexp}
\end{equation} 
where 
\begin{equation*}
L_j:=\max_{\lambda _{i}\in \Lambda _{j}}m_{i}-1,\ j\in\{1,\ldots,q\}, 
\end{equation*}
and 
\begin{equation}
Q_{jl}(t):=\sum\limits_{\substack{ \lambda _{i}\in \Lambda _j  \\ %
m_{i}\geq l+1}}\mathrm{e}^{\mathrm{i}\omega _{i}t}P_{il},\ j\in\{1,\ldots,q\}%
\text{\ and\ }l\in\{0,\ldots,L_j\}.  \label{Qjl}
\end{equation}
The formula (\ref{expexp}) is used to determine the asymptotic forms in (\ref{Ktayz0}), (\ref{KAy0}) and (\ref{KA0}).

\subsection{The asymptotic form of $\mathrm{e}^{tA}$} \label{ssexp}

By looking at the formula (\ref{expexp}), we immediately identify
\begin{equation}
	\mathrm{e}^{r_1t}\frac{t^{L_1}}{L_{1}!}Q_{1L_1}(t) \label{af}
\end{equation}
as the asymptotic form of $\mathrm{e}^{tA}$, since it is the dominant term as $t\rightarrow +\infty$. Next proposition makes precise that (\ref{af}) is indeed an asymptotic form of $\mathrm{e}^{tA}$ according to the definition (\ref{sim}) and quantifies its dominance at a finite time.

\begin{proposition}
	\label{Lemma2} We have 
	\begin{equation*}
		\mathrm{e}^{tA}\approx \mathrm{e}^{r_1t}\frac{t^{L_1}}{L_{1}!}Q_{1L_1}(t)
	\end{equation*}%
	with precision 
	\begin{equation}
		\epsilon(t):=\sum\limits_{l=0}^{L_1 -1} \frac{L_1 !}{l!}t^{l-L_1}\frac{%
			\left\Vert Q_{1l}(t)\right\Vert}{\left\Vert Q_{1L_1}(t)\right\Vert }%
		+\sum\limits_{j=2}^{q}\mathrm{e}^{\left( r_{j}-r_{1}\right)
			t}\sum\limits_{l=0}^{L_j}\frac{L_1!}{l!}t^{l-L_1}\frac{\left\Vert
			Q_{jl}(t)\right\Vert }{\left\Vert Q_{1L_1}(t)\right\Vert }.  \label{eps}
	\end{equation}
	Moreover, we have
	\begin{equation}
		\epsilon(t)\rightarrow 0,\ t\rightarrow +\infty, \label{epsilon}
	\end{equation}
	and then
	\begin{equation*}
		\mathrm{e}^{tA}\sim \mathrm{e}^{r_1t}\frac{t^{L_1}}{L_{1}!} Q_{1L_1}(t),\
		t\rightarrow +\infty.
	\end{equation*}
\end{proposition}

\begin{proof}
	By (\ref{expexp}), we can write 
	\begin{eqnarray*}
		\mathrm{e}^{tA}
		=\mathrm{e}^{r_{1}t}\frac{t^{L_1}}{L_1!}Q_{1L_1}(t)+\mathrm{e}%
		^{r_{1}t}\sum\limits_{l=0}^{L_1-1}\frac{t^{l}}{l!}Q_{1l}(t)+\sum%
		\limits_{j=2}^{q}\mathrm{e}^{r_{j}t}\sum\limits_{l=0}^{L_j}\frac{t^{l}}{l!}%
		Q_{jl}(t).
	\end{eqnarray*}
	The first part of the proposition regarding $\approx$ follows.
	
	Since, in (\ref{eps}), 
	\begin{eqnarray}
	&&\sup\limits_{t\in \mathbb{R}} \left\Vert Q_{1l} (t)\right\Vert<+\infty,\ l\in\{0,\ldots,L_1-1\}, \notag \\
	&&\sup\limits_{t\in \mathbb{R}} \left\Vert Q_{jl} (t)\right\Vert<+\infty,\ j\in\{2,\ldots,q\}\text{\ \ and\ \ }l\in\{0,\ldots,L_j\}, \notag \\
	&&\inf\limits_{t\in \mathbb{R}} \left\Vert Q_{1L_1} (t)\right\Vert>0 \notag\\
	\label{three}
	\end{eqnarray}
	(see 1) and 3) in Proposition \ref{lemmaQjl} of Appendix \ref{lastsection}), we obtain (\ref{epsilon}) and then the second part regarding $\sim$ follows. 
\end{proof}

\begin{remark}
	\label{etAAA}
	If $A$ is diagonalizable, then we have
	\begin{equation*}
		\mathrm{e}^{tA}\approx \mathrm{e}^{r_1t}Q_{10}(t)
	\end{equation*}
	with precision 
	\begin{equation*}
		\epsilon(t)=\sum\limits_{j=2}^{q}\mathrm{e}^{\left( r_{j}-r_{1}\right)
			t}\frac{\left\Vert Q_{j0}(t)\right\Vert }{\left\Vert Q_{10}(t)\right\Vert }
	\end{equation*}
	and
	\begin{equation*}
		\mathrm{e}^{tA}\sim \mathrm{e}^{r_1t}Q_{10}(t),\ t\rightarrow +\infty,
	\end{equation*}
	where
	\begin{equation*}
		Q_{j0}(t)=\sum\limits_{\lambda _{i}\in \Lambda _j}\mathrm{e}^{\mathrm{i}\omega _{i}t}P_{i0},\ j\in\{1,\ldots,q\},
	\end{equation*}
	with $P_{i0}$ the projection onto the eigenspace corresponding to $\lambda_i$.
\end{remark}

\subsection{Notations for the asymptotic form of $\mathrm{e}^{tA}u$} \label{asfetAu}
The following notations are crucial for determining the asymptotic form of $\mathrm{e}^{tA}u$. The sets $\Lambda_j$ and the matrices $P_{il}$ introduced in the previous Subsections \ref{cc1} and \ref{cc2}, respectively, are involved in defining these notations.

For $j\in\{1,\ldots,q\}$ and $u\in\mathbb{C}^{n}$, we define $\Lambda_j(u)$ and $L_j(u)$.
\begin{itemize}
	\item Let
	\begin{equation*}
		\Lambda_j(u) := \{\lambda_i \in \Lambda_j : P_{i0}u \neq 0\},
	\end{equation*}
In other words, $\Lambda_j(u)$ is obtained from $\Lambda_j$ by including only the eigenvalues $\lambda_i$ for which $u$ has a non-zero projection onto the generalized eigenspace corresponding to $\lambda_i$. 
	
	\item When $\Lambda_j(u) \neq \emptyset$, let
	\begin{equation*}
		L_j(u) := \max\{l_i(u) : \lambda_i \in \Lambda_j(u)\},
	\end{equation*}
	where (see  Appendix \ref{JCFsection})
	\begin{equation*}
		l_i(u) := \max\{ l \in \{0, \ldots, m_i - 1\} : P_{il}u \neq 0\}.
	\end{equation*}
	
	In terms of the JCF  of $A$, $l_i(u) + 1$ is the maximum $k$ such that  $u$ has  non-zero component along the $k$-th generalized eigenvector of some Jordan chain corresponding to the eigenvalue $\lambda_i$. Thus, $L_j(u) + 1$ is the maximum $k$ such that $u$ has non-zero component along the $k$-th generalized eigenvector of some Jordan chain corresponding to an eigenvalue in $\Lambda_j(u)$.
\end{itemize}

\subsubsection{Indices of dominance} \label{iod}
For $u\in \mathbb{C}^{n}\setminus\{0\}$, we define the indices of dominance $j(u)$ and $L(u)$.
\begin{itemize}
	\item Let 
	\begin{equation*}
		j(u):=\min\{j\in \left\{ 1,\ldots,q\right\}: \Lambda_j(u)\neq \emptyset\}.
	\end{equation*}
	In other words, $\Lambda_{j(u)}$ is the rightmost set $\Lambda_j$ such that $u$ has, for some $\lambda_i\in\Lambda_j$, non-zero projection onto the generalized eigenspace corresponding to $\lambda_i$. We call $j(u)$ the  \emph{primary index of dominance} of $u$.
	\item Let 
	\begin{equation*}
		L(u):=L_{j(u)}(u).
	\end{equation*}
	
	In terms of the JCF of $A$, $L(u)+1$ is the maximum $k$ such that $u$ has non-zero component along the $k$-th generalized eigenvector of some Jordan chain corresponding to an eigenvalue in $\Lambda_{j(u)}(u)$.  We call $L(u)$  the \emph{secondary index of dominance} of $u$.
\end{itemize}

In Table \ref{table1}, we summarize the notation introduced in this subsection.

\begin{table}
	\begin{tabular}{|l|l|}
		\hline
		Notation & Description \\
		\hline
		$\Lambda_j(u)$ & \makecell[l]{set of the eigenvalues in $\Lambda_j$ with non-zero projection of $u$\\ onto the generalized eigenspace}\\
		\hline
		$l_i(u)$ & maximum $l$ such that $P_{il}u\neq 0$\\
		\hline
		$ L_j(u)$ & maximum $l_i(u)$ of eigenvalues $\lambda_i\in \Lambda_j(u)$\\ 
		\hline
		\makecell[l]{Primary index\\
			of dominance $j(u)$} & $\Lambda_{j(u)}$  is the rightmost $\Lambda_j(u)\neq\emptyset$\\
		\hline
		\makecell[l]{Secondary index\\
			of dominance $L(u)$} & $L(u)$ is $L_j(u)$ for $j=j(u)$ \\
		\hline
	\end{tabular}
\caption{$\Lambda_j(u)$, $l_i(u)$, $L_j(u)$, $j(u)$ and $L(u)$: brief description.}
\label{table1}
\end{table}

\subsection{The asymptotic form of $\mathrm{e}^{tA}u$} \label{ssexpu}

Let $u\in\mathbb{C}^n\setminus\{0\}$. When we use the formula (\ref{expexp}) for $\mathrm{e}^{tA}$ in $\mathrm{e}^{tA} u$, we obtain
\begin{equation}
\mathrm{e}^{tA}u=\sum\limits_{j=1}^{q}\mathrm{e}^{r_{j}t}\sum%
\limits_{l=0}^{L_j}\frac{t^{l}}{l!}Q_{jl}(t)u. \label{expexpu}
\end{equation}
Since, in (\ref{expexpu}),
\begin{eqnarray}
	&&Q_{jl}(t)u=0,\ j\in\{1,\ldots,j(u)-1\}\text{\ \ and\ \ }l\in\{0,\ldots,L_j\}\notag\\
	&&Q_{j(u)l}(t)u=0,\ l\in\{L(u)+1,\ldots,L_{j(u)}\},\notag \\
	\label{Qjlzero}
\end{eqnarray}
(see 1) and 2) in Proposition \ref{lemma2Qjl} of Appendix \ref{lastsection}), we identify
\begin{equation}
\mathrm{e}^{r_{j\left( u\right) }t}\frac{t^{L(u)}}{%
L(u)!}Q_{j(u)L(u)}(t)u   \label{afu}
\end{equation}
as the asymptotic form of $\mathrm{e}^{tA} u$, since it is the dominant term as $t\rightarrow +\infty$. Next proposition states precisely that (\ref{afu}) is indeed an asymptotic form according to the definition (\ref{sim}) and quantifies its dominance at a finite time.
\begin{proposition}
\label{Lemma1}For $u\in \mathbb{C}^{n}\setminus\{0\}$, we have 
\begin{equation*}
\mathrm{e}^{tA}u\approx \mathrm{e}^{r_{j\left( u\right) }t}\frac{t^{L(u)}}{%
L(u)!}Q_{j(u)L(u)}(t)u
\end{equation*}%
with precision 
\begin{eqnarray}
\epsilon(t,u)&:=&\sum\limits_{l=0}^{L\left( u\right) -1} \frac{L\left(
u\right) !}{l!}t^{l-L\left( u\right) }\frac{\left\Vert
Q_{j(u)l}(t)u\right\Vert}{\left\Vert Q_{j(u)L(u)}(t)u\right\Vert } \notag \\
&&+\sum\limits_{j=j\left( u\right) +1}^{q}\mathrm{e}^{\left( r_{j}-r_{j\left( u\right) }\right)
t}\sum\limits_{l=0}^{L_j}\frac{L(u)!}{l!}t^{l-L(u)}\frac{\left\Vert
Q_{jl}(t)u\right\Vert }{\left\Vert Q_{j(u)L(u)}(t)u\right\Vert}. \notag \\
\label{epsu}
\end{eqnarray}
Moreover, we have
\begin{equation}
\epsilon(t,u)\rightarrow 0,\ t\rightarrow +\infty,  \label{epsilonu}
\end{equation}
and then
\begin{equation*}
\mathrm{e}^{tA}u\sim \mathrm{e}^{r_{j\left( u\right) }t}\frac{t^{L(u)}}{L(u)!%
}Q_{j(u)L(u)}(t)u,\ t\rightarrow +\infty.
\end{equation*}
\end{proposition}

\begin{proof}
By (\ref{expexpu}) and (\ref{Qjlzero}), we can write 
\begin{eqnarray*}
\mathrm{e}^{tA}u&=&\sum\limits_{j=j(u)}^{q}%
\mathrm{e}^{r_{j}t}\sum\limits_{l=0}^{L_j}\frac{t^{l}}{l!}Q_{jl}(t)u \\
&=&\mathrm{e}^{r_{j(u)}t}\frac{t^{L(u)}}{L(u)!}Q_{j(u)L(u)}(t)u+\mathrm{e}%
^{r_{j(u)}t}\sum\limits_{l=0}^{L(u)-1}\frac{t^{l}}{l!}Q_{j(u)l}(t)u\\
&&+\sum\limits_{j=j(u)+1}^{q}%
\mathrm{e}^{r_{j}t}\sum\limits_{l=0}^{L_j}\frac{t^{l}}{l!}Q_{jl}(t)u.
\end{eqnarray*}%
The first part of the proposition regarding $\approx$ follows.

Since, in (\ref{epsu}), 
\begin{eqnarray}
	&&\sup\limits_{t\in \mathbb{R}} \left\Vert Q_{j(u)l} (t)u\right\Vert<+\infty,\ l\in\{0,\ldots,L(u)-1\}, \notag \\
	&&\sup\limits_{t\in \mathbb{R}} \left\Vert Q_{jl} (t)u\right\Vert<+\infty,\ j\in\{j(u)+1,\ldots,q\}\text{\ \ and\ \ }l\in\{0,\ldots,L_j\}, \notag \\
	&&\inf\limits_{t\in \mathbb{R}} \left\Vert Q_{j(u)L(u)} (t)u\right\Vert>0 \notag\\
	\label{threeu}
\end{eqnarray}
(see 2) and 4) in Proposition \ref{lemmaQjl} of Appendix \ref{lastsection}), we obtain (\ref{epsilonu}) and then the second part regarding $\sim$ follows.
\end{proof}

\begin{remark} \label{etAAAu}
	If $A$ is diagonalizable, then we have
	\begin{equation*}
		\mathrm{e}^{tA}u\approx \mathrm{e}^{r_{j\left( u\right) }t}Q_{j(u)0}(t)u
	\end{equation*}
	with precision
	\begin{eqnarray*}
		\epsilon(t,u)=\sum\limits_{j=j\left( u\right) +1}^q\mathrm{e}^{\left( r_{j}-r_{j\left( u\right)
			}\right) t}\frac{\left\Vert Q_{j0}(t)u\right\Vert }{\left\Vert
			Q_{j(u)0}(t)u\right\Vert}
	\end{eqnarray*}
	and
	\begin{equation*}
		\mathrm{e}^{tA}u\sim \mathrm{e}^{r_{j\left( u\right) }t}Q_{j(u)0}(t)u,\ t\rightarrow +\infty.
	\end{equation*}
\end{remark}

The asymptotic form (\ref{afu}) of $\mathrm{e}^{tA}u$ is determined by the primary and secondary indices of dominance of $u$. In particular, the smaller the primary index of dominance $j(u)$, the higher the order of the asymptotic  form (\ref{afu}), and for a fixed primary index, the larger the secondary index of dominance $L(u)$, the higher the order of (\ref{afu}).

\begin{remark} \label{domterms}
Recall point A at the beginning of Section \ref{AF}. The contents of Propositions \ref{Lemma2} and \ref{Lemma1} are deeper than what might be apparent at first glance. We are not simply stating the obvious facts that (\ref{af}) and (\ref{afu}) are the dominant terms in $\mathrm{e}^{tA}$ and $\mathrm{e}^{tA}u$, respectively. Rather, we are asserting that they are the asymptotic forms of $\mathrm{e}^{tA}$ and $\mathrm{e}^{tA}u$ according to the definition (\ref{sim}). In other words, the relative errors of $\mathrm{e}^{tA}$ and $\mathrm{e}^{tA}u$ with respect to the asymptotic forms (\ref{af}) and (\ref{afu}), respectively, approach zero asymptotically. Proving these deeper conclusions is complicated: to obtain (\ref{epsilon}) and (\ref{epsilonu}), one needs the key facts
\begin{equation}
\inf\limits_{t\in\mathbb{R}}\Vert Q_{1L_1}(t)\Vert>0\text{\ and\ \ }\inf\limits_{t\in\mathbb{R}}\Vert Q_{j(u)L(u)}(t)u\Vert>0 \label{kfacts}
\end{equation}
in (\ref{three}) and (\ref{threeu}) preventing possible ``cancellations'' in $ Q_{1L_1}(t)$ and $Q_{j(u)L(u)}(t)u$. These key facts are included in Proposition \ref{lemmaQjl} of Appendix \ref{lastsection}  and follows from Propositions \ref{PropPilu} (linear independence of the vectors $P_{il}u$), Proposition \ref{PropPil} (linear independence of the matrices $P_{il}$) and Lemma \ref{Lemma0} (condition for time-dependent linear combinations away from zero) in Appendix \ref{JCFsection}, developed through the work in that appendix. Considering the case in which there are incommensurable imaginary parts $\omega_i$ in (\ref{Qjl}), one realizes that (\ref{kfacts}) cannot have a simple proof.

The key facts (\ref{kfacts}) are also important in the qualitative study of the asymptotic condition numbers: see below points 1-a) of Remarks \ref{KtAyz} and \ref{KtAy}, Remark \ref{KtA} and the proofs of Theorems \ref{wcy} and \ref{asworst2}.

\end{remark}

\section{The asymptotic behaviors of $K(t,y_0,\widehat{z}_0)$ and $K(t,y_0)$} \label{asymptoticbehavior}

This section and the next one form the core of our analysis. Building on the preparatory work of the previous section, we can now easily describe the asymptotic behavior of the condition numbers.

In particular, in this section we study the asymptotic behavior of the condition numbers $K(t,y_0,\widehat{z}_0)$ and $K(t,y_0)$. Their asymptotic forms are the asymptotic condition numbers $K_\infty(t,y_0,\widehat{z}_0)$ and $K_\infty(t,y_0)$. We also show that $K_\infty(t,y_0)$ coincides with the worst $K_\infty(t,y_0,\widehat{z}_0)$, by varying $\widehat{z}_0$.

\subsection{The asymptotic condition number $K_\infty\left( t,y_{0},\widehat{z}%
_{0}\right)$}  \label{ssccfirst}

We set 
\begin{equation*}
j^{\ast }:=j\left( \widehat{y}_{0}\right) =j\left( y_{0}\right)\text{\ \
and\ \ }L^{\ast }:=L\left( \widehat{y}_{0}\right) =L\left( y_{0}\right)
\end{equation*}
as well as 
\begin{equation*}
j^{\ast \ast }:=j\left( \widehat{z}_{0}\right)\text{\ \ and\ \ } L^{\ast
\ast }:=L\left( \widehat{z}_{0}\right),
\end{equation*}
i.e. $j^\ast$ and $j^{\ast\ast}$ are the primary indices of dominance, and $L^\ast$ and $L^{\ast\ast}$ are the secondary indices of dominance, of $y_0$ and $\widehat{z}_0$, respectively.

The next theorem describes the asymptotic form of $%
K\left( t,y_{0},\widehat{z}_{0}\right) $.

\begin{theorem}
\label{ThmKyz} We have 
\begin{equation*}
K\left( t,y_{0},\widehat{z}_{0}\right) \approx K_\infty\left( t,y_{0},%
\widehat{z}_{0}\right), 
\end{equation*}
where 
\begin{equation*}
K_\infty\left( t,y_{0},\widehat{z}_{0}\right):=\frac{L^{\ast}!}{%
L^{\ast\ast }!}\mathrm{e}^{\left( r_{j^{\ast \ast }}-r_{j^{\ast }}\right)
t}t^{L^{\ast \ast }-L^{\ast }}\frac{\left\Vert
Q_{j^{\ast\ast}L^{\ast\ast}}(t)\widehat{z}_{0}\right\Vert }{\left\Vert
Q_{j^{\ast}L^{\ast}}(t)\widehat{y}_{0}\right\Vert },  \label{sim11}
\end{equation*}
with precision 
\begin{equation*}
\frac{\epsilon(t,\widehat{z}_0)+\epsilon(t,\widehat{y}_0)}{%
1-\epsilon(t,\widehat{y}_0)},
\end{equation*}
whenever $\epsilon(t,\widehat{y}_0)<1$ ($\epsilon(t,\widehat{z}_0)$ and $\epsilon(t,\widehat{y}_0)$ are defined in Proposition \ref{Lemma1}). Moreover, we have 
\begin{equation*}
K\left( t,y_{0},\widehat{z}_{0}\right) \sim K_\infty\left( t,y_{0},%
\widehat{z}_{0}\right),\ t\rightarrow +\infty.  \label{simKyz}
\end{equation*}
\end{theorem}

\begin{proof}
Proposition \ref{Lemma1} states the asymptotic forms of $\mathrm{e}^{tA}\widehat{z}_{0}$ and $\mathrm{e}^{tA}\widehat{y}_{0}$ in (\ref{Ktayz0}), and quantifies how dominant they are at a finite time $t$: we have
\begin{equation*}
\mathrm{e}^{tA} \widehat{z}_{0}\approx \mathrm{e}^{r_{j^{\ast \ast }}t}\frac{%
t^{L^{\ast \ast }}}{L^{\ast \ast }!}Q_{j^{\ast\ast}L^{\ast\ast}}(t)\widehat{z%
}_{0}
\end{equation*}
with precision $\epsilon(t,\widehat{z}_0)$ and 
\begin{equation*}
\mathrm{e}^{tA} \widehat{y}_{0}\approx \mathrm{e}^{r_{j^{\ast }}t}\frac{%
t^{L^{\ast }}}{L^{\ast }!}Q_{j^{\ast}L^{\ast}}(t)\widehat{y}_{0}
\end{equation*}
with precision $\epsilon(t,\widehat{y}_0)$. The first part of the theorem regarding $\approx$ follows by using Remark \ref{Remark fg} and by bounding the relative error of the ratio $\frac{\Vert \mathrm{e}^{tA}\widehat{z}_0\Vert }{\Vert \mathrm{e}^{tA}\widehat{y}_0\Vert }$
in terms of the bounds  $\epsilon(t,\widehat{z}_0)$ and  $\epsilon(t,\widehat{y}_0)$ of the relative errors of $\Vert \mathrm{e}^{tA}\widehat{z}_0\Vert$ and $\Vert \mathrm{e}^{tA}\widehat{y}_0\Vert $
, respectively. These are relative errors with respect to the norms of the asymptotic forms. 
The second part regarding $\sim$ follows since $\epsilon(t,\widehat{y}%
_0)\rightarrow 0$ and $\epsilon(t,\widehat{z}_0)\rightarrow 0$, $%
t\rightarrow +\infty$.
\end{proof}

We define the function 
$$
t\rightarrow K_\infty\left(t,y_0,\widehat{z}_0\right),\ t\in\mathbb{R},
$$
as the \emph{asymptotic directional pointwise condition number} of the problem (\ref{due}%
).

\begin{remark}
\label{KtAyz} \quad

\begin{itemize}
\item[1.] The asymptotic directional pointwise condition number  $K_\infty\left(t,y_0,\widehat{z}_0\right)$:

\begin{itemize}
\item [a)] is bounded and away from zero, as $t$ varies, if $j^\ast=j^{\ast\ast}$ and $L^\ast=L^{\ast\ast}$: this follows by
$$
\sup\limits_{t\in \mathbb{R}} \left\Vert Q_{j^{\ast\ast}L^{\ast\ast}} (t)\widehat{z}_0\right\Vert<+\infty\text{\ \ and\ \ }\inf\limits_{t\in \mathbb{R}} \left\Vert Q_{j^{\ast}L^{\ast}} (t)\widehat{y}_0\right\Vert>0
$$
and
$$
\inf\limits_{t\in \mathbb{R}} \left\Vert Q_{j^{\ast\ast}L^{\ast\ast}} (t)\widehat{z}_0\right\Vert>0\text{\ \ and\ \ }\sup\limits_{t\in \mathbb{R}} \left\Vert Q_{j^{\ast}L^{\ast}} (t)\widehat{y}_0\right\Vert<+\infty
$$

(recall 2) and 4) in Proposition \ref{lemmaQjl} of Appendix \ref{lastsection}) ;

\item [b)] decays polynomially to zero, as $t\rightarrow +\infty$, if $j^\ast=j^{\ast\ast}$ and $L^\ast>L^{\ast\ast}$;

\item [c)] diverges polynomially to infinity, as $t\rightarrow +\infty$, if $j^\ast=j^{\ast\ast}$ and $L^\ast<L^{\ast\ast}$;

\item [d)] decays exponentially to zero, as $t\rightarrow +\infty$, if $j^\ast<j^{\ast\ast}$;

\item [e)] diverges exponentially to infinity, as $t\rightarrow +\infty$, if $j^\ast>j^{\ast\ast}$.
\end{itemize}
Hence, whether $K_\infty\left(t,y_0,\widehat{z}_0\right)$ decreases to zero, diverges to infinity, or exhibits different behavior depends on which between $y_0$ and $\widehat{z}_0$ is more dominant. It decreases to zero if $y_0$ is more dominant and diverges to infinity if $\widehat{z}_0$ is more dominant. The dominance is determined by which of $y_0$ and $\widehat{z}_0$  possesses the smaller primary index of dominance, or, in the case of equal primary indices, the larger secondary index of dominance.

\item[2.] The case $j^{\ast
}=j^{\ast \ast }=1$ and $L^{\ast}=L^{\ast \ast }=L_1$ is generic for $y_0$ and $\widehat{z}_0$. In this generic case, we have 
\begin{equation*}
K_\infty\left( t,y_{0},\widehat{z}_{0}\right) = \frac{\left\Vert
Q_{1L_1}(t)\widehat{z}_{0}\right\Vert }{\left\Vert Q_{1L_1}(t)\widehat{y}%
_{0}\right\Vert }.
\end{equation*}

\end{itemize}
\end{remark}

\subsection{The asymptotic condition number $K_\infty\left( t,y_{0}\right)$}

\bigskip

The next theorem describes the asymptotic form of  $%
K\left( t,y_{0}\right) $, worst $K(t,y_0,\widehat{z}_0)$ by varying $\widehat{z}_0$. 

\begin{theorem}
\label{ThmKy} We have 
\begin{equation*}
K\left( t,y_{0}\right) \approx K_\infty\left( t,y_{0}\right),
\end{equation*}
where 
\begin{equation*}
K_\infty\left( t,y_{0}\right) := \frac{L^{\ast }!}{L_{1}!}\mathrm{e}%
^{\left( r_{1}-r_{j^{\ast }}\right) t}t^{L_{1}-L^{\ast }}\frac{\left\Vert
Q_{1L_{1}}(t)\right\Vert }{\left\Vert Q_{j^\ast L^\ast}(t) \widehat{y}%
_{0}\right\Vert },  \label{sim2}
\end{equation*}
with precision 
\begin{equation*}
\frac{\epsilon(t)+\epsilon(t,\widehat{y}_0)}{1-\epsilon(t,%
\widehat{y}_0)}
\end{equation*}
whenever $\epsilon(t,\widehat{y}_0)<1$ ($\epsilon(t)$ and $\epsilon(t,\widehat{y}_0)$ are defined in Propositions  \ref{Lemma2}  and  \ref{Lemma1}, respectively). Moreover, we have 
\begin{equation*}
K\left( t,y_{0}\right) \sim K_\infty\left( t,y_{0}\right),\ t\rightarrow
+\infty.  \label{simKy}
\end{equation*}
\end{theorem}

\begin{proof}
Proposition \ref{Lemma2} states the asymptotic form of $\mathrm{e}^{tA}$, and quantifies how dominant it is at a finite time $t$: we have
\begin{equation*}
\mathrm{e}^{tA}\approx \mathrm{e}^{r_1t}\frac{t^{L_1}}{L_1!} Q_{1L_1}(t)
\end{equation*}
with precision $\epsilon(t)$. Proposition \ref{Lemma1} states the asymptotic form of $\mathrm{e}^{tA}\widehat{y}_0$, and quantifies how it is dominant at a finite time $t$: we have
\begin{equation*}
\mathrm{e}^{tA} \widehat{y}_{0}\approx \mathrm{e}^{r_{j^{\ast }}t}\frac{%
t^{L^{\ast }}}{L^{\ast }!}Q_{j^\ast L^{\ast}}(t)\widehat{y}_{0}
\end{equation*}
with precision $\epsilon(t,\widehat{y}_0)$. The first part of the theorem regarding $\approx$ follows by using Remark \ref{Remark fg} and by bounding the relative error of the ratio $\frac{\Vert \mathrm{e}^{tA}\Vert}{\Vert \mathrm{e}^{tA}\widehat{y}_0\Vert }$.
The second part regarding $\sim$ follows since $\epsilon(t)\rightarrow 0$ and $%
\epsilon(t,\widehat{y}_0)\rightarrow 0$, $t\rightarrow +\infty$.
\end{proof}

We define the function
$$
t\rightarrow K_\infty(t,y_0),\ t\in\mathbb{R},
$$
as the \emph{asymptotic pointwise condition number} of the problem (\ref{due}).

\begin{remark}
\label{KtAy} \quad

\begin{itemize}
\item[1.] The asymptotic pointwise condition number  $K_\infty(t,y_0)$:
\begin{itemize}
\item [a)] is bounded and away from zero as $t$ varies if $j^\ast=1$ and $L^\ast=L_1$: this follows by
$$
\sup\limits_{t\in \mathbb{R}} \left\Vert Q_{1L_1} (t)\right\Vert<+\infty\text{\ \ and\ \ }\inf\limits_{t\in \mathbb{R}} \left\Vert Q_{j^{\ast}L^{\ast}} (t)\widehat{y}_0\right\Vert>0
$$
and
$$
\inf\limits_{t\in \mathbb{R}} \left\Vert Q_{1L_1} (t)\right\Vert>0\text{\ \ and\ \ }\sup\limits_{t\in \mathbb{R}} \left\Vert Q_{j^{\ast}L^{\ast}} (t)\widehat{y}_0\right\Vert<+\infty
$$

(recall 1), 2), 3) and 4) in Proposition \ref{lemmaQjl} of Appendix \ref{lastsection});

\item [b)] diverges polynomially to infinity, as $t\rightarrow +\infty$, if $j^\ast=1$ and $L^\ast<L_1$;

\item [c)] diverges exponentially to infinity, as $t\rightarrow +\infty$, if $j^\ast>1$. 
\end{itemize}
Therefore, $K_\infty(t,y_0)$ does not diverge to infinity if and only if $y_0$ is as dominant as possible, meaning that $y_0$ has the smallest possible primary index of dominance (i.e., $j^\ast = 1$) and simultaneously the largest possible secondary index of dominance (i.e., $L^\ast = L_1$).

\item[2.] The case $j^{\ast }=1$ and $L^{\ast }=L_{1}$ is generic for $y_0$. In this generic case, we have
\begin{equation*}
K_\infty\left( t,y_{0}\right) = \frac{\left\Vert Q_{1L_{1}}(t)\right\Vert 
}{\left\Vert Q_{1L_{1}}(t)\widehat{y}_{0}\right\Vert }.
\end{equation*}

\end{itemize}
\end{remark}

\subsection{Is $K_\infty\left( t,y_{0}\right)$ the worst $%
K_\infty\left( t,y_{0},\widehat{z}_0\right)$?}

\label{wcab} By definition, $K\left(
t,y_{0}\right)$ is the worst $K\left( t,y_{0},%
\widehat{z}_0\right)$, by varying $\widehat{z}_0$. Hence, an interesting question is the following. Does this fact hold asymptotically as $t\rightarrow +\infty$? Specifically, is $K_\infty\left(t,y_{0}\right)$ the worst $K_\infty\left( t,y_{0},%
\widehat{z}_0\right)$, by varying $\widehat{z}_0$? In other words, do 
\begin{equation*}
\max\limits_{\substack{ \widehat{z}_0\in \mathbb{C}^{n}  \\ \Vert \widehat{z}%
_0\Vert =1}}\ \text{(worst case) and }\ t\rightarrow +\infty\ \text{%
(asymptotic behavior)}
\end{equation*}
commute? The answer is YES and it is given by the next theorem, which considers the ratio
\begin{equation}
\frac{K_\infty\left( t,y_{0},\widehat{z}_{0}\right)}{K_\infty\left(
t,y_{0}\right)}=\frac{L_1!}{L^{\ast\ast }!}\mathrm{e}^{\left( r_{j^{\ast
\ast }}-r_{1}\right) t}t^{L^{\ast \ast }-L_1}\frac{\left\Vert
Q_{j^{\ast\ast}L^{\ast\ast}}(t)\widehat{z}_{0}\right\Vert }{\left\Vert
Q_{1L_1}(t)\right\Vert },\label{rat}
\end{equation}
which is independent of $y_0$.
\begin{theorem} \label{wcy}
We have
\begin{equation}
\max\limits_{\substack{ \widehat{z}_0\in \mathbb{C}^{n}  \\ \Vert \widehat{z}%
_0\Vert =1}}\limsup\limits_{t\rightarrow +\infty}\frac{K_\infty\left(
t,y_{0},\widehat{z}_{0}\right)}{K_\infty\left( t,y_{0}\right)}=1.  \label{limsup}
\end{equation}
In particular:
\begin{itemize}
\item [a)] For any direction of perturbation $\widehat{z}_0$ such that $j^{\ast\ast}>1$ or $j^{\ast\ast}=1$ and $L^{\ast\ast}<L_1$, we have
$$
\lim\limits_{t\rightarrow +\infty}\frac{K_\infty\left(
t,y_{0},\widehat{z}_{0}\right)}{K_\infty\left( t,y_{0}\right)}=0.
$$
\item [b)] For any direction of perturbation $\widehat{z}_0$ such that $j^{\ast\ast}=1$ and $L^{\ast\ast}=L_1$, we have
$$ 
\frac{K_\infty\left(
t,y_{0},\widehat{z}_{0}\right)}{K_\infty\left( t,y_{0}\right)}\leq 1.
$$
\item [c)] There exists a direction of perturbation $\widehat{z}_0$, independent of $y_0$ and with $j^{\ast\ast}=1$ and $L^{\ast\ast}=L_1$, and a sequence $\{t_m\}$ with $t_m\rightarrow +\infty$, as $m\rightarrow \infty$, such that
$$
\lim\limits_{m\rightarrow \infty}\frac{K_\infty\left(
t_m,y_{0},\widehat{z}_{0}\right)}{K_\infty\left( t_m,y_{0}\right)}=1.
$$ 
\end{itemize}
\end{theorem}
\begin{proof}
Points a), b) and c) imply (\ref{limsup}).

Points a) and b) immediately follow by (\ref{rat}): for a) observe that
$$
\sup\limits_{t\in\mathbb{R}}\left\Vert Q_{j^{\ast\ast} L^{\ast\ast}}(t)\widehat{z}_0\right\Vert<+\infty\text{\ \ and\ \ }\inf\limits_{t\in\mathbb{R}}\left\Vert Q_{1L_1}(t)\right\Vert>0
$$
(see 2) and 3) in Proposition \ref{lemmaQjl} of Appendix \ref{lastsection}).

Now, we prove c). Consider a sequence $\{t_k\}$ such that $t_k\rightarrow +\infty$, $%
k\rightarrow \infty$, and a sequence $\{\widehat{z}_{0k}\}$ such that $\widehat{z}_{0k}\in\mathbb{C}^n$, $\left\Vert \widehat{z}%
_{0k}\right\Vert=1$ and
\begin{equation*}
\left\Vert Q_{1L_1}(t_k)\right\Vert=\left\Vert Q_{1L_1}(t_k)\widehat{z}%
_{0k}\right\Vert. 
\end{equation*}
By the compactness of
the unit sphere in $\mathbb{C}^n$, there exists a subsequence $\{\widehat{z}%
_{0k_m}\}$ of $\{\widehat{z}_{0k}\}$ converging to some $\widehat{z}%
_{0\infty}\in\mathbb{C}^n$ with $\left\Vert \widehat{z}_{0\infty}\right%
\Vert=1$.

We have $j\left(\widehat{z}_{0\infty}\right)=1$ and $L_1\left(\widehat{z}%
_{0\infty}\right)=L_1$. 

In fact, for any index $m$, we have
$$
\left\vert\ \left\Vert Q_{1L_1}(t_{k_m})\widehat{z}%
_{0\infty}\right\Vert-\left\Vert Q_{1L_1}(t_{k_m})\right\Vert\ \right\vert \leq \sup\limits_{t\in \mathbb{R}} \left\Vert Q_{1L_1} (t)\right\Vert\left\Vert \widehat{z}_{0\infty}-\widehat{z}_{0k_m}\right\Vert,
$$
where the right-hand side goes to zero as $m\rightarrow \infty$
(recall point 1) in Proposition \ref{lemmaQjl} of Appendix \ref{lastsection}). Therefore, there exists an index $m$ such that
$$
\left\Vert Q_{1L_1}(t_{k_m})\widehat{z}%
_{0\infty}\right\Vert\geq \frac{1}{2}\inf\limits_{t\in \mathbb{R}} \left\Vert Q_{1L_1} (t) \right\Vert,
$$ 
where the right-hand side is positive (recall point 3) in Proposition \ref{lemmaQjl} of Appendix \ref{lastsection}). Since
$
Q_{1L_1}(t_{k_m})$\ $\widehat{z}_{0\infty}\neq 0,
$
we cannot have  $j\left(\widehat{z}_{0\infty}\right)>1$, otherwise $Q_{1L_1}(t_{k_m})\widehat{z}_{0\infty}=0$ (recall point 1) in Proposition \ref{lemma2Qjl} of Appendix B). Hence, $j\left(\widehat{z}_{0\infty}\right)=1$. Moreover, since
$
Q_{1L_1}(t_{k_m})\widehat{z}_{0\infty}\neq 0,
$
we cannot have  $L_1\left(\widehat{z}_{0\infty}\right)<L_1$, otherwise $Q_{1L_1}(t_{k_m})\widehat{z}_{0\infty}=0$ (remind point 2) in Proposition \ref{lemma2Qjl} of Appendix B). Hence, $L_1\left(\widehat{z}%
_{0\infty}\right)=L_1$.

By using as a direction of perturbation $\widehat{z}_{0\infty}$, we have $%
j^{\ast\ast}=1$ and $L^{\ast\ast}=L_1$. Thus, for any index $m$, 
\begin{eqnarray*}
\left\vert \frac{K_\infty\left( t_{k_m},y_{0},\widehat{z}_{0\infty}\right)%
}{K_\infty\left( t_{k_m},y_{0}\right)}-1\right\vert&=&\left\vert \frac{%
\left\Vert Q_{1L_1}(t_{k_m})\widehat{z}_{0\infty}\right\Vert }{\left\Vert
Q_{1L_1}(t_{k_m})\right\Vert }-1\right\vert \\
&\leq& \left\Vert \widehat{z}_{0\infty}-\widehat{z}_{0k_m}\right\Vert.
\end{eqnarray*}
We conclude that 
\begin{equation*}
\lim\limits_{m\rightarrow \infty}\frac{K_\infty\left( t_{k_m},y_{0},%
\widehat{z}_{0\infty}\right)}{K_\infty\left( t_{k_m},y_{0}\right)}=1.
\end{equation*}
\end{proof}

\begin{remark}\label{realeigen1}
	
\quad
	
	\begin{itemize}
		\item [1.] When $\Lambda_1$ consists of a real eigenvalue, the point c) is modified to:
\begin{itemize}
	\item [$c^{\prime}$)] There exists a direction of perturbation $\widehat{z}_0$, independent of $y_0$ and with $j^{\ast\ast}=1$ and $L^{\ast\ast}=L_1$, such that
\begin{equation}
K_\infty\left(t,y_{0},\widehat{z}_{0}\right)=K_\infty\left(t,y_{0}\right). \label{K=K}
\end{equation}
\end{itemize}
In fact, in this case $Q_{1L_1}(t)=Q_{1L_1}$ is independent of $t$  (see (\ref{Qjl})); hence, there exists $\widehat{z}_{0}\in\mathbb{C}^n$, $\left\Vert \widehat{z}%
_{0}\right\Vert=1$, independent of $t$ such that 
\begin{equation*}
\left\Vert Q_{1L_1}\right\Vert=\left\Vert Q_{1L_1}\widehat{z}%
_{0}\right\Vert. 
\end{equation*}
By using as a direction of perturbation $\widehat{z}_{0}$, we have $%
j^{\ast\ast}=1$ and $L^{\ast\ast}=L_1$ and then  we obtain (\ref{K=K}) by (\ref{rat}).

\item [2.] Observe that
$$
		\max\limits_{\substack{ \widehat{z}_0\in \mathbb{C}^{n}  \\ \Vert \widehat{z}%
				_0\Vert =1}}\limsup\limits_{t\rightarrow +\infty}\frac{K\left(
			t,y_{0},\widehat{z}_{0}\right)}{K\left( t,y_{0}\right)}=\max\limits_{\substack{ \widehat{z}_0\in \mathbb{C}^{n}  \\ \Vert \widehat{z}%
			_0\Vert =1}}\limsup\limits_{t\rightarrow +\infty}\frac{K_\infty\left(
		t,y_{0},\widehat{z}_{0}\right)}{K_\infty\left( t,y_{0}\right)}
$$
and then
$$
\max\limits_{\substack{ \widehat{z}_0\in \mathbb{C}^{n}  \\ \Vert \widehat{z}%
		_0\Vert =1}}\limsup\limits_{t\rightarrow +\infty}\frac{K\left(
	t,y_{0},\widehat{z}_{0}\right)}{K\left( t,y_{0}\right)}=1,
$$
which, more clearly than (\ref{limsup}), shows that
$$
\max\limits_{\substack{ \widehat{z}_0\in \mathbb{C}^{n}  \\ \Vert \widehat{z}%
		_0\Vert =1}}\ \text{(worst case) and }\limsup\limits_{t\rightarrow +\infty}\text{(asymptotic behavior)}
$$ 
commute.
\end{itemize}
\end{remark}

In conclusion, the \emph{asymptotic pointwise} condition number $K_\infty(t.y_0)$, i.e., the asymptotic form of the worst $K(t,y_0,\widehat{z}_0)$ by varying  $\widehat{z_0}$, coincides with   the \emph{pointwise asymptotic} condition number, i.e., the worst asymptotic form of $K(t,y_0,\widehat{z}_0)$ by varying  $\widehat{z_0}$. 

\section{The asymptotic behavior of $K(t)$} \label{asymptoticbehavior2}

In this section, we study the asymptotic behavior of the global condition number $K(t)$. Its asymptotic form is the asymptotic condition number $K^+_\infty(t)$. We also show that $K^+_\infty(t)$ does not coincide with the worst $K_\infty(t,y_0)$, by varying $y_0$, i.e., in light of Theorem \ref{wcy}, it does not coincide with the worst $K_\infty(t,y_0,\widehat{z}_0)$, by varying  $y_0$ and $\widehat{z}_0$. The worst $K_\infty(t,y_0)$ is the asymptotic condition number $K_\infty(t)$. 

\subsection{The asymptotic condition number $K_\infty^{+}\left(t\right) $} \label{SKtA}

Next theorem describes the asymptotic form of $K\left( t\right) $, the worst $K\left( t,y_0\right) $ by varying $y_0$.

\begin{theorem}
\label{ThmK} We have 
\begin{equation*}
K\left( t\right) \approx K_\infty^{+}(t),
\end{equation*}
where 
\begin{equation}  \label{sim3}
K_\infty^{+}\left( t\right) := \frac{1}{L_1!L_{q}!}\mathrm{e}^{\left(
r_{1}-r_{q}\right) t}t^{L_{1}+L_q} \left\Vert Q_{1L_{1}}(t) \right\Vert\
\cdot \left\Vert Q_{qL_{q}}(-t)\right\Vert
\end{equation}
with precision 
\begin{equation*}
\epsilon(t)+\epsilon(t,-A)+\epsilon(t)\epsilon(t,-A),
\end{equation*}
where $\epsilon(t,-A)$ is $\epsilon(t)$ (defined in Proposition \ref{Lemma2}) for the matrix $-A$. Moreover, we have 
\begin{equation*}
K\left( t\right) \sim K^{+}_\infty(t),\ t\rightarrow +\infty.  \label{simK}
\end{equation*}
\end{theorem}

\begin{proof}
Proposition \ref{Lemma2} states the asymptotic forms of $\mathrm{e}^{tA}$ and $\mathrm{e}^{-tA}$ in (\ref{KAy0}) and how they are dominant at a finite time $t$: we have
\begin{equation*}
\mathrm{e}^{tA}\approx \mathrm{e}^{r_1t}\frac{t^{L_1}}{L_1!} Q_{1L_1}(t)
\end{equation*}
with precision $\epsilon(t)$ and 
\begin{equation*}
\mathrm{e}^{-tA} =\mathrm{e}^{t(-A)}\approx \mathrm{e}^{r_1(-A)t}\frac{t^{L_1(-A)}}{L_1(-A)!} Q_{1L_1(-A)}(t,-A)=\mathrm{e}^{-r_qt}\frac{t^{L_{q}}}{L_{q}!}(-1)^{L_q}
Q_{qL_{q}}(-t)
\end{equation*}
with precision $\epsilon(t,-A)$. In the latter, $r_1(-A)$, $L_1(-A)$ and $Q_{1L_1(-A)}(t,-A)$ are $r_1$, $L_1$ and $Q_{1L_1}(t)$ for the matrix $-A$, and we have
$$
r_1(-A)=-r_q,\ L_1(-A)=L_q\text{\ \ and\ \ }Q_{1L_1(-A)}(t,-A)=(-1)^{L_q}Q_{qL_q}(-t)
$$
(recall (\ref{q+1-}) and Proposition \ref{-A} in Appendix \ref{lastsection}). The first part of the theorem regarding $\approx$ follows by using
Remark \ref{Remark fg} and by bounding the relative error of the product $\Vert \mathrm{e}^{tA}\Vert \Vert \mathrm{e}^{t(-A)}\Vert$ in terms of the bounds  $\epsilon(t)$ and  $\epsilon(t,-A)$ of the relative errors of $\Vert \mathrm{e}^{tA}\Vert$ and $\Vert \mathrm{e}^{t(-A)}\Vert $, respectively. The second part regarding $\sim$ follows since $%
\epsilon(t)\rightarrow 0$ and $\epsilon(t,-A)\rightarrow 0$, $%
t\rightarrow +\infty$.
\end{proof}

We  define the function
$$
t\rightarrow K^{+}_\infty(t),\ t\in\mathbb{R},
$$
as the \emph{asymptotic global condition number} of the problem (\ref{due}).

\subsection{Is $K^+_\infty\left( t\right)$ the worst $K_\infty\left(
t,y_{0}\right)$?}

\label{sss}

By definition, $K\left(
t\right)$ is the worst $K\left( t,y_{0}\right)$, by varying $y_0$. Hence, as in Subsection \ref{wcab}, an interesting question is the following. Does this fact hold asymptotically as $t\rightarrow +\infty$? Specifically, is  $%
K^+_\infty\left( t\right)$ the worst $K_\infty\left(
t,y_{0}\right)$, by varying $y_0$? In other words, do 
\begin{equation*}
\max\limits_{\substack{ y_0\in \mathbb{C}^{n}  \\ y_0\neq 0}}\ \text{(worst case) and }\ t\rightarrow +\infty\ \text{(asymptotic
behavior)}
\end{equation*}
commute? Unlike the similar question in Subsection \ref{wcab}, here the answer is NO and it is given by the next theorem, which considers
\begin{equation}
K_\infty\left( t\right) := \frac{1}{L_{1}!}\mathrm{e}^{\left(
r_{1}-r_{q}\right) t}t^{L_{1}}\left\Vert
Q_{1L_{1}}(t)\right\Vert\cdot\left\Vert Q_{q0}(-t)|_{U^e_{q}}\right\Vert  \label{Kover}
\end{equation}
and the ratio
\begin{equation}
\frac{K_\infty\left( t,y_{0}\right)}{K_\infty\left( t\right)}%
=L^{\ast}!\mathrm{e}^{\left( r_{q}-r_{j^\ast}\right) t}t^{-L^{\ast}}\frac{1}{%
\left\Vert Q_{j^{\ast}L^{\ast}}(t)\widehat{y}_{0}\right\Vert\left\Vert
Q_{q0}(-t)|_{U^e_{q}}\right\Vert},  \label{rat2}
\end{equation}
where (see Appendix \ref{sssQj0}) $U^e_q$ is the subspace of $\mathbb{C}^n$ sum of the eigenspaces corresponding to eigenvalues in $\Lambda_q$ and
$Q_{q0}(-t)|_{U^e_{q}}$ is the restriction of the linear operator $Q_{q0}(-t)$ to the subspace $U^e_{q}$. 

Observe that, for $u\in\mathbb{C}^n$, we have $u\in U^e_q$ if and only if $j(u)=q$ and $L(u)=0$. For the latter, recall that $L(u)=0$ means $l_i(u)=0$ for any $\lambda_i\in\Lambda_q(u)$ (see  Table \ref{table1}). Therefore,
\begin{equation}
\text{$y_0\in U^e_q$ if and only if $j^\ast=q$ and $L^\ast=0$.}
\label{q0}
\end{equation}

\begin{theorem}
\label{asworst2} We have
\begin{equation}
\max\limits_{\substack{ y_0\in \mathbb{C}^{n}  \\ y_0\neq 0}}%
\limsup\limits_{t\rightarrow +\infty}\frac{K_\infty\left( t,y_{0}\right)}{%
K_\infty\left( t\right)}=1.   \label{limsup2}
\end{equation}
In particular:
\begin{itemize}
\item [a)] For any initial value $y_0$ such that $j^{\ast}<q$ or $j^{\ast}=q$ and $L^{\ast}>0$, i.e., $y_0\notin U^e_q$ (recall (\ref{q0})), we have
$$
\lim\limits_{t\rightarrow +\infty}\frac{K_\infty\left(
t,y_{0}\right)}{K_\infty\left( t\right)}=0.
$$
\item [b)] For any initial value $y_0$ such that $j^{\ast}=q$ and $L^{\ast}=0$, i.e., $y_0\in U^e_q$, we have
$$ 
\frac{K_\infty\left(
t,y_{0}\right)}{K_\infty\left( t\right)}\leq 1.
$$
\item [c)] There exists an initial value $y_0\in U^e_q$ and a sequence $\{t_m\}$ with $t_m\rightarrow +\infty$, as $m\rightarrow \infty$, such that
$$
\lim\limits_{m\rightarrow \infty}\frac{K_\infty\left(
t_m,y_{0}\right)}{K_\infty\left( t_m\right)}=1.
$$
\item [d)] For the initial value $y_0$ at point c), there exists a direction of perturbation $\widehat{z}_0$ with $j^{\ast\ast}=1$ and $L^{\ast\ast}=L_1$, and a subsequence $\{t_{m_s}\}$ of the sequence $\{t_{m}\}$ at point c), such that
$$
\lim\limits_{s\rightarrow \infty}\frac{K_\infty\left(
t_{m_s},y_{0},\widehat{z}_{0}\right)}{K_\infty\left( t_{m_s}\right)}=1.
$$
\end{itemize} 
\end{theorem} 
\begin{proof}
Points a), b) and c) imply (\ref{limsup2}).

Point a) follows by (\ref{rat2}): observe that
$$
\inf\limits_{t\in\mathbb{R}}\left\Vert Q_{j^{\ast}L^{\ast}}(t)\widehat{y}_{0}\right\Vert>0\text{\ \ and\ \ }\inf\limits_{t\in\mathbb{R}}\left\Vert
Q_{q0}(-t)|_{U^e_{q}}\right\Vert>0
$$
(recall point 4) in Proposition \ref{lemmaQjl} of Appendix \ref{lastsection} and Remark \ref{remQ0} of Appendix \ref{sssQj0}).

Point b) follows by (\ref{rat2}): we have
\begin{equation}
	\min\limits_{\substack{\widehat{u}\in U^e_{q} \\ \left\Vert \widehat{u}\right\Vert=1}}\left\Vert Q_{q0}(t)%
	\widehat{u}\right\Vert=\frac{1}{\left\Vert Q_{q0}(-t)|_{U^e_{q}} \right\Vert}. \label{min=1/}
\end{equation}
(recall Remark \ref{remQ00} of Appendix \ref{sssQj0}).

For the point c), we proceed as in the proof of Theorem \ref{wcy}. Consider a
sequence $\{t_k\}$ such that $t_k\rightarrow +\infty$, $k\rightarrow
\infty$, and a sequence $\{\widehat{y}_{0k}\}$, where $\widehat{y}_{0k}\in U^e_{q}$ and $\left\Vert \widehat{y}_{0k}\right\Vert=1$, such that
\begin{equation}
\left\Vert Q_{q0}(t_k)\widehat{y}_{0k}\right\Vert=\min\limits_{\substack{\widehat{u}\in U^e_{q} \\ \left\Vert \widehat{u}\right\Vert=1}}\left\Vert Q_{q0}(t_k)%
\widehat{u}\right\Vert=\frac{1}{\left\Vert Q_{q0}(-t_k)|_{U^e_{q}}
\right\Vert}  \label{K/K-1}
\end{equation}
(recall (\ref{min=1/})).  There exists a
subsequence $\{\widehat{y}_{0k_m}\}$ of $\{\widehat{y}_{0k}\}$ converging to some $%
\widehat{y}_{0\infty}\in U^e_q$ with $\left\Vert \widehat{y}_{0\infty}\right\Vert=1$. By using $%
\widehat{y}_{0\infty}$ as initial value, we have
\begin{eqnarray*}
\left\vert \frac{K_\infty\left( t_{k_m},\widehat{y}_{0\infty}\right)}{K%
_\infty\left( t_{k_m}\right)}-1\right\vert&=&\left\vert \frac{1}{%
\left\Vert Q_{q0}(t_{k_m})\widehat{y}_{0\infty}\right\Vert\left\Vert
Q_{q0}(-t_{k_m})|_{U^e_{q}}\right\Vert}-1\right\vert \\
&=&\left\vert\frac{\Vert Q_{q0}(t_{k_m})\widehat{y}_{0k_m}\Vert}{\Vert Q_{q0}(t_{k_m})\widehat{y}_{0\infty}\Vert}-1\right\vert  \ \ \ \text{(by (\ref{K/K-1}))}\\
&\leq& \frac{\left\Vert Q_{q0}(t_{k_m})\right\Vert}{\left\Vert Q_{q0}(t_{k_m})\widehat{y}_{0\infty}\right\Vert}\left\Vert \widehat{y}_{0k_m}-\widehat{y}_{0\infty}\right\Vert
\rightarrow 0,\ m\rightarrow \infty\\
&&\text{($\rightarrow 0$ by $\sup\limits_{t\in \mathbb{R}} \left\Vert Q_{q0} (t)\right\Vert<+\infty\text{\ \ and\ \ }\inf\limits_{t\in \mathbb{R}} \left\Vert Q_{q0} (t)\widehat{y}_{0\infty}\right\Vert>0$:}\\
&&\text{ see 2) and 4) in Proposition \ref{lemmaQjl} of Appendix \ref{lastsection})}.
\end{eqnarray*}
We conclude that 
\begin{equation*}
\lim\limits_{m\rightarrow \infty}\frac{K_\infty\left(
t_{k_m},\widehat{y}_{0\infty}\right)}{K_\infty\left( t_{k_m}\right)}=1.
\end{equation*}

For the point d), repeat the proof of point c) in Theorem \ref{wcy} with the sequence $\{t_k\}$  replaced by the sequence $\{t_m\}$ at point c) of this theorem. In this way, we show that there exists a direction of perturbation $\widehat{z}_0$ with $j^{\ast\ast}=1$ and $L^{\ast\ast}=L_1$, and a subsequence $\{t_{m_s}\}$ of the sequence $\{t_{m}\}$, such that
$$
\lim\limits_{s\rightarrow \infty}\frac{K_\infty\left(
	t_{m_s},\widehat{y}_{0\infty},\widehat{z}_{0}\right)}{K_\infty\left( t_{m_s},\widehat{y}_{0\infty}\right)}=1
$$
and then 
$$
\lim\limits_{s\rightarrow \infty}\frac{K_\infty\left(
	t_{m_s},\widehat{y}_{0\infty},\widehat{z}_{0}\right)}{K_\infty\left( t_{m_s}\right)}=\lim\limits_{s\rightarrow \infty}\frac{K_\infty\left(
	t_{m_s},\widehat{y}_{0\infty}\right)}{K_\infty\left( t_{m_s}\right)}\cdot \lim\limits_{s\rightarrow \infty}\frac{K_\infty\left(
	t_{m_s},\widehat{y}_{0\infty},\widehat{z}_{0}\right)}{K_\infty\left( t_{m_s},\widehat{y}_{0\infty}\right)}=1.
$$

\end{proof}

\begin{remark}\label{realeigen2}
When $\Lambda_q$ consists of a real eigenvalue, the point c) is modified to:
\begin{itemize}
	\item [$c^{\prime}$)] For any initial value $y_0\in U^e_{q}$, we have
\begin{equation*}
K_\infty\left(t,y_{0}\right)=K_\infty\left(t\right). \label{KK=KK}
\end{equation*}
\end{itemize}
This follows by (\ref{rat2}) and  the fact that
$$
Q_{q0}(-t)|_{U^e_q}=I_{U^e_q}
$$
(see Remark \ref{star} of Appendix \ref{sssQj0}).

In addition, if $\Lambda_1$ also consists of a real eigenvalue, the point d) is modified to:
\begin{itemize}
	\item [$d^{\prime}$)]There exists a direction of perturbation $\widehat{z}_0$ with $j^{\ast\ast}=1$ and $L^{\ast\ast}=L_1$ such that, for any initial value $y_0\in U^e_{q}$, we have
\begin{equation*}
K_\infty\left(t,y_{0},\widehat{z}_0\right)=K_\infty\left(t,y_{0}\right)=K_\infty\left(t\right).
\end{equation*}
\end{itemize}
See point 1 in Remark \ref{realeigen1}.

\end{remark}

\subsection{The asymptotic condition number $K_\infty\left(t\right) $} \label{ssss}

\label{SoverKtA}

Theorem \ref{asworst2} says that $K_\infty\left( t\right)$ in (\ref{Kover}),
not $K^+_\infty\left( t\right)$ in (\ref{sim3}), is the worst $K_\infty\left( t,y_0\right)$, by varying $y_0$,
i.e., the worst $K_\infty\left( t, y_0,%
\widehat{z}_0\right)$, by varying both $y_0$ and $\widehat{z}_0$.

We define the function
$$
t\rightarrow K_\infty(t),\ t\in\mathbb{R},
$$
as the \emph{global asymptotic condition number} of the problem (\ref{due}).

Therefore, the \emph{asymptotic global} condition number $K_\infty^+(t)$, i.e. the asymptotic form of the worst $K(t,y_0)$ by varying  $y_0$, does not coincide with   the \emph{global asymptotic} condition number $K_\infty(t)$, i.e. the worst asymptotic form of $K(t,y_0)$ by varying  $y_0$.

Indeed, $K^+_\infty(t)$ can be significantly larger than $K%
 _\infty(t)$: we have 
 \begin{equation*}
 	\frac{K_\infty\left( t\right)}{K^+_\infty\left( t\right)}=L_q!
 	t^{-L_q}\frac{\left\Vert Q_{q0}(-t)|_{U^e_{q}}\right\Vert}{\left\Vert
 		Q_{qL_q}(-t)\right\Vert}
 \end{equation*}
 and then 
 \begin{equation*}
 	\lim\limits_{t\rightarrow +\infty}\frac{K_\infty\left( t\right)%
 	}{K^+_\infty\left( t\right)}=0
 \end{equation*}
 for $L_q>0$ and
 \begin{equation}
 \frac{K_\infty\left( t\right)}{K^+_\infty\left(
 	t\right)}=\frac{\left\Vert Q_{q0}(-t)|_{U^e_{q}}\right\Vert}{\left\Vert
 	Q_{q0}(-t)\right\Vert}\leq 1  \label{1/+}
 \end{equation}
 for $L_q=0$.  

\begin{remark} \label{KtA}
Both the asymptotic global condition number $K^+_\infty(t)$ and the global asymptotic condition number $K_\infty(t)$:
	
	\begin{itemize}
		\item 
		are bounded and away from zero as $t$ varies if $q=1$ and $L_1=0$: this follows by
		$$
		\sup\limits_{t\in\mathbb{R}}\left\Vert Q_{1L_1}(t)\right\Vert,\sup\limits_{t\in\mathbb{R}}\left\Vert Q_{qL_q}(t)\right\Vert,\sup\limits_{t\in\mathbb{R}}\left\Vert Q_{q0}|_{U^e_q}(-t)\right\Vert<+\infty 
		$$
		and
$$
\inf\limits_{t\in\mathbb{R}}\left\Vert Q_{1L_1}(t)\right\Vert,\inf\limits_{t\in\mathbb{R}}\left\Vert Q_{qL_q}(t)\right\Vert,\inf\limits_{t\in\mathbb{R}}\left\Vert Q_{q0}|_{U^e_q}(-t)\right\Vert>0 
$$		
(see Proposition \ref{lemmaQjl} of Appendix \ref{lastsection} and Remark \ref{remQ0} of Appendix \ref{sssQj0})	
		\item  diverge polynomially to infinity, as $t\rightarrow +\infty$, if $q=1$ and $L_1>0$;
		
		\item diverge exponentially to infinity, as $t\rightarrow +\infty$, if $q>1$.
	\end{itemize}
	Therefore, $K_\infty^+(t)$ and $K_\infty(t)$ do not diverge to infinity if and only if all the eigenvalues of $A$ have the same real part, i.e., they lie in a vertical line of the complex plane, and are non-defective. Observe that in this case $K(t)=K^+_\infty(t)$ (we have $\epsilon(t)=\epsilon(t,-A)=0$ in Theorem \ref{ThmK})  and $K^+_\infty(t)=K_\infty(t)$ (we have $U^e_q=\mathbb{C}^n$ in (\ref{1/+})). 
	
\end{remark}

In the next example, we illustrate the difference between $K^+_\infty(t)$ and $K_\infty(t)$.
\begin{example}
 Consider the non-diagonalizable matrix
$$
A_1=V
\left[
\begin{array}{rr}
1 & 1\\
0 & 1
\end{array}
\right]V^{-1}
$$
and the diagonalizable matrix
$$
A_2=V
\left[
\begin{array}{rr}
1 & 0\\
0 & -1
\end{array}
\right]V^{-1},
$$
where
$$
V=\left[
\begin{array}{rr}
1 & \frac{1}{2}\\
1 &  1
\end{array}
\right]=\left[v^{(1)}\ v^{(2)}\right]\text{\ \ and\ \ }
V^{-1}=2\left[
\begin{array}{rr}
1 &  -\frac{1}{2}\\
-1 & 1
\end{array}
\right]=\left[\begin{array} {c} w^{(1)}\\w^{(2)}\end{array}\right]
$$
with $v^{(1)}$ and $v^{(2)}$ the columns of $V$, and $w^{(1)}$ and $w^{(2)}$ the rows of $V^{-1}$.

\underline{The matrix $A_1$} \newline

For the matrix $A_1$, we have
$$
\Lambda=\{1\},\ q=1{\ \ and\ \ }L_1=L_q=1.
$$
The asymptotic global and global asymptotic condition numbers are
\begin{equation*}
K^+_\infty(t)=t^2\left\Vert Q_{11}(t)\right\Vert \left\Vert Q_{11}(-t)\right\Vert\text{\ and\ \ }K_\infty(t)=t\left\Vert Q_{11}(t)\right\Vert \left\Vert Q_{10}(-t)|_{U^e_{1}}\right\Vert,
\end{equation*}
where 
$$
Q_{11}(t)=Q_{11}(-t)=P_{11}=v^{(1)}w^{(2)}
$$
(recall (\ref{Qjl}) and, in Appendix \ref{JCFsection}, point 2 in Remark \ref{mi}) and 
$$
Q_{10}(-t)|_{U^e_{1}}=I_{U^e_{1}}\\
$$
(recall Remark \ref{star} of Appendix \ref{sssQj0}). Thus, for the Euclidean norm as vector norm, we have, since $\Vert P_{11}\Vert_2=\Vert v^{(1)}\Vert_2\Vert w^{(2)}\Vert_2=4$, 
\begin{equation*}
K^+_\infty(t)=16t^2\text{\ and\ \ }K_\infty(t)=4t.
\end{equation*}

In Figure \ref{KKoverline}, we see $K(t)$ (blue dashed line), $K^{+}_\infty(t)$ (red dashed line) and $K_\infty(t)$ (red dashed line), $t\in[0,100]$, in logarithmic scale. Only at  the beginning, $K(t)$ is  distinguishable from its asymptotic form $K^+_\infty(t)$.

\begin{figure}[tbp]
	\includegraphics[width=1\textwidth]{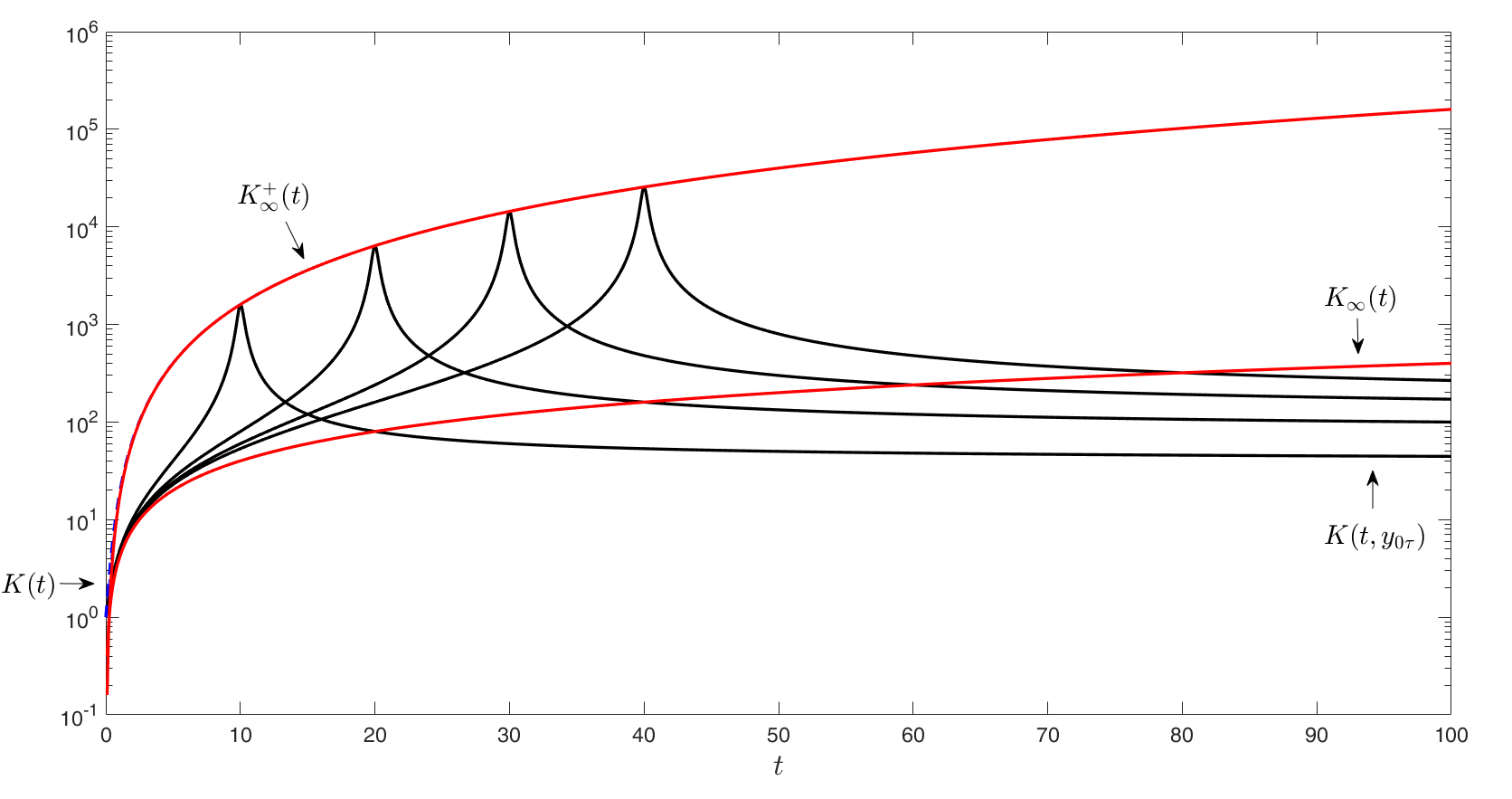}
	\caption{The condition numbers $K(t)$, $K^+_\infty(t)$, $K_\infty(t)$ and $K(\tau,y_{0\tau})$ for the matrix $A_1$. The vector norm is the Euclidean norm.}
	\label{KKoverline} 
\end{figure}

Figure \ref{KKoverline} confirms that the asymptotic global condition number $K^+_{\infty}(t)$ is not the worst asymptotic form of $K(t,y_0)$,  by varying $y_0$. In fact, in the figure we also see, for each $\tau\in\{10,20,30,40\}$,  $K(t,y_{0\tau})$ (black solid lines) for an initial value $y_{0\tau}$ such that $K(\tau,y_{0\tau})=K(\tau)$. To obtain this, we take the initial value $y_{0\tau}$ such that
$$
\Vert \mathrm{e}^{-\tau A} \Vert_2=\frac{1}{\Vert \mathrm{e}^{\tau A} \widehat{y}_{0\tau}\Vert_2};
$$
this is obtained with $y_{0\tau}=\mathrm{e}^{-\tau A}x_{0\tau}$, where $x_{0\tau}$ is such that
$$
\Vert \mathrm{e}^{-\tau A} \Vert_2=\frac{\Vert \mathrm{e}^{-\tau A}x_{0\tau}\Vert_2}{\Vert x_{0\tau}\Vert_2}.
$$

Although $K(t,y_{0\tau})$ touches $K^+_\infty(t)$ at $t=\tau$, i.e., $K(t,y_{0\tau})$ is the worst $K(t,y_0)$ by varying $y_0$ at $t=\tau$, it asymptotically  falls below $K_\infty(t)$, which is much smaller than $K^+_\infty(t)$.

This last fact confirms that $K_\infty(t)$ is the worst $K_\infty(t,y_0)$, by varying  $y_0$. Indeed, we have $K_{\infty}(t)=K_\infty(t,y_0)$ for $y_0\in U^e_1=\mathrm{span}(v^{(1)})$
(see c') in Remark \ref{realeigen2}).\newline 

\underline{The matrix $A_2$}\newline 

For the matrix $A_2$, we have
$$
\Lambda=\{1,-1\},\ q=2\text{\ \ and\ \ }L_1=L_q=0.
$$
The asymptotic global and global asymptotic condition numbers are
\begin{equation*}
K^+_\infty(t)=\mathrm{e}^{2t}\left\Vert Q_{10}(t)\right\Vert \left\Vert Q_{20}(-t)\right\Vert\text{\ and\ \ }K_\infty(t)=\mathrm{e}^{2t}\left\Vert Q_{10}(t)\right\Vert \left\Vert Q_{20}(-t)|_{U^e_2}\right\Vert,
\end{equation*}
where
\begin{equation*}
Q_{10}(t)=P_{10}=v^{(1)}w^{(1)},\ Q_{20}(-t)=P_{20}=v^{(2)}w^{(2)}
\end{equation*}
(see (\ref{Qjl}) and, in Appendix \ref{JCFsection}, Subsection \ref{diagonalizable}) and
\begin{equation*}
Q_{20}(-t)|_{U^e_2}=I_{U^e_2}
\end{equation*}
(see Remark \ref{star} of Appendix \ref{sssQj0}). Thus, for the Euclidean norm as vector norm, we have, since $\Vert P_{10}\Vert_2=\Vert v^{(1)}\Vert_2\Vert w^{(1)}\Vert_2=\sqrt{10}$ and $\Vert P_{20}\Vert_2=\Vert v^{(2)}\Vert_2\Vert w^{(2)}\Vert_2=\sqrt{10}$, 
\begin{equation*}
K^+_\infty(t)=10\mathrm{e}^{2t}\text{\ and\ \ }K_\infty(t)=\sqrt{10}\mathrm{e}^{2t}.
\end{equation*}

In Figure \ref{KKoverline2}, for the matrix $A_2$, we reproduce everything shown in Figure \ref{KKoverline} for the matrix $A_1$. Now, $t\in[0,5]$ and $\tau\in\{1,2,3,4\}$. The same behaviour of $K(t,y_{0\tau})$ appears.

\begin{figure}[tbp]
\includegraphics[width=1\textwidth]{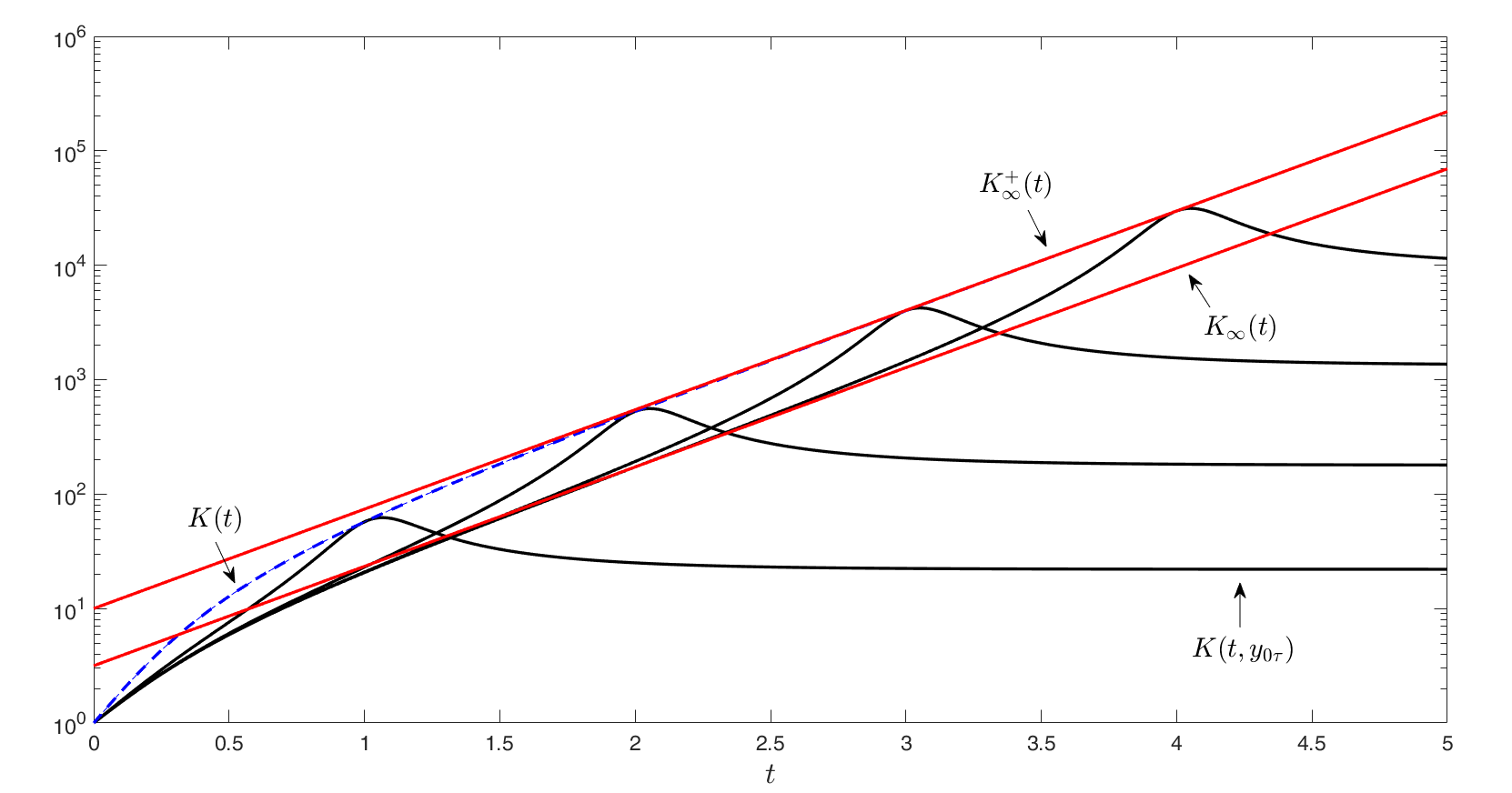}
\caption{The condition numbers $K(t)$, $K^+_\infty(t)$, $K_\infty(t)$ and $K(\tau,y_{0\tau})$ for the matrix $A_2$. The vector norm is the Euclidean norm.}
\label{KKoverline2}
\end{figure}

\end{example}

\section{Final considerations about the asymptotic condition numbers}  \label{sscclast}

Suppose that the matrix $A$ in the ODE (\ref{ode}) does not have all eigenvalues with same real part, meaning $q > 1$ (recall that $q$ is the number of distinct real parts of the eigenvalues of $A$).

In this case, recall Remark \ref{KtA}, the global asymptotic condition number  $K_\infty\left( t\right)$ exponentially diverges. Hence, the relative error $\varepsilon$ of the initial value is exponentially magnified in the relative error $\delta(t)$ of the solution, in the worst case for the initial value and the
perturbation of the initial value.

However, this conclusion is too pessimistic, since this exponential magnification of $\varepsilon$ in $\delta(t)$ appears only in a non-generic case.
In fact, as illustrated in points 1-a) and 2 of Remark \ref{KtAyz}, in the generic case $j^\ast=j^{\ast\ast}=1$ and $L^\ast=L^{\ast\ast}=L_1$ for $y_0$ and $\widehat{z}_0$, we have
\begin{equation}
\frac{\delta(t)}{\varepsilon}=K\left( t,y_{0},\widehat{z}_{0}\right)\sim K_\infty\left( t,y_{0},\widehat{z}_{0}\right) = \frac{\left\Vert
Q_{1L_1}(t)\widehat{z}_{0}\right\Vert }{\left\Vert Q_{1L_1}(t)\widehat{y}%
_{0}\right\Vert },\ t\rightarrow +\infty.  \label{Kgeneric}
\end{equation}
and \emph{$K_\infty\left( t,y_{0},\widehat{z}_{0}\right)$ remains bounded as well as away from zero as $t$ varies}.

We have introduced the three asymptotic condition numbers $K_\infty(t,y_0,\widehat{z}_0)$, $K_\infty(t,$\ $y_0)$ and $K_\infty(t)$. The most important is $K_\infty(t,y_0)$. In fact, $K_\infty(t,y_0,\widehat{z}_0)$ is given in terms of the (in general) unknown direction of perturbation $\widehat{z}_0$ and $K_\infty(t)$ is a worst asymptotic magnification factor $\frac{\delta(t)}{\varepsilon}$ which applies in a non-generic case.

As illustrated in points 1-a) and 2 of Remark \ref{KtAy}, in the generic case $j^\ast=1$ and $L^\ast=L_1$ for $y_0$, this most important asymptotic condition number $K_\infty(t,y_0)$, representing the asymptotic magnification factor $\frac{\delta(t)}{\varepsilon}$ for the worst perturbation of $y_0$, is given by 
\begin{equation}
K_\infty\left( t,y_{0}\right) = \frac{\left\Vert Q_{1L_1}(t)\right\Vert 
}{\left\Vert Q_{1L_1}(t)\widehat{y}_{0}\right\Vert } \label{Kgeneric2}
\end{equation}
and \emph{it remains bounded as well as away from zero as $t$ varies}.

\subsection{The RLGE condition}\label{RLGE}

In this final subsection, we consolidate the conclusions (\ref{Kgeneric}) and (\ref{Kgeneric2}) in a theorem that provides further details. 

Let $u\in\mathbb{C}^n$. We say that $u$ satisfies the \emph{Rightmost Last Generalized Eigenvector (RLGE) condition} if $j(u)=1$ and $L(u)=L_1$, i.e., $u$ is as dominant as possible.

In terms of the JCF of $A$, this means that there is a non-zero component of $u$ along the last generalized eigenvector in some of the longest Jordan chains  of the rightmost eigenvalues of $A$ (recall Subsection \ref{iod}). Observe that satisfying the RLGE condition is a generic case for $u$.

Here is the theorem regarding the asymptotic condition numbers $K_\infty(t,y_0,\widehat{z}_0)$ and $K_\infty(t,y_0)$, where we set $Q_1(t):=Q_{1L_1}(t)$ to simplify the notation.
\begin{theorem} \label{RLGE1}
	If $y_0$ and $\widehat{z}_0$ satisfy the RLGE condition, then
	$$
	K_\infty\left( t,y_{0},\widehat{z}_{0}\right)=\frac{\left\Vert
		Q_{1}(t)\widehat{z}_{0}\right\Vert }{\left\Vert Q_{1}(t)\widehat{y}%
		_{0}\right\Vert }$$
		and
		$$
		K_\infty\left( t,y_{0}\right)=\frac{\left\Vert
		Q_{1}(t)\right\Vert }{\left\Vert Q_{1}(t)\widehat{y}%
		_{0}\right\Vert },
	$$
	where
	$$
	Q_1(t):=Q_{1L_1}(t)=\sum\limits_{\substack{ \lambda _{i}\in \Lambda _1  \\ %
			m_{i}=L_1+1}}\mathrm{e}^{\mathrm{i}\omega _{i}t}P_{iL_1}.
	$$

	In terms of the JCF of $A$, we have
	$$
	Q_{1}(t)=\sum\limits_{\substack{\lambda _{i}\in \Lambda _1\\ j^\prime\in\{1,\ldots,d_i\}\\ m_{ij^\prime}=L_1+1}}\mathrm{e}^{\mathrm{i}\omega _{i}t}v^{%
			\left( i,j^\prime,1\right) }w^{\left( i,j^\prime,L_1+1\right) }
	$$
	and
		$$
	Q_{1}(t)u=\sum\limits_{\substack{\lambda _{i}\in \Lambda _1\\ j^\prime\in\{1,\ldots,d_i\}\\ m_{ij^\prime}=L_1+1}}\mathrm{e}^{\mathrm{i}\omega _{i}t}\alpha_{ij^\prime L_1+1}(u)
	v^{\left( i,j^\prime,1\right) },\ u\in\mathbb{C}^n,
	$$
	where (see Appendix \ref{JCFsection}):
	\begin{itemize}
		\item the sum
	$$
	\sum\limits_{\substack{\lambda _{i}\in \Lambda _1\\ j^\prime\in\{1,\ldots,d_i\}\\ m_{ij^\prime}=L_1+1}}
	$$
	is over  the longest Jordan chains of the rightmost eigenvalues of $A$, whose length is $L_1+1$;	here, $d_i$ is the number of Jordan chains of $\lambda_i$ (equal to the geometric multiplicity of $\lambda_i$), $j^\prime$ is an index of Jordan chain of $\lambda_i$ and $m_{ij^\prime}$ is the length of the $j^\prime$-th Jordan chain of $\lambda_i$;
	\item $v^{\left( i,j^\prime,1\right) }$ are the eigenvectors of these longest Jordan chains (they are column vectors);
	\item $\alpha_{ij^\prime L_1+1}(u)$, $u\in\mathbb{C}^n$, are the components of $u$ along the last generalized eigenvectors of these longest Jordan chains;
	\item $w^{( i,j^\prime,L_1+1)}$ are the last left generalized eigenvectors  corresponding  to these longest Jordan chains (they are row vectors).
\end{itemize}
\end{theorem}

\begin{proof}
Recall (\ref{Kgeneric}) and (\ref{Kgeneric2}). By the definition (\ref{Qjl}), we have
$$
Q_1(t)=Q_{1L_1}(t)=\sum\limits_{\substack{ \lambda _{i}\in \Lambda _1  \\ %
		m_{i}\geq L_1+1}}\mathrm{e}^{\mathrm{i}\omega _{i}t}P_{iL_1}=\sum\limits_{\substack{ \lambda _{i}\in \Lambda _1  \\ %
		m_{i}= L_1+1}}\mathrm{e}^{\mathrm{i}\omega _{i}t}P_{iL_1}
$$
since $L_1+1=\max\limits_{\lambda _{i}\in \Lambda _{1}}m_{i}$.

For the part regarding the JCF, by using
$$
P_{iL_1}=\sum\limits
_{\substack{ j^\prime=1  \\ m_{ij^\prime}= L_1+1}}^{d_{i}}v^{%
	\left( i,j^\prime,1\right) }w^{\left( i,j^\prime,L_1+1\right) }\text{\ \ and\ \ }P_{iL_1}u=\sum\limits
	_{\substack{ j^\prime=1  \\ m_{ij^\prime}=L_1+1}}^{d_{i}}\alpha_{ij^\prime L_1+1}(u) v^{\left( i,j^\prime,1\right) }
$$
(see Proposition \ref{lastexample} in Appendix \ref{JCFsection}), we obtain
\begin{eqnarray*}
	Q_{1}(t)&=&\sum\limits_{\substack{ \lambda _{i}\in \Lambda _1  \\ %
					m_{i}=L_1+1}}\mathrm{e}^{\mathrm{i}\omega _{i}t}\sum\limits
	_{\substack{ j^\prime=1  \\ m_{ij^\prime}= L_1+1}}^{d_{i}}v^{%
			\left( i,j^\prime,1\right) }w^{\left( i,j^\prime,L_1+1\right) }\\
			&=&\sum\limits_{\substack{\lambda _{i}\in \Lambda _1\\ j^\prime\in\{1,\ldots,d_i\}\\ m_{ij^\prime}=L_1+1}}\mathrm{e}^{\mathrm{i}\omega _{i}t}v^{%
				\left( i,j^\prime,1\right) }w^{\left( i,j^\prime,L_1+1\right) } \label{QQQ}
\end{eqnarray*}
and
\begin{eqnarray*}
	Q_{1}(t)u&=&\sum\limits_{\substack{ \lambda _{i}\in \Lambda _1   \\ %
			m_{i}=L_1+1}}\mathrm{e}^{\mathrm{i}\omega _{i}t}\sum\limits
	_{\substack{ j^\prime=1  \\ m_{ij^\prime}=L_1+1}}^{d_{i}}\alpha_{ij^\prime L_1+1}(u) v^{\left( i,j^\prime,1\right) }\\
	&=&\sum\limits_{\substack{\lambda _{i}\in \Lambda _1\\ j^\prime\in\{1,\ldots,d_i\}\\ m_{ij^\prime}=L_1+1}}\mathrm{e}^{\mathrm{i}\omega _{i}t}\alpha_{ij^\prime L_1+1}(u)
	v^{\left( i,j^\prime,1\right) }.  \notag
\end{eqnarray*}
\end{proof}

The previous theorem has the following corollary, which considers a particular situation regarding $\Lambda_1$.
\begin{corollary} \label{cor}
Suppose that  $\Lambda_1$ consists of only a real eigenvalue $\lambda_1$. For $P_{1L_1}y_0\neq 0$ and $P_{1L_1}\widehat{z}_0\neq 0$, we have 
$$
K_\infty\left( t,y_{0},\widehat{z}_{0}\right)=K_\infty\left(y_{0},\widehat{z}_{0}\right):=\frac{\Vert P_{1L_1}\widehat{z}_{0}\Vert}{\Vert P_{1L_1}\widehat{y}_{0}\Vert}
$$
and
$$
K_\infty\left( t,y_{0}\right)=K_\infty\left(y_{0}\right):=\frac{\Vert P_{1L_1}\Vert}{\Vert P_{1{L_1}}\widehat{y}_{0}\Vert}.
$$

Now, in addition, suppose that there is a unique Jordan chain of $\lambda_1$ of maximum length $L_1+1$, say the $j^\prime$-th of the Jordan chains of $\lambda_1$.
For $w^{(1,j^\prime,L_1+1)}y_0\neq 0$ and $w^{(1,j^\prime,L_1+1)}\widehat{z}_0\neq 0$, we have
	$$
	K_\infty\left(y_{0},\widehat{z}_{0}\right)=\frac{\vert w^{(1,j^\prime,L_1+1)}\widehat{z}_{0}\vert}{\vert w^{(1,j^\prime,L_1+1)}\widehat{y}_{0}\vert}
	$$
and
$$
K_\infty\left(y_{0}\right)=\frac{\Vert w^{(1,j^\prime,L_1+1)}\Vert}{\vert w^{(1,j^\prime,L_1+1)}\widehat{y}_{0}\vert},
$$
where $\Vert w^{(1,j^\prime,L_1+1)} \Vert$ is the induced norm of the row $w^{(1,j^\prime,L_1+1)}$.
\end{corollary}
\begin{proof}
For the first part, observe that the RLGE condition for $u\in\mathbb{C}^n$, i.e., $j(u)=1$ and $L(u)=L_1$, is equivalent to $P_{1L_1}u\neq 0$, by recalling the definition of $j(u)$ and $L(u)$ in Subsection \ref{iod}. Moreover, we have
$
Q_{1}(t)=P_{1L_1}.
$

For the second part, observe that the expressions for $Q_1(t)$ and $Q_1(t)u$, $u\in\mathbb{C}^n$, given in Theorem \ref{RLGE1} become
$$
Q_1(t)=P_{1L_1}=v^{(1,j^\prime,1)}w^{(1,j^\prime,L_1+1)}
$$
and
\begin{equation}
Q_1(t)u=P_{1L_1}u=\alpha_{1j^\prime L_1+1}(u)v^{(1,j^\prime,1)}. \label{Lambda1real}
\end{equation}
Moreover, the RLGE condition for $u$ is
$$
\alpha_{1j^\prime L_1+1}(u)=w^{(1,j^\prime,L_1+1)}u\neq 0
$$
by (\ref{Lambda1real}). Therefore, for $w^{(1,j^\prime,L_1+1)}y_0\neq 0$ and $w^{(1,j^\prime,L_1+1)}\widehat{z}_0\neq 0$, we have
$$
K_\infty\left(y_{0},\widehat{z}_{0}\right)=\frac{\Vert \alpha_{1j^\prime L_1+1}(\widehat{z}_0)v^{(1,j^\prime,1)}\Vert}{\Vert \alpha_{1j^\prime L_1+1}(\widehat{y}_0)v^{(1,j^\prime,1)}\Vert}=\frac{\vert \alpha_{1j^\prime L_1+1}(\widehat{z}_0)\vert}{\vert \alpha_{1j^\prime L_1+1}(\widehat{y}_0)\vert}=\frac{\vert w^{(1,j^\prime,L_1+1)}\widehat{z}_{0}\vert}{\vert w^{(1,j^\prime,L_1+1)}\widehat{y}_{0}\vert}
$$
and
\begin{eqnarray*}
	K_\infty\left(y_{0}\right)&=&\frac{\Vert v^{(1,j^\prime,1)}\Vert \Vert w^{(1,j^\prime,L_1+1)}\Vert}{\Vert \alpha_{1j^\prime L_1+1}(\widehat{y}_0)v^{(1,j^\prime,1)}\Vert}=\frac{\Vert w^{(1,j^\prime,L_1+1)}\Vert}{\vert \alpha_{1j^\prime L_1+1}(\widehat{y}_0)\vert}=\frac{\Vert w^{(1,j^\prime,L_1+1)}\Vert}{\vert w^{(1,j^\prime,L_1+1)}\widehat{y}_{0}\vert}
\end{eqnarray*}
since
$$
\Vert P_{1L_1}\Vert=\Vert v^{(1,j^\prime,1)}\Vert \Vert w^{(1,j^\prime,L_1+1)}\Vert.
$$
\end{proof}

\begin{remark}
	If $\lambda_1$ in the corollary is non-defective, i.e., $L_1=0$, then $P_{1L_1}=P_{10}$ is the projection onto the eigenspace corresponding to $\lambda_i$.
	
	If $\lambda_1$ in the corollary is simple, i.e., $d_1=1$ and $L_1=0$, then $w^{(1,j^\prime,L_1+1)}=w^{(1,1,1)}$ is a left eigenvector corresponding to $\lambda_i$.  
\end{remark}
The next example considers the asymptotic condition number $K_\infty(y_0)$ in the particular situation considered in Corollary \ref{cor}.

\begin{example}

Consider the non-diagonalizable matrix
$$
A_1=V
\left[ 
\begin{array}{rrr}
	0 & 1 & 0 \\ 
	0 & 0 & 0 \\ 
	0 & 0 & -1%
\end{array}%
\right]V^{-1}
$$
and the diagonalizable matrix
$$
A_2=V
\left[ 
\begin{array}{rrr}
	-1 & 0 & 0 \\ 
	0 & 0 & 0 \\ 
	0 & 0 & -1%
\end{array}%
\right]V^{-1},
$$
where
\begin{equation*}
	V=\left[ 
	\begin{array}{rrr}
		0 & 0 & 1 \\ 
		1 & 0 & -1 \\ 
		-1 & 1 & 0%
	\end{array}%
	\right] \text{\ \ and\ \ }V^{-1}=\left[ 
	\begin{array}{ccc}
		1 & 1 & 0 \\ 
		1 & 1 & 1 \\ 
		1 & 0 & 0%
	\end{array}%
	\right]=\left[\begin{array} {c} w^{(1)}\\w^{(2)}\\ w^{(3)}\end{array}\right]
\end{equation*}
with $w^{(1)}$, $w^{(2)}$ and $w^{(3)}$ the rows of $V^{-1}$. 

For both matrices $A_1$ and $A_2$, $\Lambda_1$ consists of the real eigenvalue $0$ and for it there is a unique Jordan chain of length $L_1+1=2$ for $A_1$, and $L_1+1=1$ for $A_2$.

By Corollary \ref{cor}, for both matrices, the RLGE condition is $w^{\left(2\right) }y_0\neq 0$ and we have
\begin{equation*}
	K_{\infty }(y_{0})=\frac{\Vert w^{(2)}\Vert}{\vert w^{(2)}\widehat{y}_{0}\vert}=\frac{\Vert w^{(2)}\Vert \Vert y_0\Vert}{\vert w^{(2)}y_{0}\vert}.
\end{equation*}
Suppose that the $1$-norm is used as vector norm. We have $w^{(2)}=[\ 1\ \ 1\ \ 1\ ]$ and then
$$
\Vert w^{(2)} \Vert=\left\Vert \left(w^{(2)}\right)^T\right\Vert_\infty=1.
$$
Therefore, 
$$
K_{\infty }(y_{0})=\frac{\Vert y_0\Vert_1}{\vert y_{01}+y_{02}+y_{03}\vert}=\frac{\vert y_{01}\vert+\vert y_{02}\vert +\vert y_{03}\vert}{\vert y_{01}+y_{02}+y_{03}\vert}
$$

In Figure \ref{FigureJ}, we see, for three different initial values $y_{0}$, $K(t,y_{0})$ (blue line) and $K_\infty(y_0)$ (red line) for both matrices $A_1$ and $A_2$: $A_1$ on the left for $t\in[0,50]$  and $A_2$ on the right for $t\in[0,15]$.

We can observe a slower approach to $K_\infty(y_{0})$ for the matrix $A_1$ than for the matrix $A_2$. Indeed, for $A_1$ we have a $\frac{1}{t}$-convergence of $K(t,y_0)$ to $K_\infty(y_0)$ as $t\rightarrow +\infty$, instead of the exponential $\mathrm{e}^{-t}$-convergence valid for $A_2$ (see Theorem \ref{ThmKy} and Propositions \ref{Lemma2} and \ref{Lemma1} for $\epsilon(t)$ and $\epsilon(t,\widehat{y}_0)$).

\begin{figure}[tbp]
	\centering
	\begin{subfigure}[b]{1\textwidth}
		\centering
		\includegraphics[width=0.8\textwidth]{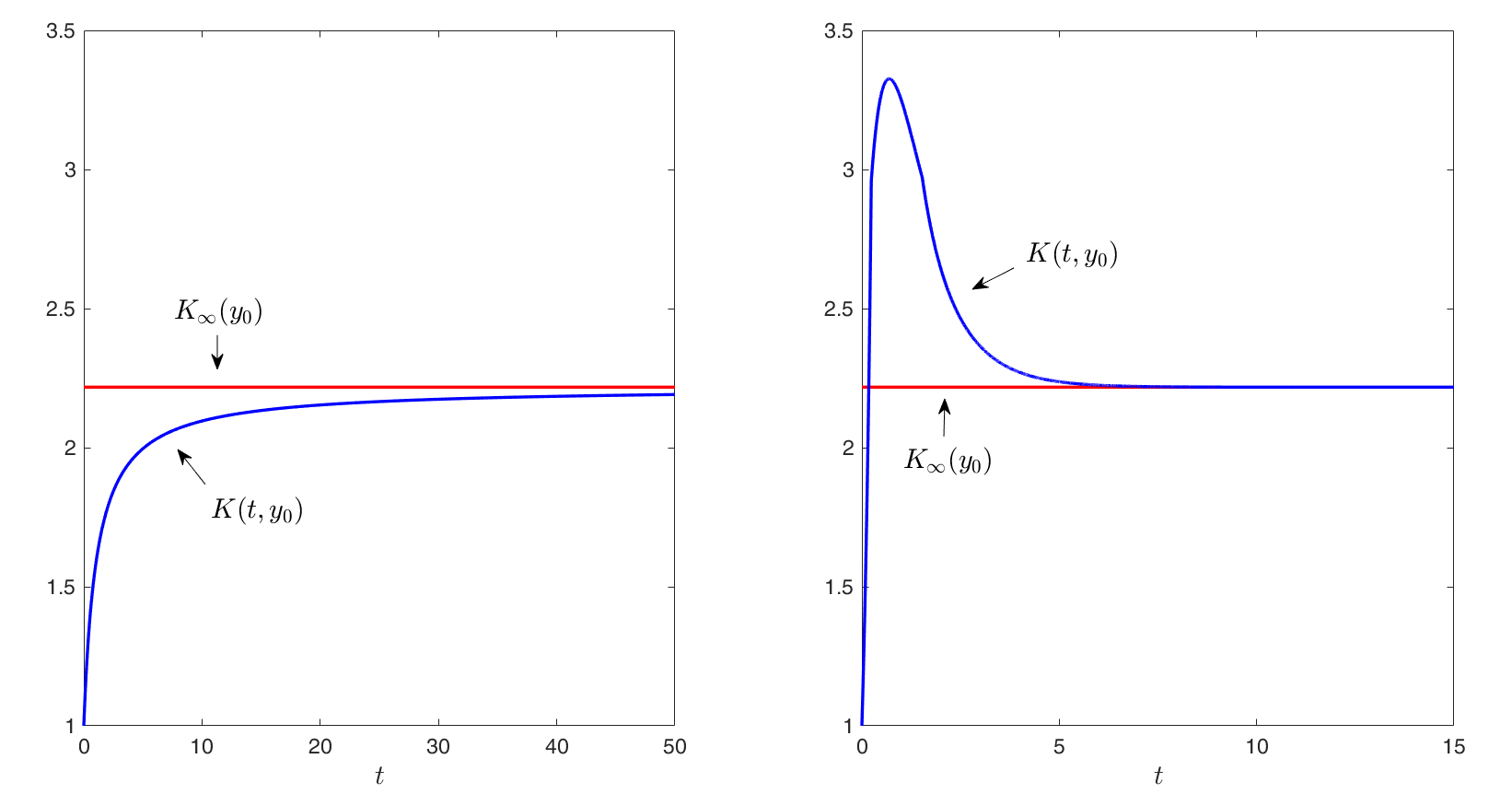}
		\caption{$y_0=(-0.8, -2.9, 1.4)$}
	\end{subfigure}
	\hfill  
	\begin{subfigure}[b]{1\textwidth}
		\centering
		\includegraphics[width=0.8\textwidth]{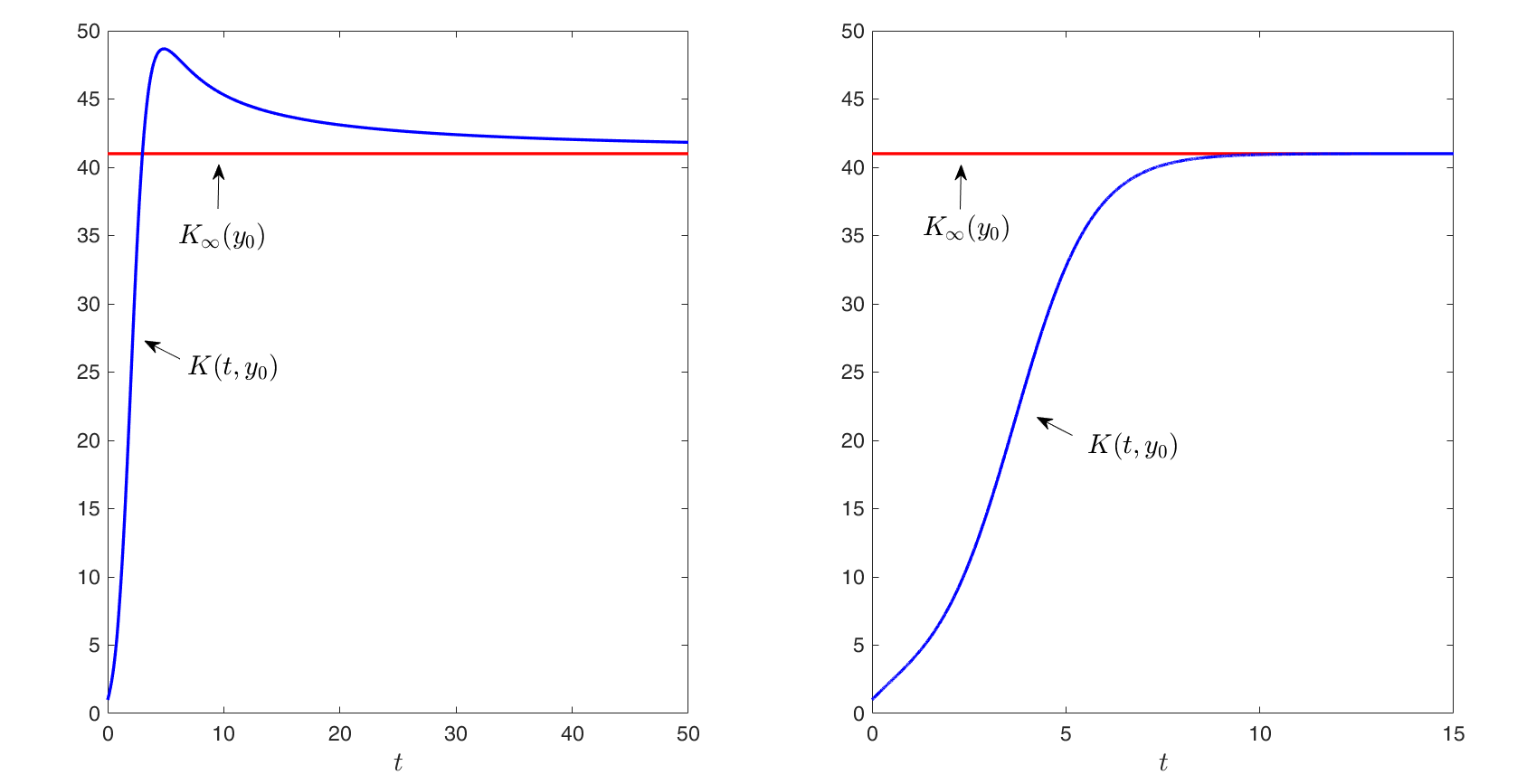}
		\caption{$y_0=(-1,1,0.05)$}
	\end{subfigure}
	
	\begin{subfigure}[b]{1\textwidth}
		\centering
		\includegraphics[width=0.8\textwidth]{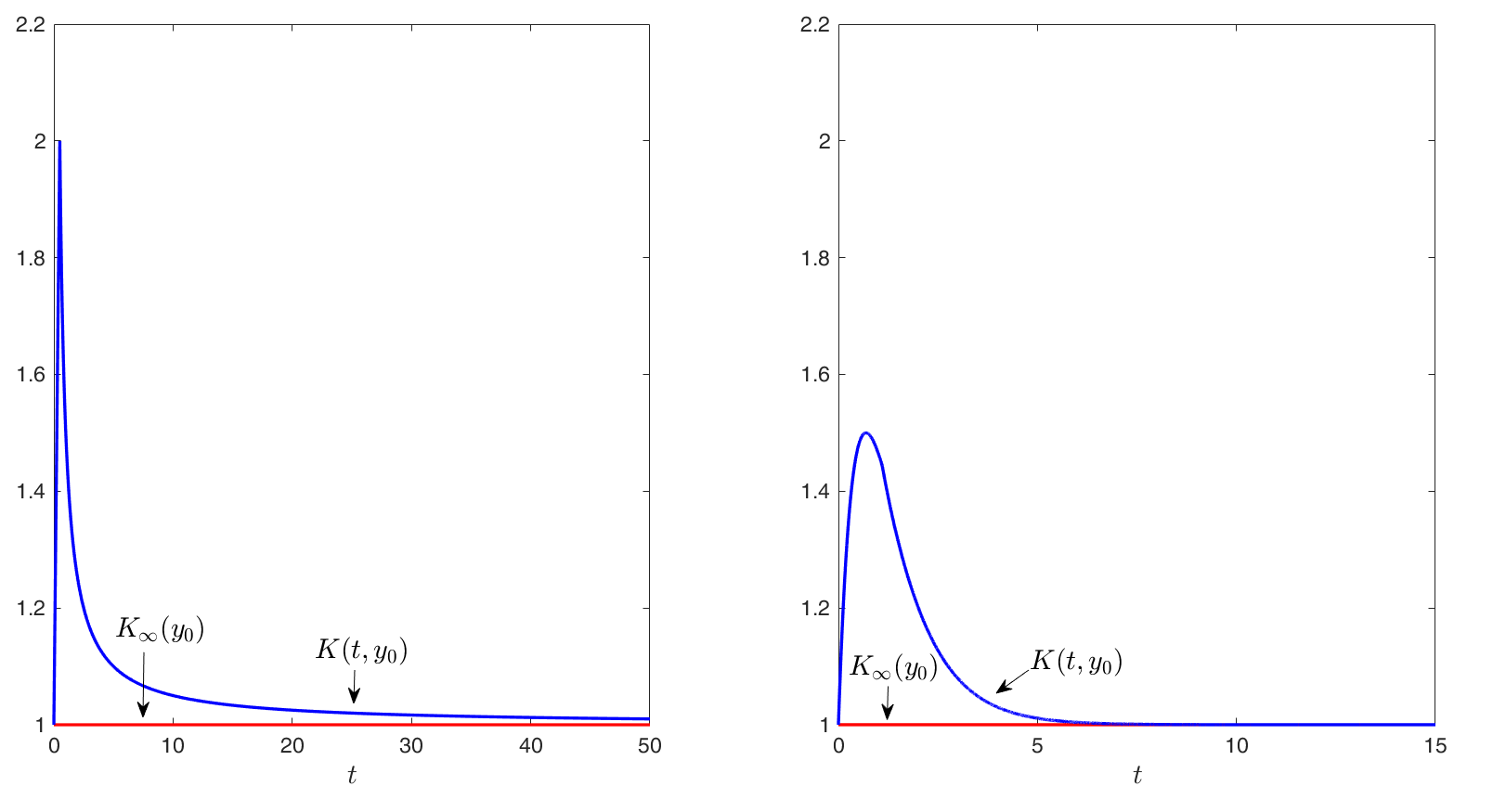}
		\caption{$y_0=(1,2,3)$}
	\end{subfigure}
	\hfill
	\caption{Condition numbers $K(t,y_0)$ and $K_\infty(y_0)$ for $A_1$ (left, $t\in[0,50]$) and $A_2$ (right, $t\in[0,15]$). The vector norm is the $1$-norm.}
	\label{FigureJ}
\end{figure}
\end{example}

\section{Conclusion}

In the present paper, we have considered how a perturbation in the initial value $y_0$ of the ODE (\ref{ode}) is propagated to the solution over a long time, by measuring the perturbation with a normwise relative error. In other words, we have studied the long time relative conditioning of the problem (\ref{due}).

We have defined:
\begin{itemize}
	\item a \emph{directional pointwise condition number} $K(t,y_0,\widehat{z}_0)$ such that
$$
K(t,y_0,\widehat{z}_0)=\frac{\delta(t)}{\varepsilon},
$$
where $\widehat{z}_0$ is the direction of the perturbation, $\varepsilon$ is the normwise relative error of the perturbed initial value and $\delta(t)$ is the normwise relative error of the perturbed solution;
\item  a \emph{pointwise condition number} $K(t,y_0)$, worst $K(t,y_0,\widehat{z}_0)$ by varying $\widehat{z}_0$;
\item  a \emph{global condition number} $K(t)$,  worst $K(t,y_0,\widehat{z}_0)$ by varying both $\widehat{z}_0$ and $y_0$.
\end{itemize}

Here is a summary of the results obtained in the present paper.

\begin{itemize}
	\item Determination of the \emph{asymptotic condition numbers} $K_\infty(t,y_0,\widehat{z}_0)$, $K_\infty(t,y_0)$ and $K_\infty^+(t)$, i.e., the asymptotic forms of $K(t,y_0,\widehat{z}_0)$, $K(t,y_0)$ and $K(t)$, respectively. See Theorems \ref{ThmKyz}, \ref{ThmKy} and \ref{ThmK}.
	
	\item \enquote{Asymptotic} and \enquote{worst} commute for $\widehat{z}_0$, but not for $y_0$.  More precisely, we have what follows.
	
	\begin{itemize}
		\item $K_\infty(t,y_0)$, the asymptotic form of the worst $K(t,y_0,\widehat{z}_0)$ by varying $\widehat{z}_0$, coincides with the worst asymptotic form $K_\infty(t,y_0,\widehat{z}_0)$  by varying $\widehat{z}_0$. See Theorem \ref{wcy}.
		
		\item On the contrary, $K_\infty^+(t)$, the asymptotic form of the worst $K(t,y_0)$ by varying  $y_0$, does not coincide with the worst asymptotic form $K_\infty(t,y_0)$ by varying $y_0$. This worst asymptotic $K_\infty(t,y_0)$, the asymptotic condition number $K_\infty(t)$, is smaller than $K_\infty^+(t)$ and it is determined in Theorem \ref{asworst2}.
	\end{itemize} 

\item Considering the worst-case scenario for $y_0$ and $\widehat{z}_0$ is too pessimistic. More precisely, we have what follows.
\begin{itemize}

\item The asymptotic condition number $K_\infty(t)$ is an exponentially diverging function of $t$ (provided that the eigenvalues of $A$ do not all have the same real part).

\item On the contrary, in the generic case where $y_0$ and $\widehat{z}_0$ satisfy the \emph{Rightmost Last Generalized Eigenvector} (RLGE) condition introduced in Subsection \ref{RLGE}, the asymptotic condition numbers $K_\infty(t,y_0,\widehat{z}_0)$ and $K_\infty(t,y_0)$ remain bounded as well as away from zero by varying $t$.  Expressions for $K_\infty(t,y_0,\widehat{z}_0)$ and $K_\infty(t,y_0)$ valid in the RLGE condition are presented in Theorem \ref{RLGE1} and Corollary \ref{cor}.
\end{itemize}

\end{itemize}

The present paper has two sequels. The first is \cite{M2}, which develops the results of the present paper (valid for a \emph{complex} ODE (\ref{ode})) in depth for a \emph{real} ODE (\ref{ode}) in a generic case for $A$, $y_0$ and $\widehat{z}_0$. In particular, the case where $A$ has complex conjugate rightmost eigenvalues is analyzed in depth. The second is \cite{M3}, which presents extensive experimental tests, applications to real-world systems, and other issues, such as the non-asymptotic behavior of the condition numbers and how rapidly the asymptotic behavior sets in, also in relation to the non-normality of the matrix $A$. Papers \cite{M2} and \cite{M3} are further developments that cannot be included here for space constraints.

Beyond \cite{M2} and \cite{M3}, the present paper may open the door to several important extensions, such as:
\begin{itemize}
	\item generalizations to infinite-dimensional settings, such as partial differential equations or delay differential equations, where the finite sum (\ref{expexp}) is replaced by a series;
	\item 	the treatment of componentwise relative errors, observables and structured perturbations, where we study the relative conditioning of the problem
	$$
	y_0\mapsto C\mathrm{e}^{tA}y_0,
	$$
	where $C$ is a matrix, and perturbations of $y_0$ are of the form $B\xi_0$, with $B$ a matrix defining the perturbation structure and $\xi_0$ a vector of free parameters within that structure, belonging to a space whose dimension is smaller than $n$, the dimension of the space of $y_0$;
	\item the numerical computation of all quantities involved in the asymptotic condition number $K_\infty(t,y_0)$ (see its expression in Theorem \ref{ThmKy});
	\item backward error analysis of numerical methods for the ODE (\ref{ode}), where a numerical solution is interpreted as an exact perturbed solution, as in \cite{Ma2021}, where the case of normal $A$ is considered.
	\item the analysis of the transient behavior of the relative conditioning of the problem (\ref{due}), where pseudospectra are not applicable, since they cannot provide a lower bound for $\Vert \mathrm{e}^{tA}\widehat{y}_0\Vert$ in (\ref{KAy0}), and different tools must be developed.
\end{itemize}

\bigskip

\noindent {\bf Acknowledgements:} the research was supported by the INdAM
Research group GNCS (Gruppo Nazionale di Calcolo Scientifico).

\newpage

\appendix
\clearpage
\pagenumbering{arabic}
\setcounter{page}{1}

\section{Linear Algebra Results} \label{JCFsection}

We introduce linear algebra notations and results necessary for the analysis developed in the present paper. Although the content pertains to the well-known topic of the Jordan Canonical Form (JCF) and matrix functions, it addresses very specific aspects that are either not known or not sufficiently detailed for the purposes of the present paper.

The main goal of the section is to derive a formula for the matrix exponential $\mathrm{e}^{tA}$, where it is explicitly
identified how $\mathrm{e}^{tA}$ depends on $t$. Since it is
based on the JCF of $A$, first we revise such a form.

\subsection{The JCF of the matrix $A$} \label{JCFsubsection}

Let $A\in\mathbb{C}^{n\times n}$ and let $\lambda _{1},\ldots,\lambda _{p}$
be the distinct eigenvalues of $A$. The matrix $A$ is similar to a matrix $%
J\in \mathbb{C}^{n\times n}$, called a \emph{Jordan Canonical Form (JCF)} of 
$A$, with the following structure.
\begin{itemize}
	\item The matrix $J$ is block-diagonal with $p$ blocks $J^{\left( 1\right)
	},\ldots ,J^{\left( p\right) }$ called \emph{Jordan blocks}:%
	\begin{equation*}
		J=\mathrm{diag}\left(J^{(1)},\ldots,J^{(p)}\right)\in\mathbb{C}^{n\times n}.
	\end{equation*}%
	\item For any $i\in\{1,\ldots ,p\}$, the Jordan block $J^{\left( i\right) }$ has dimension $%
	\nu _{i}$, where $\nu _{i}$ is the algebraic multiplicity of $\lambda _{i}$,
	and it is block-diagonal with $d_{i}$ blocks $J^{\left( i,1\right) },\ldots
	,J^{\left( i,d_{i}\right) }$ called \emph{Jordan mini-blocks}, where $d_{i}$
	is the geometric multiplicity of $\lambda _{i}$: 
	\begin{equation*}
		J^{(i)}=\mathrm{diag}\left(J^{\left( i,1\right) },\ldots
		,J^{\left( i,d_{i}\right)}\right)\in\mathbb{C}^{\nu_i\times\nu_i}.
	\end{equation*}%
	\item For any $i\in\{1,\ldots ,p\}$ and for any $j\in\{1,\ldots ,d_{i}\}$, the Jordan mini-block $%
	J^{\left( i,j\right) }$ is upper bidiagonal: 
	\begin{equation*}
		J^{\left( i,j\right) }=\left[ 
		\begin{array}{ccccc}
			\lambda _{i} & 1 &  &  &  \\ 
			& . & 1 &  &  \\ 
			&  & . & 1 &  \\ 
			&  &  & . & 1 \\ 
			&  &  &  & \lambda _{i}%
		\end{array}%
		\right]\in\mathbb{C}^{m_{ij}\times m_{ij}} .
	\end{equation*}%
	The dimension of $J^{\left( i,j\right) }$ is denoted by $m_{ij}$.
\end{itemize}
There is a unique JCF of $A$, except for permutations of the blocks and permutations of the mini-blocks within the blocks.

For any $i\in\{1,\ldots,p\}$, we have
$$
\sum\limits_{j=1}^{d_{i}}m_{ij}=\nu _{i}
$$
and we call
\begin{equation}
	m_{i}:=\max_{j\in\{1,\ldots ,d_{i}\}}m_{ij}  \label{ascent}
\end{equation}%
the \emph{ascent} of $\lambda _{i}$.

An eigenvalue $\lambda_i$, $i\in\{1,\ldots,p\}$, is called \emph{defective} if $d_i<\nu_i$ and \emph{non-defective} if $d_i=\nu_i$. Clearly, $\lambda_i$ is non-defective if and only if $m_i=1$, i.e., $J^{(i)}$ is diagonal. The matrix $A$ is diagonalizable if and only if all the eigenvalues $\lambda_1,\ldots,\lambda_p$ are non-defective, i.e., $J$ is diagonal.

An eigenvalue $\lambda_i$, $i\in\{1,\ldots,p\}$, is called \emph{simple} if $d_i=\nu_i=1$, i.e., $\lambda_i$ has a unique mini-block of dimension $1$.

\subsubsection{Three-index notation}

\begin{figure}[tbp]
	\includegraphics[width=0.5\textwidth]{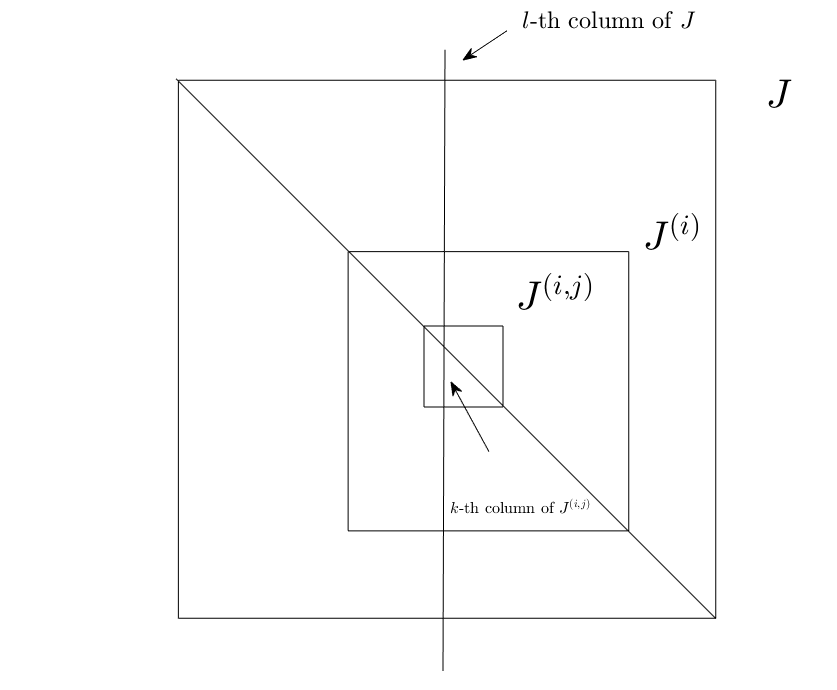}
	\caption{The index $l$ and the three indices $(i,j,k)$.}
	\label{threeindeces} 
\end{figure}

In the following, we denote an index $l\in\{1,\ldots,n\}$ by three indices $%
\left( i,j,k\right) $, where (see Figure \ref{threeindeces}):

\begin{itemize}
	\item $i\in\{1, \ldots, p\}$ is the index of the block $J^{(i)}$ traversed by the $l$-th column (or row) of $J$;
	
	\item $j\in \left\{ 1,\ldots,d_{i}\right\} $ is the index of the mini-block $%
	J^{(i,j)}$ of $J^{(i)}$ traversed by the $l$-th column (or row) of $J$;
	
	\item $k\in \left\{ 1,\ldots,m_{ij}\right\} $ is the index of the column (or row) of
	the mini-block $J^{(i,j)}$ included in the $l$-th column (or row) of $J$.
\end{itemize}

Observe that the triples $(i,j,k)$ appear lexicographically ordered when the index $l$ moves from $1$ to $n$.

\subsubsection{The Jordan basis}

Since $A$ is similar to $J$, there exists $V\in \mathbb{C}^{n\times n}$
non-singular such that 
\begin{equation}
	J=V^{-1}AV.  \label{JCF}
\end{equation}%
Let
$$
V=\left[v^{(1)}\ \cdots\ v^{(n)}\right].
$$
The $n$ columns $v^{(1)},\ldots, v^{(n)}$ of $V$ constitute a \emph{Jordan basis} of the
space $\mathbb{C}^{n}$.
In the three-index notation, the $n$ columns $%
v^{(l)}$, $l\in\{1,\ldots,n\}$, are denoted by 
\begin{equation}
	v^{\left( i,j,k\right) },\ i\in\{1,\ldots ,p\},\ j\in\{1,\ldots ,d_{i}\}\text{\ and\ }%
	k\in\{1,\ldots ,m_{ij}\}.  \label{JB}
\end{equation}
The vectors (\ref{JB}) appear as columns of $V$ in lexicographic order. 

For $i\in\{1,\ldots,p\}$, the vectors
\begin{equation*}
	v^{\left( i,j,1\right) },\ j\in\{1,\ldots ,d_{i}\},
\end{equation*}
are \emph{eigenvectors} corresponding to the eigenvalue $\lambda_i$: they constitute a basis for the \emph{eigenspace} corresponding to the eigenvalue $\lambda_i$. The vectors
\begin{equation*}
	v^{\left( i,j,k\right) },\ j\in\{1,\ldots ,d_{i}\}\text{\ and\ }%
	k\in\{1,\ldots ,m_{ij}\},
\end{equation*}
are \emph{generalized eigenvectors} corresponding to the eigenvalue $\lambda_i$: they constitute a basis for the \emph{generalized eigenspace} corresponding to the eigenvalue $\lambda_i$.

\subsubsection{Jordan chains}
The Jordan basis (\ref{JB}) is partitioned into the \emph{Jordan chains} 
\begin{eqnarray*}
	&&\left(v^{\left( i,j,k\right)}\right)_{k=1,\ldots,m_{ij}}=\left(v^{\left( i,j,1\right)},v^{\left( i,j,2\right)},\ldots,v^{\left(
		i,j,m_{ij}\right)}\right)\\
	&&i\in\{1,\ldots,p\}\text{\ and\ }j\in\{1,\ldots,d_i\}.
\end{eqnarray*}

The elements of a
Jordan chain satisfy 
\begin{eqnarray}  \label{ge}
	&&\left(A-\lambda_iI\right)v^{\left( i,j,1\right)}=0  \notag \\
	&&\left(A-\lambda_iI\right)v^{\left( i,j,k+1\right)}=v^{\left(
		i,j,k\right)},\ k\in\{1,\ldots,m_{ij}-1\}.  \notag \\
\end{eqnarray}
The chain stops with $v^{\left( i,j,m_{ij}\right)}$ since the system 
\begin{equation*}
	\left(A-\lambda_iI\right)x=v^{\left( i,j,m_{ij}\right)}
\end{equation*}
has no solution.

Observe that (\ref{JCF}) is equivalent to have (\ref{ge}), for all $i\in\{1,\ldots,p\}$ and $j\in\{1,\ldots,d_i\}$, and that there is a correspondence one-to-one between mini-blocks and Jordan chains.
\subsubsection{The matrix $zA$}

The next proposition describes the JCF and a Jordan basis of the matrix $%
z A$, where $z\in\mathbb{C}\setminus\{0\}$, in terms of the JCF
and a Jordan basis of $A$.

\begin{proposition}
	\label{Vfor-A} 
	Let $z\in\mathbb{C}\setminus\{0\}$. The
	distinct eigenvalues of $z A$ are $z\lambda_i$,\ $i\in\{1,\ldots,p\}$.
	For any $i\in\{1,\ldots,p\}$, the number and the dimensions of the mini-blocks corresponding to the eigenvalue $z\lambda_i$ in the JCF of $z A
	$ are equal to the number and the dimensions of the mini-blocks corresponding to the eigenvalue $\lambda_i$ in the JCF of $A$. A Jordan basis of $zA$ is 
	\begin{equation}
		z^{-(k-1)}v^{\left( i,j,k\right) },\ i\in\{1,\ldots ,p\},\ j\in\{1,\ldots ,d_{i}\}%
		\text{\ and\ }k\in\{1,\ldots ,m_{ij}\}.  \label{vv}
	\end{equation}
\end{proposition}
\begin{proof}
	Let $i\in\{1,\ldots,p\}$. Given a Jordan chain 
	\begin{equation*}
		\left(v^{\left( i,j,k\right)}\right)_{k=1,\ldots,m_{ij}}  \label{chainsstart}
	\end{equation*}
	of $A$ corresponding to the eigenvalue $\lambda_i$, we have that
	\begin{equation}
		\left(z^{-(k-1)}v^{\left( i,j,k\right) }\right)_{k=1,\ldots,m_{ij}} \label{chains}
	\end{equation}
	is a Jordan chain of $z A$ corresponding to the eigenvalue $z\lambda_i$.
	
	In fact,  by (\ref{ge})    we have
	\begin{eqnarray*}
		&&\left(z A-z\lambda_iI\right)v^{\left( i,j,1\right) }=0\\
		&&\left(z A-z\lambda_iI\right)z^{-k}v^{\left( i,j,k+1\right)
		}=z^{-(k-1)}v^{\left( i,j,k\right) },\ k\in\{1,\ldots,m_{ij}-1\}.
	\end{eqnarray*}
	Moreover, the system
	$$
	\left(z A-z\lambda_iI\right)x=z^{-(m_{ij}-1)}v^{\left( i,j,m_{ij}\right)}
	$$
	has no solution; otherwise, the system
	$$
	\left(A-\lambda_iI\right)y=v^{\left( i,j,m_{ij}\right)}
	$$
	would have the solution $y = z^{m_{ij}} x$.
	
	By exchanging the role of $A$ and $zA$, i.e., we consider $zA$ and $z^{-1}(zA)$, we also see that, given a Jordan chain 
	\begin{equation*}
		\left(u^{\left( i,j,k\right)}\right)_{k=1,\ldots,m_{ij}} 
	\end{equation*}
	of $zA$ corresponding to the eigenvalue $z\lambda_i$, 
	\begin{equation*}
		\left(z^{k-1}u^{\left( i,j,k\right) }\right)_{k=1,\ldots,m_{ij}}
	\end{equation*}
	is a Jordan chain of $z^{-1}(zA)=A$ corresponding to the eigenvalue $\lambda_i$. Here, $j$ is the index of the particular Jordan chain of $zA$ considered, corresponding to the eigenvalue $z\lambda_i$, and $m_{ij}$ denotes its length. This shows that the Jordan chains of $zA$ are all of type (\ref{chains}). Consequently, the number and the lengths of the Jordan chains corresponding to the eigenvalue $z\lambda_i$ of $z A$ are equal to the number and the lengths of the Jordan chains corresponding to the eigenvalue $\lambda_i$ of $A$. Therefore, the number and the dimensions of the mini-blocks corresponding to the eigenvalue $z\lambda_i$ of $z A
	$ are equal to the number and the dimensions of the mini-blocks corresponding to the eigenvalue $\lambda_i$ of $A$.
	
	A Jordan basis for $z A$ is given by collecting all the Jordan chains (\ref{chains}). Thus, we obtain the Jordan basis (\ref{vv}) for $z A$.
\end{proof}
\subsubsection{The real case}

When $A$ is a real matrix, the distinct complex eigenvalues $%
\lambda_1,\ldots,\lambda_p$ of $A$ are divided in real eigenvalues and
complex conjugate pairs of eigenvalues. The next proposition describes the JCF and a
Jordan basis of $A$, when $A$ is real, in terms of these real
eigenvalues and complex conjugate pairs of eigenvalues.

Here and in the following, for a vector or matrix $Z$, $\overline{Z}$
denotes the vector or matrix given by the complex conjugates of the elements
of $Z$.

\begin{proposition}
	\label{Prop22} Assume $A\in\mathbb{R}^{n\times n}$. For $i_{1},i_{2}\in\{1,%
	\ldots ,p\}$ with $i_{1}\neq i_{2}$ such that $\lambda _{i_{2}}=\overline{%
		\lambda _{i_{1}}}$, i.e., $\lambda _{i_{1}}$ and $\lambda _{i_{2}}$ form a
	complex conjugate pair of eigenvalues, we have 
	\begin{equation}
		\nu _{i_{2}}=\nu _{i_{1}},\ d_{i_{2}}=d_{i_{1}}\text{\ and\ }
		m_{i_{2}j}=m_{i_{1}j},\ j\in\{1,\ldots,d_{i_{2}}\}.  \label{R2}
	\end{equation}
	(Indeed, in (\ref{R2}) we should say that there exists an ordering of the mini-blocks corresponding to $\lambda_{i_2}$ such that
	$m_{i_2j}=m_{i_1j}$, $j\in\{1,\ldots,d_{i_2}\}$). Moreover, there exists a Jordan basis of $A$ such that:
	\begin{itemize}
		\item for $i\in\{1,\ldots,p\}$ such that $\lambda_{i}\in\mathbb{R}$, we have 
		\begin{equation*}
			v^{\left( i,j,k\right) }\in\mathbb{R}^n,\ j\in\{1,\ldots ,d_{i}\}\text{\ and\ }%
			k\in\{1,\ldots ,m_{ij}\};  \label{conjugate}
		\end{equation*}
		
		\item for $i_{1},i_{2}\in\{1,\ldots ,p\}$ with $i_{1}\neq i_{2}$ such that $%
		\lambda _{i_{2}}=\overline{\lambda _{i_{1}}}$, we have 
		\begin{equation}
			v^{\left( i_{2},j,k\right) }=\overline{v^{\left( i_{1},j,k\right) }},\
			j\in\{1,\ldots ,d_{i_{2}}\}\text{\ and\ }k\in\{1,\ldots
			,m_{i_{2}j}\}.  \label{R3}
		\end{equation}
		Observe that $d_{i_2}=d_{i_1}$ and $m_{i_2j}=m_{i_1j}$ hold in (\ref{R3}).
	\end{itemize}
\end{proposition}

\begin{proof}
	Consider a complex conjugate pair given by $\lambda_{i_1}$ and $\lambda _{i_{2}}=\overline{%
		\lambda _{i_{1}}}$. Given a Jordan chain 
	\begin{equation*}
		\left(v^{\left( i_1,j,k\right)}\right)_{k=1,\ldots,m_{i_1j}}  \label{chains2start}
	\end{equation*}
	of $A$ corresponding to the eigenvalue $\lambda_{i_1}$, we have that
	\begin{equation}
		\left(\overline{v^{\left( i_1,j,k\right) }}\right)_{k=1,\ldots ,m_{i_1j}}
		\label{chains2}
	\end{equation}
	is a Jordan chain of $A$ corresponding to the eigenvalue $\overline{\lambda_{i_1}}=\lambda_{i_2}$.
	
	In fact, by conjugating both sides in (%
	\ref{ge}), we have
	\begin{eqnarray*}
		&&\left(A-\overline{\lambda_{i_1}}I\right)\overline{v^{\left( i_1,j,1\right) }}=0\\
		&&\left(A-\overline{\lambda_{i_1}}I\right)\overline{v^{\left( i_1,j,k+1\right)
		}}=\overline{v^{\left( i_1,j,k\right) }},\ k\in\{1,\ldots,m_{i_1j}-1\}.
	\end{eqnarray*}
	Moreover, the system
	$$
	\left(A-\overline{\lambda_{i_1}}I\right)x=\overline{v^{\left( i_1,j,m_{i_1j}\right)}}
	$$
	has no solution; otherwise, the system
	$$
	\left(A-\lambda_{i_1}I\right)y=v^{\left( i_1,j,m_{i_1j}\right)}
	$$
	would have the solution $y =\overline{x}$.
	
	By exchanging the role of $\lambda_{i_1}$ and $\lambda_{i_2}=\overline{\lambda_{i_1}}$, i.e., we consider $\lambda_{i_2}$ and $\overline{\lambda_{i_2}}=\lambda_{i_1}$, we also see that, given a Jordan chain 
	\begin{equation*}
		\left(u^{\left( i_2,j,k\right)}\right)_{k=1,\ldots,m_{i_2j}} 
	\end{equation*}
	corresponding to the eigenvalue $\lambda_{i_2}$, 
	\begin{equation*}
		\left(\overline{u^{\left( i_2,j,k\right) }}\right)_{k=1,\ldots,m_{i_2j}}
	\end{equation*}
	is a Jordan chain corresponding to the eigenvalue $\overline{\lambda_{i_2}}=\lambda_{i_1}$. This shows that the Jordan chains corresponding to $\lambda_{i_2}=\overline{\lambda_{i_1}}$ are all of type (\ref{chains2}). This implies (\ref{R2}). By collecting all the Jordan chains (\ref{chains2}) corresponding to eigenvalues $\lambda_{i_2}$, we obtain a Jordan
	basis satisfying  (\ref{R3}).
	
	The result regarding a real eigenvalue $\lambda_i$ follows from the next
	two facts:
	
	\begin{itemize}
		\item the eigenspace of $\lambda_i$ in $\mathbb{C}^n$ has a basis of real
		eigenvectors;
		
		\item if $v^{(i,j,k)}$ is real and the linear system 
		\begin{equation*}
			(A-\lambda_iI)x=v^{(i,j,k)}
		\end{equation*}
		has a solution in $\mathbb{C}^n$, then it also has a solution in $\mathbb{R}%
		^n$.
	\end{itemize}
	Thus, we can have Jordan chains corresponding to $\lambda_i$ constituted by real
	generalized eigenvectors.
\end{proof}
In the following, in case of a real matrix $A$, we assume to have a Jordan
basis as that described in the previous Proposition \ref{Prop22}.

\subsection{The matrices $V^{(i)}$, $W^{(i)}$, $N^{(i,l)}$ and the vector $\alpha^{(i)}(u)$}

The formula of our interest for the matrix exponential $\mathrm{e}^{tA}$ is constructed by using the matrices $V^{(i)}$, $W^{(i)}$ and $N^{(i,l)}$ now introduced. Here,  $i\in\{1,\ldots,p\}$ and  $l\in\{0,\ldots,m_{i}-1\}$ (remind that $m_i$ is the ascent of $\lambda_i$ defined in (\ref{ascent})). We also introduce the vector $\alpha^{(i)}(u)$ of components of $u\in\mathbb{C}^n$ in the Jordan basis.

Recall that $V$ is the matrix whose columns (\ref{JB}) constitute a Jordan basis.

\begin{figure}[tbp]
	\includegraphics[width=0.6\textwidth]{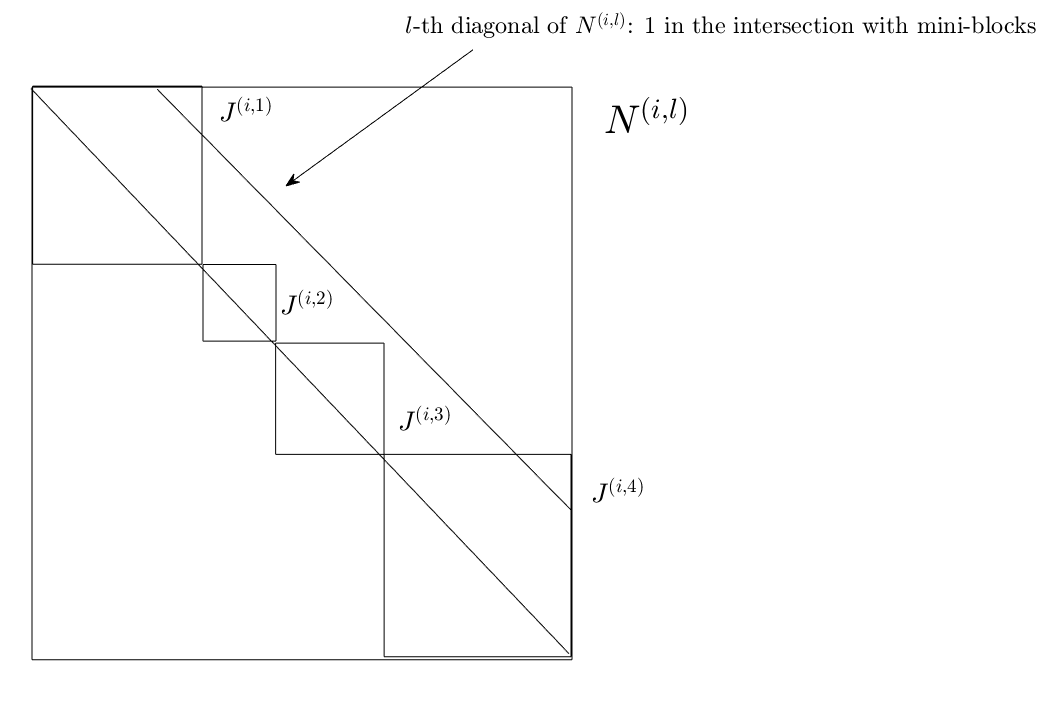}
	\caption{The matrix $N^{(i,l)}$.}
	\label{matrixN} 
\end{figure}

\begin{itemize}
	\item For $i\in\{1,\ldots ,p\}$ and  $j\in\{1,\ldots ,d_{i}\}$, let
	\begin{equation*}
		V^{(i,j)}:=\left[ v^{(i,j,1)}\ \cdots \ v^{(i,j,m_{ij})}\right] \in \mathbb{C%
		}^{n\times m_{ij}}
	\end{equation*}
	and let
	\begin{equation}
		V^{(i)}:=\left[ V^{(i,1)}\ \cdots \ V^{(i,d_{i})}\right] \in \mathbb{C}%
		^{n\times \nu _{i}}.  \label{this1}
	\end{equation}%
	Observe that
	\begin{equation*}
		V=\left[ V^{(1)}\ \cdots \ V^{(p)}\right].
	\end{equation*}
	
	\item Let $W:=V^{-1}$ and let 
	\begin{equation}
		w^{\left( i,j,k\right) },\ i\in\{1,\ldots ,p\},\ j\in\{1,\ldots ,d_{i}\}\text{\ and\ }%
		k\in\{1,\ldots ,m_{ij}\},  \label{ws}
	\end{equation}
	be the rows of $W$ in the three-index notation. They appear in $W$ in lexicographic order. The $n$ row vectors in (\ref{ws}) are called \emph{left generalized eigenvectors} of $A$, whereas, as we have already seen, the $n$  column vectors in (\ref{JB}) are called \emph{(right) generalized eigenvectors} of $A$. For $i\in\{1,\ldots ,p\}$ and $j\in\{1,\ldots ,d_{i}\}$, let
	\begin{equation*}
		W^{(i,j)}:=\left[ 
		\begin{array}{c}
			w^{(i,j,1)} \\ 
			\vdots \\ 
			w^{(i,j,m_{ij})}%
		\end{array}%
		\right] \in \mathbb{C}^{m_{ij}\times n}
	\end{equation*}
	and let 
	\begin{equation}
		W^{(i)}:=\left[ 
		\begin{array}{c}
			W^{(i,1)} \\ 
			\vdots \\ 
			W^{(i,d_{i})}%
		\end{array}%
		\right] \in \mathbb{C}^{\nu _{i}\times n}.  \label{this}
	\end{equation}%
	Observe that
	\begin{equation*}
		W=\left[ 
		\begin{array}{c}
			W^{(1)} \\ 
			\vdots \\ 
			W^{(p)}%
		\end{array}%
		\right].
	\end{equation*}

	\item For $u\in\mathbb{C}^n$, let 
	\begin{equation*}
		\alpha(u):=Wu.
	\end{equation*}
	Observe that $\alpha(u)$ is the vector of the components of $u$ in the
	Jordan basis. In the three-index notation, the components of $\alpha(u)$
	are 
	\begin{equation*}
		\alpha_{ijk}(u),\ i\in\{1,\ldots ,p\},\ j\in\{1,\ldots ,d_{i}\}\text{\ and\ }k\in\{1,\ldots
		,m_{ij}\},
	\end{equation*}
	$\alpha_{ijk}(u)$ being the component of $u$ along $v^{(i,j,k)}$.
	For $i\in\{1,\ldots ,p\}$ and $j\in\{1,\ldots ,d_{i}\}$, let
	\begin{equation*}
		\alpha^{(i,j)}(u):=\left[ 
		\begin{array}{c}
			\alpha_{ij1}(u) \\ 
			\vdots \\ 
			\alpha_{ijm_{ij}}(u)%
		\end{array}%
		\right] \in \mathbb{C}^{m_{ij}}  \label{this3}
	\end{equation*}
	and let 
	\begin{equation*}
		\alpha^{(i)}(u):=\left[ 
		\begin{array}{c}
			\alpha^{(i,1)}(u) \\ 
			\vdots \\ 
			\alpha^{(i,d_{i})}(u)%
		\end{array}%
		\right] \in \mathbb{C}^{\nu _{i}}. \label{this2}
	\end{equation*}%
	Observe that
	\begin{equation*}
		\alpha(u):=\left[ 
		\begin{array}{c}
			\alpha^{(1)}(u) \\ 
			\vdots \\ 
			\alpha^{(p)}(u)%
		\end{array}%
		\right].
	\end{equation*}
	
	\item For $i\in\{1,\ldots,p\}$, $l\in\{0,\ldots,m_{i}-1\}$ and $j\in\{1,\ldots,d_{i}\}$, let 
	\begin{equation}
		N^{(i,l,j)}:=\left\{
		\begin{array}{l}
			\left[ 
			\begin{array}{ccccccc}
				0 & . & 0 & 1 & 0 & . & 0 \\ 
				& . & . & . & . & . & . \\ 
				&  & . & . & . & . & 0 \\ 
				&  &  & . & . & . & 1 \\ 
				&  &  &  & . & . & 0 \\ 
				&  &  &  &  & . & . \\ 
				&  &  &  &  &  & 0%
			\end{array}%
			\right] \in \mathbb{C}^{m_{ij}\times m_{ij}}\text{ if }l\leq m_{ij}-1\\
			\\
			0\in \mathbb{C}^{m_{ij}\times m_{ij}}\text{ if }m_{ij}\leq
			l\leq m_i-1,
		\end{array}
		\right.
		\label{DA}
	\end{equation}%
	where the upper diagonal of elements equal to $1$ is the $l$-th upper
	diagonal.
	For $i\in\{1,\ldots,p\}$ and $l\in\{0,\ldots,m_i-1\}$, let 
	\begin{equation}
		N^{(i,l)}:=\mathrm{diag}\left( N^{(i,l,1)},\ldots
		,N^{(i,l,d_i)}\right) \in \mathbb{C}^{\nu _{i}\times \nu _{i}}.  \label{A}
	\end{equation}%
	Observe that the matrix $N^{(i,l)}$ has the same dimensions $\nu
	_{i}\times \nu _{i}$ of the Jordan block $J^{(i)}=\mathrm{diag}\left(
	J^{(i,1)},\ldots ,J^{(i,d_{i})}\right) $ and the diagonal blocks $%
	N^{(i,l,1)},$ $\ldots,$\ $N^{(i,l,d_{i})}$ have the same
	dimensions of the Jordan mini-blocks $J^{(i,1)},\ldots ,$ $J^{(i,d_{i})}$, respectively. Thus, as it is illustrated in Figure \ref{matrixN}, the matrix $N^{(i,l)}$ has $1$ in the intersection of the $l$-th upper diagonal with the (frame of the) Jordan mini-blocks and $0$ in all the other
	places.
\end{itemize}

\begin{remark}\label{N}
	For $i\in\{1,\dots,p\}$ and $j\in\{1,\ldots,d_i\}$, we have
	$$
	J^{(i,j)}=\lambda_i I_{m_{ij}}+N^{(i,1,j)}\text{\ \ and\ \ }J^{(i)}=\lambda_i I_{\nu_{i}}+N^{(i,1)}.
	$$
	Moreover, the matrices $N^{(i,l,j)}$ and $N^{(i,l)}$, $l\in\{0,\ldots,m_i-1\}$, are powers of the nilpotent matrices $N^{(i,1,j)}$ and $N^{(i,1)}$, respectively: we have 
	$$
	N^{(i,l,j)}=\left(N^{(i,1,j)}\right)^l\text{\ \ and\ \ }N^{(i,l)}=\left(N^{(i,1)}\right)^l.
	$$
\end{remark}

\subsubsection{The matrix $z A$}

When we replace $A$ by $z A$, $z\in\mathbb{C}\setminus\{0\}$, the number and the
dimensions of blocks and mini-blocks remain the same (recall Proposition \ref%
{Vfor-A}). Moreover, a Jordan basis for $z A$ is given in (\ref{vv}). Therefore, we know how the matrix $V$ is transformed  by replacing $A$ by $z A$. The next proposition states this and, in addition, how the matrix $W=V^{-1}$ is transformed.

\begin{proposition}
	\label{PropWfor-A} 
	Let $z\in\mathbb{C}\setminus\{0\}$. The matrix $V(zA)$
	corresponding to $z A$ has columns 
	\begin{equation}
		z^{-(k-1)}v^{\left( i,j,k\right) },\ i\in\{1,\ldots ,p\},\ j\in\{1,\ldots ,d_{i}\}
		\text{\ and\ }k\in\{1,\ldots ,m_{ij}\},  \label{w-11}
	\end{equation}
	and the matrix $W(zA)$
	corresponding to $z A$ has rows 
	\begin{equation}
		z^{k-1}w^{\left( i,j,k\right) },\ i\in\{1,\ldots ,p\},\ j\in\{1,\ldots ,d_{i}\}
		\text{\ and\ }k\in\{1,\ldots ,m_{ij}\}.  \label{w-1}
	\end{equation}
\end{proposition}

\begin{proof}
	The columns of the matrix $V(zA)$ are given in (\ref{vv}). For $i_1,i_2\in\{1,\ldots,p\}$, $j_1\in\{1,\ldots ,d_{i_1}\}$, $j_2\in\{1,\ldots ,d_{i_2}\}$, $k_1\in\{1,\ldots ,m_{i_1j_1}\}$ and $k_2\in\{1,\ldots ,m_{i_2j_2}\}$, we have 
	\begin{equation*}
		z^{k_1-1}w^{\left( i_1,j_1,k_1\right) }z^{-(k_2-1)}v^{\left(
			i_2,j_2,k_2\right) }=\left\{ 
		\begin{array}{l}
			1\ \ \text{if}\ \ (i_1,j_1,k_1)=(i_2,j_2,k_2) \\ 
			0\ \ \text{otherwise}%
		\end{array}
		\right.
	\end{equation*}
	since $W$, whose rows are (\ref{ws}), is the inverse of $V$, whose columns are (\ref{JB}). This shows that the matrix of rows (\ref{w-1}) is the inverse of the matrix of columns (\ref{w-11}).
\end{proof}

\subsubsection{The real case}

The next proposition describes, when the matrix $A$ and the vector $u$ are
real, the matrices $V^{(i)}$ and $W^{(i)}$ and the vector $\alpha^{(i)}(u)$
in terms of the real eigenvalues and the complex
conjugate pairs of eigenvalues of $A$.

\begin{proposition}
	\label{conjVW} Assume $A\in\mathbb{R}^{n\times n}$ and $u\in\mathbb{R}^n$. Moreover, assume we have a Jordan basis as that described in Proposition \ref{Prop22}.
	
	For $i\in\{1,\ldots,p\}$ such that $\lambda_i\in\mathbb{R}$, we have 
	\begin{equation*}
		V^{(i)}\in\mathbb{R}^{n\times\nu_i},\ W^{(i)}\in\mathbb{R}^{\nu_i\times n}%
		\text{\ \ and\ \ }\alpha^{(i)}(u)\in\mathbb{R}^{\nu_i}.
	\end{equation*}
	
	For $i_1,i_2\in\{1,\ldots,p\}$ with $i_1\neq i_2$ such that $\lambda_{i_2}=%
	\overline{\lambda_{i_1}}$, we have 
	\begin{equation*}
		V^{(i_2)}=\overline{V^{(i_1)}},\ W^{(i_2)}=\overline{W^{(i_1)}}\text{\ \
			and\ \ }\alpha^{(i_2)}(u)=\overline{\alpha^{(i_1)}(u)}.
	\end{equation*}
\end{proposition}

\begin{proof}
	Since we have a Jordan basis as that described in
	Proposition \ref{Prop22}, we have $V^{(i)}\in\mathbb{R}^{n\times\nu_i}$ for  $i\in\{1,\ldots,p\}$ such that $\lambda_i\in\mathbb{R}$ and  $V^{(i_2)}=\overline{V^{(i_1)}}$ for $i_1,i_2\in\{1,\ldots,p\}$ with $i_1\neq i_2$ such that $\lambda_{i_2}=%
	\overline{\lambda_{i_1}}$.
	
	Now, we prove that for $i\in\{1,\ldots,p\}$ such that $\lambda_i\in\mathbb{R}$ we have $\alpha^{(i)}(u)\in\mathbb{R}$, and for $i_{1},i_{2}\in\{1,\ldots ,p\}$ such that $\lambda
	_{i_{2}}=\overline{\lambda _{i_{1}}}$ with $\mathrm{Im}\left(\lambda_{i_1}\right)>0$ we have $\alpha ^{(i_{2})}(u)=%
	\overline{\alpha ^{(i_{1})}(u)}$. By conjugating both sides of 
	\begin{equation*}
		u=V\alpha (u)=\sum\limits_{i=1}^{p}V^{(i)}\alpha ^{(i)}(u),
	\end{equation*}%
	we obtain 
	\begin{equation*}
		u=\sum\limits_{i=1}^{p}\overline{V^{(i)}}\ \overline{\alpha ^{(i)}(u)}.
	\end{equation*}
	Therefore, by separating real eigenvalues and complex conjugate pairs of eigenvalues, we have%
	\begin{equation*}
		u=\sum\limits_{i=1}^{p}V^{(i)}\alpha ^{(i)}(u)=\sum\limits_{\lambda _{i}\in \mathbb{R}}V^{(i)}\alpha
		^{(i)}(u)+\sum\limits _{\substack{ \lambda _{i_{1}}\in \mathbb{C}  \\ 
				\mathrm{Im}\left( \lambda _{i_{1}}\right) >0}}\left( V^{(i_{1})}\alpha
		^{(i_{1})}(u)+V^{(i_{2})}\alpha ^{(i_{2})}(u)\right),
	\end{equation*}%
	as well as%
	\begin{equation*}
		u=\sum\limits_{i=1}^{p}\overline{V^{(i)}}\ \ \overline{\alpha ^{(i)}(u)}%
		=\sum\limits_{\lambda _{i}\in \mathbb{R}}V^{(i)}\ \overline{\alpha ^{(i)}(u)}%
		+\sum\limits_{\substack{ \lambda _{i_{1}}\in \mathbb{C}  \\ \mathrm{Im}%
				\left( \lambda _{i_{1}}\right) >0}}\left( V^{(i_{2})}\overline{\alpha
			^{(i_{1})}(u)}+V^{(i_{1})}\overline{\alpha ^{(i_{2})}(u)}\right)
	\end{equation*}%
	by recalling that $V^{(i)}\in\mathbb{R}^{n\times \nu_i}$ for $\lambda_i\in\mathbb{R}$ and $V^{(i_{2})}=\overline{V^{(i_{1})}}$ for $\lambda_{i_1}\in\mathbb{C}$ with $\mathrm{Im}\left(\lambda_i\right)>0$. Since $u$ can be expressed as a linear combination of the Jordan
	basis in a unique manner, we obtain $\alpha ^{(i)}(u)=\overline{\alpha
		^{(i)}(u)}$, i.e., $\alpha^{(i)}(u)\in\mathbb{R}$, for $\lambda_{i}\in\mathbb{R}$ and $\alpha ^{(i_{2})}(u)=\overline{\alpha
		^{(i_{1})}(u)}$ for $\lambda _{i_{1}}\in\mathbb{C}$ with $\mathrm{Im}%
	\left( \lambda _{i_{1}}\right) >0$.
	
	Finally, we show that, for $i_{1},i_{2}\in\{1,\ldots ,p\}$ such that $\lambda
	_{i_{2}}=\overline{\lambda _{i_{1}}}$, we have $W^{(i_{2})}=\overline{%
		W^{(i_{1})}}$. This includes the case $i_{1}=i_{2}$ and $\lambda _{i_{1}}\in 
	\mathbb{R}$, for which we can conclude that $W^{(i_{1})}$ is a real matrix. We have, with $e^{\left( 1\right)},\ldots,e^{\left(
		n\right)} $ the real vectors of the canonical basis of $\mathbb{C}^n$, 
	\begin{eqnarray*}
		W^{(i_{2})} &=&W^{(i_{2})}\left[ e^{\left( 1\right) }\ \cdots\ e^{\left(
			n\right) }\right] =\left[ \alpha ^{\left( i_{2}\right) }\left( e^{\left(
			1\right) }\right) \ \cdots\ \alpha ^{\left( i_{2}\right) }\left(
		e^{(n)}\right) \right] \\
		&=&\left[ \overline{\alpha ^{\left( i_{1}\right) }\left( e^{\left( 1\right)
			}\right)} \ \cdots \ \overline{\alpha ^{\left( i_{1}\right) }\left( e^{(n)}\right)} \right]=\overline{\left[ \alpha ^{\left( i_{1}\right) }\left( e^{\left( 1\right)
			}\right) \ \cdots \ \alpha ^{\left( i_{1}\right) }\left( e^{(n)}\right) \right]
		}\\
		&=&\overline{W^{(i_{1})}\left[ e^{\left( 1\right) }\ \cdots \ e^{\left(
				n\right) }\right] }=\overline{W^{(i_{1})}}.
	\end{eqnarray*}
\end{proof}

\subsection{The formula for $\mathrm{e}^{tA}$} \label{AetA}

Next proposition provides the announced formula for the matrix exponential $\mathrm{e}%
^{tA}$. For sake of generality, we consider a
matrix function $f(A)$, where $f(z)$, $z\in\mathcal{D}\subseteq \mathbb{C}$,
is an analytic complex function of $z$. The domain $\mathcal{D}$ of $f$ is an open subset of $\mathbb{C}$ and we assume that the eigenvalues $\lambda_1,\ldots,\lambda_p$ of $A$ are contained in $\mathcal{D}$.

\begin{proposition}
	\label{matrixexponential} \label{one} We have 
	\begin{equation}
		f(A)=\sum\limits_{i=1}^{p}\sum\limits_{l=0}^{m_{i}-1}\frac{f^{(l)}(\lambda_i)%
		}{l!}P_{il},  \label{f(A)}
	\end{equation}%
	where, for $i\in\{1,\ldots ,p\}$ and $l\in\{0,\ldots,m_{i}-1\}$, 
	\begin{equation}
		P_{il}:=V^{\left( i\right) }N^{(i,l)}W^{\left( i\right) }\in\mathbb{C}%
		^{n\times n}  \label{Pil}
	\end{equation}%
	with $V^{\left( i\right) }$, $W^{\left( i\right) }$ and $N^{(i,l)}$
	defined in (\ref{this1}), (\ref{this}) and (\ref{A}), respectively.
\end{proposition}
\begin{proof}
	By recalling the definition of matrix function by the JCF (see, e.g., \cite%
	{Higham2008}), we have 
	\begin{equation}
		f(A)=Vf(J)W=\sum\limits_{i=1}^{p}V^{\left( i\right) }f(J^{(i)})W^{\left(
			i\right) },  \label{exp0}
	\end{equation}%
	where
	\begin{equation}
		f\left(J\right)=\mathrm{diag}\left( f\left(J^{(1)}\right),\ldots
		,f\left(J^{(p)}\right)\right) \in \mathbb{C}^{n\times n}\notag
	\end{equation}
	and 
	\begin{equation}
		f\left(J^{(i)}\right)=\mathrm{diag}\left( f\left(J^{(i,1)}\right),\ldots
		,f\left(J^{(i,d_i)}\right)\right) \in \mathbb{C}^{\nu _{i}\times \nu _{i}},\
		i\in\{1,\ldots,p\},  \notag
	\end{equation}%
	with 
	\begin{eqnarray*}
		&&f\left(J^{(i,j)}\right):=\left[ 
		\begin{array}{ccccc}
			f(\lambda_i) & f^\prime(\lambda_i) & \frac{f^{\prime\prime}(\lambda_i)}{2} & 
			\cdots & \frac{f^{(m_{ij}-1)}(\lambda_i)}{\left( m_{ij}-1\right) !} \\ 
			& . & . & . & \vdots \\ 
			&  & . & . & \frac{f^{\prime\prime}(\lambda_i)}{2} \\ 
			&  &  & . & f^\prime(\lambda_i) \\ 
			&  &  &  & f(\lambda_i)%
		\end{array}%
		\right] \in \mathbb{C}^{m_{ij}\times m_{ij}} \\
		&&j\in\{1,\ldots,d_i\}.
	\end{eqnarray*}
	We can write 
	\begin{equation*}
		f\left(J^{(i,j)}\right)=\sum\limits_{l=0}^{m_{ij}-1}\frac{f^{(l)}(\lambda_i)%
		}{l!}N^{(i,l,j)}=\sum\limits_{l=0}^{m_{i}-1}\frac{f^{(l)}(\lambda_i)%
		}{l!}N^{(i,l,j)},
	\end{equation*}%
	where the matrices $N^{(i,l,j)}$, $l\in\{0,\ldots ,m_{i}-1\}$, are defined in (\ref%
	{DA}).
	
	We conclude that, for\ $i\in\{1,\ldots,p\}$, we have 
	\begin{eqnarray*}
		f\left(J^{(i)}\right)&=&\mathrm{diag}\left( \sum\limits_{l=0}^{m_{i}-1}\frac{%
			f^{(l)}(\lambda_i)}{l!}N^{(i,l,1)},\ldots ,\sum\limits_{l=0}^{m_{i}-1}\frac{%
			f^{(l)}(\lambda_i)}{l!}N^{(i,l,d_{i})}\right) \\
		&=&\sum\limits_{l=0}^{m_{i}-1}\frac{f^{(l)}(\lambda_i)}{l!}N^{(i,l)}.
	\end{eqnarray*}%
	Now (\ref{f(A)}) follows by (\ref{exp0}).
\end{proof}

For the case $f(A)=\mathrm{e}^{tA}$ of our interest, formula (\ref{f(A)}) becomes 
\begin{equation*}
	\mathrm{e}^{tA}=\sum\limits_{i=1}^{p}\mathrm{e}^{\lambda
		_{i}t}\sum\limits_{l=0}^{m_{i}-1}\frac{t^{l}}{l!}P_{il}.  \label{etA}
\end{equation*}

\subsubsection{The matrices $P_{i0}$}  \label{matricesPi00}

For $i\in\{1,\ldots,p\}$, since $N^{(i,0)}=I_{\nu_i}$, we have
\begin{equation}
	P_{i0}=V^{\left( i\right) }W^{\left( i\right) }=\sum\limits_{j=1}^{d_i}\sum\limits_{k=1}^{m_{ij}}v^{(i,j,k)}w^{(i,j,k)}. \label{Pi0}
\end{equation}

The matrix $P_{i0}$ is the projection onto the generalized eigenspace corresponding to the eigenvalue $\lambda _{i}$, i.e., the subspace spanned by the generalized eigenvectors corresponding to the eigenvalue $\lambda_i$. In fact, we have, for $u\in\mathbb{C}^n$,
\begin{equation}
P_{i0}u=V^{\left( i\right) }W^{\left( i\right) }u=V^{\left( i\right) }\alpha^{\left( i\right) }(u)=\sum\limits_{j=1}^{d_i}\sum\limits_{k=1}^{m_{ij}}\alpha_{ijk}(u)v^{(i,j,k)}. \label{Pi0u}
\end{equation}

\subsubsection{The case $A$ diagonalizable} \label{diagonalizable}

When $A$ is diagonalizable, the formula (\ref{f(A)}) simplifies to the well-known formula
\begin{equation}
	f(A)=\sum\limits_{i=1}^{p}f(\lambda _{i})P_{i0}.  \label{fdiag}
\end{equation}
where, for $i\in\{1,\ldots,p\}$,
\begin{equation*}
	P_{i0}=\sum\limits_{j=1}^{d_i}v^{(i,j,1)}w^{(i,j,1)}.
\end{equation*}
In this case of $A$ diagonalizable, the matrix $P_{i0}$ is the projection onto the eigenspace corresponding to the eigenvalue $\lambda _{i}$, i.e., the subspace spanned by the eigenvectors corresponding to the  eigenvalue $\lambda_i$.

For the case $f(A)=\mathrm{e}^{tA}$ of our interest, formula (\ref{fdiag}) becomes
\begin{equation*}
	\mathrm{e}^{tA}=\sum\limits_{i=1}^{p}\mathrm{e}^{\lambda _{i}t}P_{i0}.
	\label{etAdiag}
\end{equation*}

\subsection{The matrices $P_{il}$} \label{matricesPi0}

In this subsection, we see some properties of the matrices $P_{il}$
defined in (\ref{Pil}).

\subsubsection{$P_{il}$ by a Cauchy integral}

The matrices $P_{il}$ are introduced in Subsection \ref{AetA} via the JCF of $A$. Alternatively, they can be expressed, independently of the JCF,  by a Cauchy integral, as stated in the next proposition.

\begin{proposition}
	For $i\in\{1,\ldots,p\}$ and $l\in\{0,\ldots,m_i-1\}$, we have
	\begin{equation}
		P_{il}
		=\frac{1}{2\pi\mathrm{i}}\int\limits_{\Gamma_i} (z-\lambda_i)^l(zI-A)^{-1}dz=(A-\lambda_iI)^lP_{i0}, \label{CI}
	\end{equation}
	where $\Gamma_i$ is a positively oriented simple closed contour in the complex plane enclosing $\lambda_i$ and excluding the other eigenvalues of $A$. 
\end{proposition}

\begin{proof}
	Let
	$$
	B_{il}
	=\frac{1}{2\pi\mathrm{i}}\int\limits_{\Gamma_i} (z-\lambda_i)^l(zI-A)^{-1}dz.
	$$
	Since $A=VJW$, we have, for $z\in\Gamma_i$,
	$$
	(zI-A)^{-1}=(V(zI-J)W)^{-1}=V(zI-J)^{-1}W
	$$
	with
	$$
	zI-J=\mathrm{diag}\left(zI_{\nu_{i_1}}-J^{(1)},\ldots,zI_{\nu_p}-J^{(d_i)}\right).
	$$
	Thus
	\begin{eqnarray}
		B_{il}&=&V\left(\frac{1}{2\pi\mathrm{i}}\int\limits_{\Gamma_i}(z-\lambda_i)^l(zI-J)^{-1}dz\right)W\notag \\
	&=&\sum\limits_{i^\prime=1}^p V^{(i^\prime)}\left(\frac{1}{2\pi\mathrm{i}}\int\limits_{\Gamma_i}(z-\lambda_i)^l\left(zI_{\nu_{i^\prime}}-J^{(i^\prime)}\right)^{-1}\right)W^{(i^\prime)}\notag \\
	&=&\sum\limits_{i^\prime=1}^p V^{(i^\prime)}\left(\frac{1}{2\pi\mathrm{i}}\int\limits_{\Gamma_i}(z-\lambda_i)^l\left((z-\lambda_{i^\prime})I_{\nu_{i^\prime}}-N^{(i^\prime,1)}\right)^{-1}dz\right)W^{(i^\prime)}\notag \\
	&&\text{(see Remark \ref{N})}\notag \\
	&=&\sum\limits_{i^\prime=1}^p V^{(i^\prime)}\left(\frac{1}{2\pi\mathrm{i}}\int\limits_{\Gamma_i}(z-\lambda_i)^{l}\sum\limits_{r=0}^{m_{i^\prime}-1}(z-\lambda_{i^\prime})^{-r-1}N^{(i^\prime,r)}dz\right)W^{(i^\prime)}\notag \\
	&&\text{(use Neumann series and see Remark \ref{N})}\notag \\
	&=&\sum\limits_{i^\prime=1}^p \left(\sum\limits_{r=0}^{m_{i^{\prime}}-1}\frac{1}{2\pi\mathrm{i}}\int\limits_{\Gamma_i}(z-\lambda_i)^{l}(z-\lambda_{i^\prime})^{-r-1}dz\right)V^{(i^\prime)}N^{(i^\prime,r)}W^{(i^\prime)}\notag \\		
	&=&\left(\sum\limits_{r=0}^{m_{i}-1}\frac{1}{2\pi\mathrm{i}}\int\limits_{\Gamma_i}(z-\lambda_i)^{l-r-1}dz\right)V^{(i)}N^{(i,r)}W^{(i)}\label{sumB}\\
	&&\text{(for $i^\prime\neq i$, the integral is zero since the integrand is analytic on and inside $\Gamma_i$)}\notag
	\end{eqnarray}
	For $k$ integer, by the residue theorem we have
	$$
	\frac{1}{2\pi\mathrm{i}}\int\limits_{\Gamma_i}(z-\lambda_i)^kdz=0
	$$
	unless $k=-1$, in which case it equals $1$. Hence, only the term with $l-r-1=-1$,
	i.e., $r=l$, is not zero in (\ref{sumB}). Therefore, we obtain
	$$
	B_{il}=V^{(i)}N^{(i,l)}W^{(i)}=P_{il}.
	$$

For the second equality, we prove that, 
\[
P_{il+1}=(A-\lambda_i I)P_{il},\ l\in\{0,\ldots,m_{i}-2\}.
\]
For $z\in\Gamma_i$, we have
\[
(A-\lambda_i I)(zI-A)^{-1}
=\bigl((z-\lambda_i)I-(zI-A)\bigr)(zI-A)^{-1}
=(z-\lambda_i)(zI-A)^{-1}-I
\]
and then, for $l\in\{0,\ldots,m_{i}-2\}$,
\[
(z-\lambda_i)^l(A-\lambda_i I)(zI-A)^{-1}
=(z-\lambda_i)^{l+1}(zI-A)^{-1}-(z-\lambda_i)^l I.
\]
Therefore, we obtain 
\begin{eqnarray}
(A-\lambda_i I)P_{il}&=&\frac{1}{2\pi \mathrm{i}}\int_{\Gamma_i}(z-\lambda_i)^l(A-\lambda_i I)(zI-A)^{-1}dz\notag \\
&=&\frac{1}{2\pi \mathrm{i}}\int\limits_{\Gamma_i}(z-\lambda_i)^{l+1}(zI-A)^{-1}dz \notag \\
&&-\frac{1}{2\pi \mathrm{i}}\int\limits_{\Gamma_i}(z-\lambda_i)^ldzI \label{twoints} \\
&=&P_{il+1} \notag
\end{eqnarray}
since the integral in (\ref{twoints}) is zero.
\end{proof}

\subsubsection{$P_{il}$ by right and left generalized eigenvectors}

Unpacking the definition (\ref{Pil}), we can express the matrices $P_{il}$  in terms of right and left generalized eigenvectors corresponding to the eigenvalue $\lambda_i$.

\begin{proposition}  \label{lastexample}
	For $i\in\{1,\ldots,p\}$ and $l\in\{0,\ldots,m_i-1\}$, we have 
	\begin{equation}
		P_{il}=\sum\limits
		_{\substack{ j=1  \\ m_{ij}\geq l+1}}^{d_{i}}\sum\limits_{k=1}^{m_{ij}-l}v^{%
			\left( i,j,k\right) }w^{\left( i,j,l+k\right) }.  \label{Pil=}
	\end{equation}
	Moreover, for $u\in\mathbb{C}^n$, we have
	\begin{equation}
		P_{il}u=\sum\limits
		_{\substack{ j=1  \\ m_{ij}\geq l+1}}^{d_{i}}\sum\limits_{k=1}^{m_{ij}-l}\alpha_{ij\,l+k}(u)v^{%
			\left( i,j,k\right) }.  \label{Pilu=}
	\end{equation}
\end{proposition}

\begin{figure}[tbp]
	\includegraphics[width=0.6\textwidth]{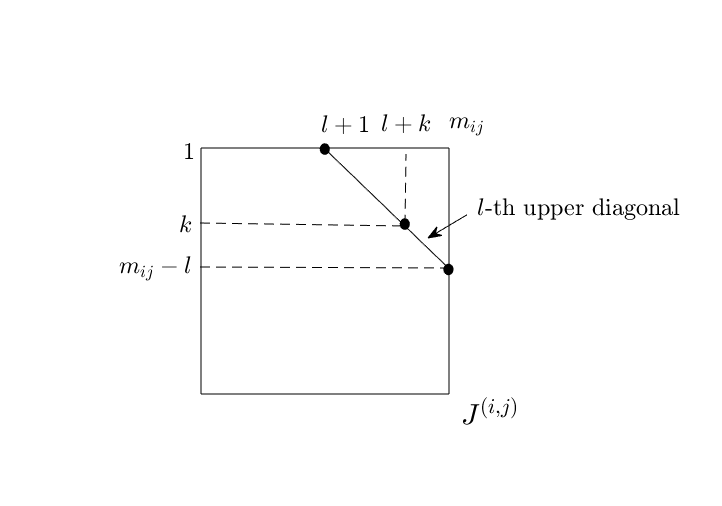}
	\caption{Indices $j$ and $k$ in (\ref{Pil=}) and (\ref{Pilu=}): $j$ is an index of mini-blocks $J^{(i,j)}$ of $J^{(i)}$ having the $l$-th upper diagonal and $k$ is an index of row of elements in the $l$-th upper diagonal.}
	\label{Pil==} 
\end{figure}

The indices $j$ in the outer sum
$$
\sum\limits
_{\substack{ j=1  \\ m_{ij}\geq l+1}}^{d_{i}}
$$
in (\ref{Pil=}) and (\ref{Pilu=}) are the indices $j$ of the mini-blocks $J^{(i,j)}$ of $J^{(i)}$ having the $l$-th upper diagonal. The indices $k$ in the inner sum 
$$
\sum\limits_{k=1}^{m_{ij}-l}
$$
are the row indices of the elements of the mini-block $J^{(i,j)}$ on this $l$-th upper diagonal (see Figure \ref{Pil==}). The column indices $l+k$ of these elements appear as third indices in $w^{(i,j,l+k)}$ in (\ref{Pil=}) and  $\alpha_{ij\,l+k}(u)$ in (\ref{Pilu=}).

In other words, with reference to the previous Figure \ref{matrixN}, in the double sum in (\ref{Pil=}) and (\ref{Pilu=}) we are summing over all elements in the intersection of the $l$-th diagonal of $N^{(i,l)}$ with the mini-blocks: the three-index notation $(i,j,k)$ of the row index of these elements appears in $v^{(i,j,k)}$ and the three-index notation $(i,j,l+k)$ of the column index appears in $w^{(i,j,l+k)}$ for  (\ref{Pil=}) and  in $\alpha_{ijl+k}(u)$ for (\ref{Pilu=}).

\begin{proof}
	Regarding (\ref{Pil=}), by the definition (\ref{Pil}) we have
	\begin{eqnarray*}
		P_{il}&=&\sum\limits_{j=1}^{d_{i}}V^{%
			\left( i,j\right) }N^{(i,l,j)}W^{\left( i,j\right) }=\sum\limits_{\substack{ j=1  \\ m_{ij}\geq l+1}}^{d_{i}}V^{\left( i,j\right) }N^{(i,l,j)}W^{\left( i,j\right)}\\
		&=&\sum\limits
		_{\substack{ j=1  \\ m_{ij}\geq l+1}}^{d_{i}}\sum\limits_{k=1}^{m_{ij}-l}v^{%
			\left( i,j,k\right) }w^{\left( i,j,l+k\right) }
	\end{eqnarray*}
	by recalling the form (\ref{DA}) of $N^{(i,l,j)}$: $N^{(i,l,j)}$ has entries $1$ at the indices $(k,l+k)$, $k\in\{1,\ldots,m_{ij}-l\}$, and entries $0$ at the other indices.

	Regarding (\ref{Pilu=}), we have
	\begin{equation*}
		P_{il}u=\sum\limits
		_{\substack{ j=1  \\ m_{ij}\geq l+1}}^{d_{i}}\sum\limits_{k=1}^{m_{ij}-l}v^{%
			\left( i,j,k\right) }\underset{=\alpha_{ij\,l+k}(u)}{\underbrace{w^{\left( i,j,l+k\right) }u}}.
	\end{equation*}

\end{proof}

\begin{remark}\label{mi}
	\quad
	\begin{itemize}

		\item [1.] For $l=0$, formulas (\ref{Pil=}) and (\ref{Pilu=}) reduce to formulas (\ref{Pi0}) and (\ref{Pi0u}), respectively.
		\item[2.] If there exists a unique $j\in\{1,\ldots,d_i\}$ such that $m_i=m_{ij}$, then
		\begin{equation*}
			P_{i\,m_i-1}=v^{%
				\left( i,j,1\right) }w^{\left( i,j,m_i\right) }
		\end{equation*}
		and, for $u\in\mathbb{C}^n$, 
		\begin{equation*}
			P_{i\,m_i-1}u=\alpha_{ijm_{i}}(u)v^{%
				\left( i,j,1\right) }.
		\end{equation*}
	\end{itemize}
	
\end{remark}

\subsubsection{The matrix $z A$}

The next proposition shows how the matrices $P_{il}$ are transformed when the matrix $A$ is replaced by $z A$, $z\in\mathbb{C}\setminus\{0\}$.
Observe that, when $A$ is replaced by $z A$, the number and the
dimensions of blocks and mini-blocks remain the same (see Proposition \ref%
{Vfor-A}). Thus, the indices $i$ and $l$ for the matrices $P_{il}(zA)$ corresponding to $z A$ are the same indices $i\in\{1,\ldots,p\}$ and $l\in\{0,\ldots,m_i-1\}$ for the matrices $P_{il}$ corresponding to $A$.

\begin{proposition}\label{PropalphalPil}
	Let $z\in\mathbb{C}\setminus\{0\}$. We have 
	\begin{equation*}
		P_{il}(zA)=z^{l}P_{il},\ i\in\{1,\ldots,p\}\text{\ and\ }l\in\{0,\ldots,m_{i}-1\}.  \label{alphalPil}
	\end{equation*}
\end{proposition}

\begin{proof}
	Recall Proposition \ref{PropWfor-A}. The columns of the matrix $V(zA)$ corresponding to $z A$ are given in (\ref{w-11}) and the rows of the matrix $W(zA)$ corresponding to $z A$ are given in (\ref{w-1}). Now, use the formula (\ref{Pil=}). 
\end{proof}
\subsubsection{The real case}

The next proposition describes the matrices $P_{il}$, when $A$ is real, in terms of real eigenvalues and complex conjugate pairs of eigenvalues of $A$.

\begin{proposition}
	\label{complex} Assume $A\in\mathbb{R}^{n\times n}$.  Moreover, assume we have a Jordan basis as that described in Proposition \ref{Prop22}.
	
	For $i\in\{1,\ldots ,p\}$ such that $\lambda _i\in\mathbb{R}$, we have%
	\begin{equation*}
		P_{il}\in\mathbb{R}^{n\times n},\ l\in\{0,\ldots,m_i-1\}.
	\end{equation*}
	
	For $i_{1},i_{2}\in\{1,\ldots ,p\}$ with $i_1\neq i_2$ such that $\lambda _{i_{2}}=%
	\overline{\lambda _{i_{1}}}$, we have%
	\begin{equation*}
		P_{i_{2}l}=\overline{P_{i_{1}l}},\ l\in\{0,\ldots ,m_{i_2}-1\}.
	\end{equation*}
	Observe that $m_{i_2}=m_{i_1}$ holds by Proposition \ref{Prop22}.
\end{proposition}
\begin{proof}
	We prove that, for $i_{1},i_{2}\in\{1,\ldots ,p\}$ such that $\lambda _{i_{2}}=%
	\overline{\lambda _{i_{1}}}$, we have%
	\begin{equation*}
		P_{i_{2}l}=\overline{P_{i_{1}l}},\ l\in\{0,\ldots ,m_{i_2}-1\}.
	\end{equation*}
	This includes the case $i_1=i_2$ and $\lambda_{i_1}\in\mathbb{R}$, for which
	we can conclude that $P_{i_1l}$, $%
	l\in\{0,\ldots,m_{i_1}-1\} $, is a real matrix.
	
	Since (see Proposition \ref{Prop22})
	\begin{equation*}
		d_{i_{2}}=d_{i_{1}},\ m_{i_{2}j}=m_{i_{1}j}\text{\ for\ }j\in\{1,\ldots ,d_{i_{1}}\},\text{\ and\ }m_{i_2}=m_{i_1}
	\end{equation*}%
	we obtain,  for $l\in\{0,\ldots,m_{i_2}-1\}$, 
	\begin{equation*}
		N^{(i_{2},l)}=N^{(i_{1},l)}
	\end{equation*}%
	and then (see Proposition \ref{conjVW}) 
	\begin{equation*}
		P_{i_{2}l}=V^{\left( i_{2}\right) }N^{(i_{2},l)}W^{\left( i_{2}\right)
		} =\overline{V^{\left( i_{1}\right) }}N^{(i_{1},l)}\overline{W
			^{\left( i_{1}\right) }}=\overline{V^{\left( i_{1}\right) }N^{(i_{1},l)}W^{\left( i_{1}\right) }} =\overline{P_{i_{1}l}}.
	\end{equation*}
	
\end{proof}

\subsubsection{The kernel of the matrices $P_{il}$}

Next three propositions concern the kernel of the matrices $P_{il}$. The first proposition describes this kernel in terms of components in the Jordan basis.

\begin{proposition}\label{eleven}
	\label{Prop1} For $i\in\{1,\ldots ,p\}$ and $l\in\{0,\ldots ,m_{i}-1\}$, we have%
	\begin{eqnarray*}
		&&\mathrm{ker}\left(P_{il}\right)=\left\{u\in\mathbb{C}^n:\alpha_{ijk}(u)=0\ 
		\text{for all }j\in\{1,\ldots ,d_{i}\}\ \text{with $m_{ij}\geq l+1$}\right. \\
		&&\ \ \ \ \ \ \ \ \ \ \ \ \ \ \ \ \ \ \ \ \ \ \ \ \ \left.\ \text{and for
			all }k\in\{l+1,\ldots ,m_{ij}\}\right\}.
	\end{eqnarray*}
\end{proposition}
\begin{remark}\label{REM1}
	With reference to the previous Figure \ref{matrixN} and \ref{Pil==}, Proposition \ref{eleven} states that the kernel of $P_{il}$ is constituted by the vectors $u\in\mathbb{C}^n$ with zero component $\alpha_{ijk}(u)$ along all the generalized eigenvectors $v^{(i,j,k)}$ such that, in the three-index notation, $(i,j,k)$ is an index column of elements in the intersection of the $l$-th diagonal of $N^{(i,l)}$ with mini-blocks.
\end{remark}     
\begin{proof}
	For $u\in\mathbb{C}^n$, by Proposition \ref{lastexample} we have 
	\begin{eqnarray*}
		P_{il}u=0 &\Leftrightarrow &\sum\limits
		_{\substack{ j=1  \\ m_{ij}\geq l+1}}^{d_{i}}\sum\limits_{k=1}^{m_{ij}-l}\alpha_{ij\,l+k}(u)v^{%
			\left( i,j,k\right) }=0 \\
		&\Leftrightarrow &\alpha _{ijk}\left( u\right) =0\ \text{for\ all }%
		j\in\{1,\ldots ,d_{i}\}\ \text{with $m_{ij}\geq l+1$}\\
		&&\text{and for all $k\in\{l+1,\ldots,m_{ij}\}$,}
	\end{eqnarray*}
where the second $\Leftrightarrow$ follows by the linear independence of the vectors $v^{(i,j,k)}$ in (\ref{JB}) of the Jordan basis.
\end{proof}
Next two propositions are immediate consequences of the previous proposition.
\begin{proposition}
	\label{Prop2A} For $i\in\{1,\ldots ,p\}$ and $l_{1},l_{2}\in\{0,\ldots ,m_{i}-1\}$, with $%
	l_{1}<l_{2}$, we have%
	\begin{equation*}
		\mathrm{ker}\left(P_{il_{1}}\right)\subseteq \mathrm{ker}\left(P_{il_{2}}%
		\right).
	\end{equation*}
\end{proposition}
\begin{proposition}
	\label{Prop3A} For $i\in\{1,\ldots ,p\}$ and $u\in \mathbb{C}^{n}$, we have%
	\begin{equation*}
		u\in\mathrm{ker}\left(P_{i0}\right)\ \ \Leftrightarrow \ \ \alpha ^{\left(
			i\right) }\left( u\right) =0.
	\end{equation*}
\end{proposition}
Proposition \ref{Prop3A} confirms our previous observation in Subsection \ref{matricesPi00} that $P_{i0}$ is the projection on the generalized eigenspace corresponding to $\lambda_i$. 

\subsection{The index $l_i(u)$}

\label{sliu}

In view of the Propositions \ref{Prop2A} and \ref{Prop3A}, for $i\in\{1,\ldots ,p\}$ and $u\in 
\mathbb{C}^{n}$ such that $\alpha^{(i)}(u)\neq 0$, we define the
index 
\begin{equation}
	l_{i}\left( u\right) :=\max \left\{ l\in \left\{ 0,\ldots ,m_{i}-1\right\}
	:u\notin \mathrm{ker}\left(P_{il}\right)\right\}.  \label{liu}
\end{equation}

Indeed, if  $\alpha^{(i)}(u)=0$, then
$$
\left\{ l\in \left\{ 0,\ldots ,m_{i}-1\right\}
:u\notin \mathrm{ker}\left(P_{il}\right)\right\}=\emptyset
$$
by Propositions \ref{Prop2A} and \ref{Prop3A}. In this case, the index $l_i(u)$ cannot be defined. On the other hand,  if  $\alpha^{(i)}(u)\neq 0$, then
$$\left\{ l\in \left\{ 0,\ldots ,m_{i}-1\right\}
:u\notin \mathrm{ker}\left(P_{il}\right)\right\}\neq\emptyset,
$$
since $u\notin \mathrm{ker}(P_{i0})$ by Proposition \ref{Prop3A}. In this case, the index $l_i(u)$ can be defined. By Proposition \ref{Prop2A}, we have
$$
u\notin\mathrm{ker}(P_{il}),\ l\in\{0,\ldots,l_i(u)\}
$$
and, by the definition (\ref{liu}) of $l_i(u)$, we have
$$
u\in\mathrm{ker}(P_{il}),\ l\in\{l_i(u)+1,\ldots,m_i-1\}.
$$

Observe that if $\lambda_i$ is a non-defective eigenvalue, i.e., $m_i=1$, and $\alpha^{(i)}(u)\neq 0$, then $l_{i}\left(
u\right) =0$.

 The following proposition relates $l_i(u)$ to the components of $u$ in the Jordan basis.

\begin{proposition} \label{number14}
	\label{propliu} For $i\in\{1,\ldots ,p\}$ and $u\in \mathbb{C}^{n}$ such that $%
	\alpha ^{\left( i\right) }\left( u\right) \neq 0$, we have 
	\begin{eqnarray}
		&&l_{i}\left( u\right)+1\notag\\
		&&=\max \left\{ k\in \left\{ 1,\ldots ,m_{i}\right\}
		:\alpha _{ijk}\left( u\right) \neq 0\text{ for some }j\in\{1,\ldots,d_{i}\}\text{\ with $m_{ij}\geq k$} \right\}.\notag
	\end{eqnarray}
\end{proposition}
In other words, the proposition states that $l_i(u)+1$ is the maximum index $k$ such that $u$ has non-zero component $\alpha_{ijk}(u)$ along the $k$-th generalized eigenvector $v^{(i,j,k)}$ of some Jordan chain
$$
(v^{(i,j,k)})_{k=1,\ldots,m_{ij}},\ j\in\{1,\ldots,d_i\},
$$
corresponding to the eigenvalue $\lambda_i$.

\begin{proof}
	Let
$$
k_{\max}=\max \left\{ k\in \left\{ 1,\ldots ,m_{i}\right\}
:\alpha _{ijk}\left( u\right) \neq 0\text{ for some }j\in\{1,\ldots,d_{i}\}\text{\ with $m_{ij}\geq k$} \right\}.
$$

	Let $l=l_i(u)+1$. Since $u\in\mathrm{ker}(P_{il})$, we have $\alpha_{ijk}(u)=0$, for any  $j\in\{1,\ldots,d_{i}\}$ and for any $k\in\{l+1,\ldots,m_{ij}\}$. Therefore $k_{\max}\leq l=l_i(u)+1$.
	
	Let $l=l_i(u)$. Since $u\notin\mathrm{ker}(P_{il})$, we have $\alpha_{ijk}(u)\neq 0$, for some  $j\in\{1,\ldots,d_{i}\}$ and for some $k\in\{l+1,\ldots,m_{ij}\}$. Therefore $k_{\max}\geq l=l_i(u)+1$.
	
	It follows that $k_{\max}=l_i(u)+1$.
\end{proof}

\subsection{Linear independence}

This subsection deals with the linear independence of the vectors $P_{il}u$
and the matrices $P_{il}$.

\begin{proposition}
	\label{PropPilu} \label{Prop} For $u\in\mathbb{C}^n$, the vectors 
	\begin{equation}
		P_{il}u,\ i\in\{1,\ldots,p\}\text{\ such that $\alpha^{(i)}(u)\neq 0$ and\ }%
		l\in\{0,\ldots,l_{i}\left( u\right)\},  \label{Pu}
	\end{equation}
	are linearly independent in the vector space $\mathbb{C}^n$.
\end{proposition}

\begin{proof}
	Consider a zero linear combination of the vectors (\ref{Pu}): 
	\begin{equation}
		0 =\sum\limits_{\substack{i=1 \\ \alpha^{(i)}(u)\neq 0}}^{p}\sum\limits_{l=0}^{l_{i}\left( u\right)
		}c_{il}P_{il}u=\sum\limits_{\substack{i=1 \\ \alpha^{(i)}(u)\neq 0}}^{p}V^{\left( i\right)
		}\sum\limits_{l=0}^{l_{i}\left( u\right) }c_{il}N^{(i,l)}\alpha
		^{\left( i\right) }\left( u\right),  \label{zlc}
	\end{equation}
	where the second equality follows by (\ref{Pil}). Since $V$ has linearly independent columns, we obtain, for any $i\in\{1,\ldots,p\}$ such that $\alpha^{(i)}(u)\neq 0$, 
	\begin{eqnarray*}
		0 &=&\sum\limits_{l=0}^{l_{i}\left( u\right) }c_{il}N^{(i,l)}\alpha
		^{\left( i\right) }\left( u\right) =\sum\limits_{l=0}^{l_{i}\left( u\right)
		}c_{il}\mathrm{diag}\left( N^{(i,l,1)},\ldots ,N^{(i,l,d_{i})}\right) \alpha
		^{\left( i\right) }\left( u\right) \\
		&=&\left( \sum\limits_{l=0}^{l_{i}\left( u\right) }c_{il}N^{(i,l,1)}\alpha
		^{\left( i,1\right) }\left( u\right) ,\ldots ,\sum\limits_{l=0}^{l_{i}\left(
			u\right) }c_{il}N^{(i,l,d_{i})}\alpha ^{\left( i,d_{i}\right) }\left(
		u\right) \right)
	\end{eqnarray*}%
	and then
	\begin{equation*}
		\sum\limits_{l=0}^{l_{i}\left( u\right) }c_{il}N^{(i,l,j)}\alpha ^{\left(
			i,j\right) }\left( u\right) =0,\ j\in\{1,\ldots,d_i\}.
	\end{equation*}%

	Now, fix $i\in\{1,\ldots,p\}$ such that $\alpha^{(i)}(u)\neq 0$. By Proposition \ref{propliu}, there exists $%
	\overline{j}\in\{1,\ldots,d_{i}\}$ with $m_{i%
		\overline{j}}\geq l_{i}\left( u\right)+1$  such that $\alpha _{i\overline{j}\,l_{i}\left(
		u\right)+1 }\left( u\right) \neq 0$ and $\alpha_{i\overline{j}k}(u)=0$ for $k\in\{l_i(u)+2,\ldots,m_{i%
		\overline{j}}\}$. Thus%
	\begin{eqnarray*}
		&&0=\sum\limits_{l=0}^{l_{i}\left( u\right) }c_{il}N^{(i,l,\overline{j}%
			)}\alpha ^{\left( i,\overline{j}\right) }\left( u\right) \\
		&&{\footnotesize =\left[ 
		\begin{array}{ccccccccc}
			c_{i0} & c_{i1} & \cdot & c_{il_{i}\left( u\right)-1 } & c_{il_{i}\left(
				u\right)} & | & 0 & \cdot & 0 \\ 
			0 & c_{i0} & c_{i1} & \cdot & c_{il_{i}\left( u\right)-1 } & | & 
			c_{il_{i}\left( u\right) } & \cdot & \cdot \\ 
			\cdot & \cdot & \cdot & \cdot & \cdot & | & c_{il_{i}\left( u\right)-1 } & 
			\cdot & 0 \\ 
			\cdot & \cdot & \cdot & c_{i0} & c_{i1} & | & \cdot & \cdot & 
			c_{il_{i}\left( u\right)} \\ 
			\cdot & \cdot & \cdot & \cdot & c_{i0} & | & c_{i1} & \cdot & 
			c_{il_{i}\left( u\right)-1 } \\ 
			- & - & - & - & - &  & - & - & - \\ 
			\cdot & \cdot & \cdot & \cdot & \cdot & | & c_{i0} & \cdot & \cdot \\ 
			\cdot & \cdot & \cdot & \cdot & \cdot & | & \cdot & \cdot & c_{i1} \\ 
			0 & \cdot & \cdot & \cdot & \cdot & | & \cdot & 0 & c_{i0}%
		\end{array}%
		\right] \left[ 
		\begin{array}{c}
			\alpha _{i\overline{j}1}\left( u\right) \\ 
			\alpha _{i\overline{j}2}\left( u\right) \\ 
			\\ 
			\cdot \\ 
			\alpha _{i\overline{j}l_{i}\left( u\right)}\left( u\right) \\ 
			\alpha _{i\overline{j}l_{i}\left( u\right)+1}\left( u\right) \\ 
			- \\ 
			0 \\ 
			\cdot \\ 
			0%
		\end{array}%
		\right]}
	\end{eqnarray*}
	and then 
	\begin{eqnarray*}
		&&0=\left[ 
		\begin{array}{ccccc}
			c_{i0} & c_{i1} & \cdot & c_{il_{i}\left( u\right)-1 } & c_{il_{i}\left(
				u\right)} \\ 
			0 & c_{i0} & c_{i1} & \cdot & c_{il_{i}\left( u\right)-1 } \\ 
			\cdot & \cdot & \cdot & \cdot & \cdot \\ 
			0 & 0 & \cdot & c_{i0} & c_{i1} \\ 
			0 & 0 & \cdot & 0 & c_{i0} \\ 
			&  &  &  & 
		\end{array}%
		\right] \left[ 
		\begin{array}{c}
			\alpha _{i\overline{j}1}\left( u\right) \\ 
			\alpha _{i\overline{j}2}\left( u\right) \\ 
			\cdot \\ 
			\alpha _{i\overline{j}l_{i}\left( u\right)}\left( u\right) \\ 
			\alpha _{i\overline{j}\,l_{i}\left( u\right)+1}\left( u\right) \\ 
		\end{array}%
		\right].
	\end{eqnarray*}
Since $\alpha _{i\overline{j}\,l_{i}\left( u\right)+1}\left( u\right)\neq
	0$, we obtain $c_{il}=0$ successively for $l=0,1,\ldots,l_{i}\left( u\right) $.
	
	We have proved that the coefficients of the zero linear combination (\ref{zlc}) are all zero. 
\end{proof}

\begin{proposition}
	\label{PropPil} The matrices
	\begin{equation}
		P_{il},\ i\in\{1,\ldots,p\}\text{\ and \ }l\in\{0,\ldots,m_{i}-1\},  \label{P}
	\end{equation}
	are linearly independent in the vector space $\mathbb{C}^{n\times n}$.
\end{proposition}

\begin{proof}
	Consider a zero linear combination of the matrices (\ref{P}): 
	\begin{equation*}
		0 =\sum\limits_{i=1}^{p}\sum\limits_{l=0}^{m_{i}-1 }c_{il}P_{il}.
	\end{equation*}
	Consider $u\in\mathbb{C}^n$ such that $\alpha_{ijk}(u)\neq 0$ for $%
	i\in\{1,\ldots,p\}$, $j\in\{1,\ldots,d_i\}$ and $k\in\{1,\ldots,m_{ij}\}$. For example, one can take $%
	u=Ve$, where $e=(1,\ldots,1)$, and then  $\alpha(u)=e$. We have $l_i(u)=m_i-1$
	for $i\in\{1,\ldots,p\}$ by Proposition \ref{number14}. Since we have
	\begin{equation*}
		0 =\sum\limits_{i=1}^{p}\sum\limits_{l=0}^{m_{i}-1 }c_{il}P_{il}u=\sum\limits_{i=1}^{p}\sum\limits_{l=0}^{l_{i}(u) }c_{il}P_{il}u
	\end{equation*}
	and $P_{il}u$, $i\in\{1,\ldots,p\}$ and $l\in\{0,\ldots,l_i(u)\}$, are linearly independent by the previous Proposition \ref{PropPilu}, we obtain $c_{il}=0$, $i\in\{1,\ldots,p\}$ and 
	$l\in\{0,\ldots,m_i-1\}$.
\end{proof}

\subsubsection{Linear combinations depending on $t$}

In the main paper, we deal with linear combinations of some of the vectors $P_{il}u$ in (\ref{Pu}), or some of the
matrices $P_{il}$ in (\ref{P}), whose coefficients depend on time $t$. For a fixed $t$, if these coefficients are not all zero, then
such a linear combination is not zero since the vectors and the matrices are
linearly independent: recall the two previous Propositions \ref{PropPilu} and \ref{PropPil}.

However, in the analysis of the asymptotic behavior of the condition numbers,
we require that the norm of the linear combination is away from zero,
uniformly with respect to $t$. This is the content of the next lemma.

\begin{lemma}
	\label{Lemma0} Let $\mathcal{V}$ be a vector space over $\mathbb{C}$ equipped
	with the norm $\left\Vert \ \cdot \ \right\Vert $, let $a_{1},\ldots ,a_{K}\in \mathcal{V}
	$ and let $f:I\rightarrow \mathbb{C}^{K}$, where $I$ is an
	arbitrary set. For any $t\in I$, the components $f_{1}(t),\ldots ,f_{K}(t)$ of $f(t)$ act as
	the coefficients in the linear combination
	$$
	\sum\limits_{k=1}^{K}f_{k}\left( t\right)
	a_{k}
	$$
	of $a_{1},\ldots ,a_{K}$.
	
	If $a_{1},\ldots ,a_{K}$ are linearly independent and%
	\begin{equation}
		\inf_{t\in I}\Vert f\left( t\right)\Vert_\infty >0, \label{last}
	\end{equation}%
	then%
	\begin{equation*}
		\inf_{t\in I}\left\Vert \sum\limits_{k=1}^{K}f_{k}\left( t\right)
		a_{k}\right\Vert >0. 
	\end{equation*}
\end{lemma}

\begin{proof}
	Suppose $a_1,\ldots,a_K$ linearly independent and (\ref{last}). Let 
	\begin{equation*}
		C=\left\{ z\in \mathbb{C}^{K}:\left\Vert z\right\Vert _{\infty }=1\right\} .
	\end{equation*}%
	Consider the function $g: \mathbb{C}^{K}\rightarrow \mathbb{R}$ given by%
	\begin{equation*}
		g\left( z\right) =\left\Vert \sum\limits_{k=1}^{K}z_{k}a_{k}\right\Vert ,\
		z\in \mathbb{C}^K.
	\end{equation*}%
	Since $C$ is a compact subset of $\mathbb{C}^{K}$ and $g$ is a continuous
	function, the extreme value theorem says that 
	\begin{equation*}
		m=\inf_{z\in C}g\left( z\right) =g\left( z^\ast\right)
	\end{equation*}%
	for some $z^\ast\in C$. Since $a_{1},\ldots ,a_{K}$ are linearly
	independent, we have $m>0$.
	
	Now, let $t\in I$. We have 
	\begin{equation*}
		\frac{f(t)}{\Vert f(t)\Vert_\infty}=\left(\frac{f_{1}\left( t\right) }{\Vert f(t)\Vert_\infty},\ldots,\frac{%
			f_{K}\left( t\right) }{\Vert f(t)\Vert_\infty }\right)\in C,
	\end{equation*}%
	which implies 
	\begin{equation*}
		g\left(\frac{f(t)}{\Vert f(t)\Vert_\infty}\right)=\left\Vert \sum\limits_{k=1}^{K}\frac{f_{k}\left( t\right) }{\Vert f(t)\Vert_\infty}a_{k}\right\Vert \geq m
	\end{equation*}%
	and then%
	\begin{equation*}
		\left\Vert \sum\limits_{k=1}^{K}f_{k}\left( t\right) a_{k}\right\Vert \geq
		m\Vert f(t)\Vert_\infty.
	\end{equation*}
	We conclude that 
	\begin{equation*}
		\inf_{t\in I}\left\Vert \sum\limits_{k=1}^{K}f_{k}\left( t\right)
		a_{k}\right\Vert\geq m\inf_{t\in I}\Vert f(t)\Vert_\infty >0.
	\end{equation*}
\end{proof}

In the main paper, we apply this lemma to the case where $a_1, \ldots, a_K$ are some of the linearly independent vectors $P_{il}u$ in (\ref{Pu}), or some of the linearly independent matrices $P_{il}$ in (\ref{P}), and $f_1(t),\ldots,f_K(t)$ are of the form $\mathrm{e}^{\mathrm{i}\omega t}$, with  $\omega\in\mathbb{R}$, and then $\Vert f(t)\Vert_\infty=1$.

\section{Properties of the matrices $Q_{jl}(t)$}  \label{lastsection}

In our study of asymptotic forms and asymptotic condition numbers, the following properties of the matrices $Q_{jl}(t)$ defined in (\ref{Qjl}) of the main paper are fundamental. We collect these properties in four propositions.

The first proposition relates the set $\Lambda_j(u)$ and the index $L_j(u)$, defined in Subsection \ref{asfetAu} of the main paper, to the condition $%
Q_{jl}(t)u=0$.

\begin{proposition}
	\label{lemma2Qjl} Let $u\in\mathbb{C}^n$. We have:
	
	\begin{itemize}
		\item[1)] $Q_{jl}(t)u=0$ for $j\in\{1,\ldots,q\}$ with $\Lambda_j(u)=\emptyset$
		and $l\in\{0,\ldots,L_j\}$;

		\item[2)] $Q_{jl}(t)u=0$ for $j\in\{1,\ldots,q\}$ with $\Lambda_j(u)\neq
		\emptyset$ and $l\in\{L_j(u)+1,\ldots,L_j\}$.
	\end{itemize}
\end{proposition}
\begin{proof}
	
	In 1) and 2), we prove that $P_{il}u=0$, for any $\lambda_i\in\Lambda_j$ with $m_i\geq l+1$ (see (\ref{Qjl}) of the main paper).
	
	Proof of point 1). Consider $j\in\{1,\ldots,q\}$ with $\Lambda_j(u)=\emptyset$
	and $l\in\{0,\ldots,L_j\}$. By Proposition \ref{Prop3A} in Appendix \ref{JCFsection}, $\Lambda_j(u)=\emptyset$ implies, for any $\lambda_i\in\Lambda_j$, $P_{i0}u=0$ and then, for any $\lambda_i\in\Lambda_j$ we have $P_{ik}u=0$, $k\in\{0,\ldots,m_i-1\}$, by Proposition \ref{Prop2A} in Appendix \ref{JCFsection}. In particular, for any $\lambda_i\in\Lambda_j$ with $m_i\geq l+1$, we have $P_{il}u=0$.

	Proof of point 2).  Consider $j\in\{1,\ldots,q\}$ with $\Lambda_j(u)\neq
	\emptyset$ and $l\in\{L_j(u)+1,\ldots,L_j\}$. For any $%
	\lambda_i\in\Lambda_j$ with $m_i\geq l+1$ and $P_{i0}u=0$, we have $P_{il}u=0$,  since $P_{ik}u=0$, $k\in\{0,\ldots,m_i-1\}$. For any $\lambda_i\in\Lambda_j$ with $m_i\geq l+1$ and $P_{i0}u\neq 0$, we also have $P_{il}u=0$, since $P_{ik}u=0$, $k\in\{l_i(u)+1,\ldots,m_i-1\}$, by the definition of $l_i(u)$ in Appendix \ref{JCFsection}, and $l\geq L_j(u)+1\geq l_i(u)+1$.
\end{proof}

The second proposition states that the matrices $Q_{jl}(t)$ and their actions
on vectors remain bounded and away from zero, by varying $t$.

\begin{proposition}
	\label{lemmaQjl} Let $j\in\{1,\ldots,q\}$, $l\in\{0,\ldots,L_j\}$ and $u\in\mathbb{C}^n$%
	. We have:
	
	\item[1)] $\sup\limits_{t\in \mathbb{R}} \left\Vert Q_{jl} (t)\right\Vert<+\infty$;
	
	\item[2)] 
	$\sup\limits_{t\in \mathbb{R}} \left\Vert Q_{jl} (t)u\right\Vert<+\infty$;
	
	\item[3)] 
	$
	\inf\limits_{t\in \mathbb{R}} \left\Vert Q_{jl} (t)\right\Vert>0$;

	\item[4)] $\inf\limits_{t\in \mathbb{R}} \left\Vert Q_{jl} (t)u\right\Vert>0$ if $\Lambda_j(u)\neq \emptyset$ and $l\leq L_j(u)$.
\end{proposition}

\begin{proof}
	The points 1) and 2) are trivial: by the definition (\ref{Qjl}) in the main paper, we have
	\begin{equation*}
		\left\Vert Q_{jl} (t)\right\Vert\leq \sum\limits_{\substack{ \lambda _{i}\in
				\Lambda _j  \\ m_{i}\geq l+1}}\left \Vert P_{il}\right\Vert\text{\ \ and\ \ }\left\Vert Q_{jl} (t)u\right\Vert\leq \sum\limits_{\substack{ \lambda
				_{i}\in \Lambda _j  \\ m_{i}\geq l+1}}\left \Vert P_{il}u\right\Vert. 
	\end{equation*}
	
	Lemma \ref{Lemma0} in Appendix \ref{JCFsection}, as applied to the linear combination (\ref{Qjl}), and Proposition \ref%
{PropPil} in Appendix \ref{JCFsection} imply 3).
	
	Finally, suppose $\Lambda_j(u)\neq \emptyset$ and $l\leq
	L_j(u)$. We obtain 
	\begin{equation}
		Q_{jl}(t)u=\sum\limits_{\substack{ \lambda _{i}\in \Lambda _j  \\ %
				m_{i}\geq l+1}}\mathrm{e}^{\mathrm{i}\omega _{i}t}P_{il}u=\sum\limits
		_{\substack{ \lambda _{i}\in \Lambda _j(u)  \\ m_{i}\geq l+1}}\mathrm{e}^{%
			\mathrm{i}\omega _{i}t}P_{il}u=\sum\limits_{\substack{ \lambda _{i}\in
				\Lambda _j(u)  \\ m_i\geq l+1 \\ l_i(u)\geq l}}\mathrm{e}^{\mathrm{i}\omega _{i}t}P_{il}u, \label{lc2}
	\end{equation}
	where the second equality holds since, for any $\lambda_i\in\Lambda_j$ with $m_i\geq l+1$ and $P_{i0}u=0$, i.e., $\alpha^{(i)}(u)=0$ by Proposition \ref{Prop3A} in Appendix \ref{JCFsection}, we have $P_{il}u=0$ by Propositions \ref{Prop2A} in Appendix \ref{JCFsection}; and the third equality holds since, for any $\lambda_i\in\Lambda_j$ with $m_i\geq l+1$, $P_{i0}u\neq 0$ and $l_i(u)<l$, we have $P_{il}u=0$ by the definition of $l_i(u)$ in Appendix \ref{JCFsection}. Now, Lemma \ref{Lemma0} in Appendix \ref{JCFsection}, as applied to the linear combination (\ref{lc2}), and Proposition \ref{PropPilu} in Appendix \ref{JCFsection} imply 4).
\end{proof}

The third proposition states how the matrices $Q_{jl}(t)$ are transformed when we
replace the matrix $A$ by $-A$. Observe that the matrix $-A$ has opposite eigenvalues, i.e., eigenvalues with opposite imaginary and real parts, with respect to the matrix $A$ and the dimensions of blocks and mini-blocks in the JCF of $-A$ are the same as in the JCF of $A$: see Proposition \ref{Vfor-A} in Appendix \ref{JCFsection} with $z=-1$. Therefore, we have the same number $q$ of different real parts for the eigenvalues of $-A$ and $A$. Moreover, the set $\Lambda_j(-A)$ and the numbers $r_j(-A)$ and $L_j(-A)$, $j\in\{1,\ldots,q\}$, corresponding to $-A$ are
\begin{equation}
	\Lambda_j(-A)=-\Lambda_{q+1-j},\ r_j(-A)=-r_{q+1-j}\text{\ \ and\ \ }L_j(-A)=L_{q+1-j}, \label{q+1-}
\end{equation}
where $\Lambda_{q+1-j}$, $r_{q+1-j}$ and $L_{q+1-j}$ correspond to $A$.   The indices $j$ and $l$ for the matrices $Q_{jl}(t,-A)$ corresponding to $-A$ range over  $j\in\{1,\ldots,q\}$ and $l\in\{0,\ldots,L_{q+1-j}\}$, respectively.
\begin{proposition}\label{-A}
	\label{-QA} We have 
	\begin{equation*}
		Q_{jl}(t,-A)=(-1)^lQ_{q+1-j l}(-t),\ j\in\{1,\ldots,q\}\text{\ and\ }l\in\{0,\ldots,L_{q+1-j}\}.  \label{alphalQjl}
	\end{equation*}
\end{proposition}
\begin{proof}
	The matrix $Q_{jl}(t,-A)$, $j\in\{1,\ldots,q\}$ and $l\in\{0,\ldots,L_{q+1-j}\}$, is given by
	\begin{equation*}
		Q_{jl}(t,-A)=\sum\limits_{\substack{ \lambda _{i}\in \Lambda _{q+1-j}  \\ %
				m_{i}\geq l+1}}\mathrm{e}^{\mathrm{i}(-\omega _{i})t}P_{il}(-A).
	\end{equation*}
	Now, use Proposition \ref{PropalphalPil} in Appendix \ref{JCFsection} with $z=-1$.
\end{proof}

In the case of a real matrix $A$, the fourth proposition shows how to rewrite the expression (\ref{Qjl}) in the main paper that defines $Q_{jl}(t)$ in terms of the real eigenvalues and complex conjugate pairs of eigenvalues of $A$. For
a matrix $Z$, we denote by $\mathrm{Re}(Z)$ and $\mathrm{Im}(Z)$ the matrices given by the real parts and imaginary parts, respectively, of the elements of $Z$.
\begin{proposition}
	\label{Pireal} Assume $A\in \mathbb{R}^{n\times n}$. For $j\in\{1,\ldots ,q\}$ and 
	$l\in\{0,\ldots,L_j\}$, we have 
	\begin{equation}
		Q_{jl}(t)=\sum\limits_{\substack{ \lambda _{i}\in \Lambda _{j}  \\ \lambda
				_{i}\mathrm{\ is\ real} \\ m_i\geq l+1}}P_{il} +2\sum\limits _{\substack{ %
				\lambda _{i}\in \Lambda _{j}  \\ \omega _{i}>0 \\ m_i\geq l+1}}\mathrm{Re}%
		\left( \mathrm{e}^{\mathrm{i}\omega _{i}t}P_{il} \right) .  \label{expd0}
	\end{equation}
\end{proposition}

\begin{proof}
	For $\lambda_i\in\Lambda_j$ such that $\lambda_i$ is real, we have $\omega_i=0$. For a complex conjugate pair $\lambda_{i_1},\lambda_{i_2}\in \Lambda_j$ with $\lambda _{i_{2}}=\overline{\lambda _{i_{1}}}$ and $\omega_{i_1}>0$,  by Proposition \ref{complex} in Appendix \ref{JCFsection} we have 
	\begin{eqnarray*}
		\mathrm{e}^{\mathrm{i}\omega _{i_{2}}t}P_{i_{2}l} &=&\mathrm{e}^{-\mathrm{i}\omega _{i_{1}}t}\overline{P_{i_{1}l}}=\overline{\mathrm{e}^{\mathrm{i}
				\omega _{i_{1}}t}P_{i_{1}l}}
	\end{eqnarray*}%
	and then 
	\begin{equation*}
		\mathrm{e}^{\mathrm{i}\omega _{i_{1}}t}P_{i_{1}l}+\mathrm{e}^{\mathrm{i}
			\omega _{i_{2}}t}P_{i_{2}l} =2\mathrm{Re}\left( \mathrm{e}^{\mathrm{i}\omega
			_{i_{1}}t}P_{i_{1}l} \right) .
	\end{equation*}
\end{proof}

\begin{remark}
	In (\ref{expd0}) the sum 
	\begin{equation*}
		\sum\limits_{\substack{ \lambda_i\in \Lambda _{j}  \\ \lambda _{i}\text{ is
					real} \\ m_i \geq l+1}} P_{il}
	\end{equation*}
	has zero or one term $P_{il}$, which is real: see Proposition \ref{complex} in Appendix \ref{JCFsection}. Moreover, in the other sum, each term
	\begin{equation*}
		\mathrm{Re}\left( \mathrm{e}^{\mathrm{i}\omega _{i}t}P_{il}
		\right)=\cos\omega_i t\cdot \mathrm{Re}(P_{il})-\sin\omega_i t\cdot \mathrm{%
			Im}(P_{il})
	\end{equation*}
	is a periodic function of $t$ of period $\frac{2\pi}{\omega_i}$.
\end{remark}

\section{The linear operators $Q_{j0}(t)|_{U_j}$ and $Q_{j0}(t)|_{U^e_j}$} \label{sssQj0}

The contents of this appendix are used in Subsections \ref{sss} and \ref {SoverKtA} of the main paper.  

For $j\in\{1,\ldots,q\}$, we have (see (\ref{Qjl}) in the main paper)
\begin{equation}
	Q_{j0}(t)=\sum\limits_{\lambda _{i}\in \Lambda _j }\mathrm{e}^{\mathrm{i}
		\omega _{i}t}P_{i0},  \label{Qj01}
\end{equation}
where $P_{i0}$ is the projection onto the generalized eigenspace corresponding to the eigenvalue $\lambda_i$ and it is given by (see Subsection \ref{matricesPi00} of Appendix \ref{JCFsection})
\begin{equation}
P_{i0}u=\sum\limits_{j^\prime=1}^{d_i}\sum\limits_{k=1}^{m_{ij}}\alpha_{ij^\prime k}(u)v^{(i,j^\prime,k)},\ u\in\mathbb{C}^n. \label{Pi0u1}
\end{equation}

Let $U_{j}$ be the sum of the generalized eigenspaces corresponding to eigenvalues in $\Lambda_j$ and let $U_{j}^e$ be the sum of the eigenspaces corresponding to eigenvalues in $\Lambda_j$. By (\ref{Pi0u1}), we see that, for $\lambda_i\in\Lambda_j$, 
$$
P_{i0}|_{U_j}:U_{j}\rightarrow U_{j}\text{\ \ and\ \ }P_{i0}|_{U^e_j}:U^e_{j}\rightarrow U^e_{j}
$$
and  then
$$
Q_{j0}(t)|_{U_{j}}:U_{j}\rightarrow U_{j}\text{\ \ and\ \ }Q_{j0}(t)|_{U^e_{j}}:U^e_{j}\rightarrow U^e_{j}.
$$
\begin{remark}\label{star}
When $\Lambda_j$ consists of a real eigenvalue $\lambda_i$, $U_j$ and $U_j^e$ are the generalized eigenspace and the eigenspace, respectively, corresponding to $\lambda_i$ and we have $Q_{j0}(t)=P_{i0}$ and 
$$
Q_{j0}(t)|_{U_{j}}=P_{i0}|_{U_j}=I_{U_{j}}\text{\ \ and\ \ }Q_{j0}(t)|_{U^e_{j}}=P_{i0}|_{U^e_j}=I_{U^e_{j}}.
$$
\end{remark}

The next proposition states that the linear operators $Q_{j0}(t)|_{U_{j}}$ and $Q_{j0}(t)|_{U^e_{j}}$ are invertible.
\begin{proposition}\label{Qj0}
	Let $j\in\{1,\ldots,q\}$. The linear operators $Q_{j0}(t)|_{U_{j}}:U_{j}\rightarrow U_{j}$ and $Q_{j0}(t)|_{U^e_{j}}:U^e_{j}\rightarrow U^e_{j}$ are invertible and their inverses are
	\begin{equation*}
		\left(Q_{j0}(t)|_{U_j}\right)^{-1}=Q_{j0}(-t)|_{U_{j}}\text{\ \ and\ \ }\left(Q_{j0}(t)|_{U^e_j}\right)^{-1}=Q_{j0}(-t)|_{U^e_{j}}.
	\end{equation*}
\end{proposition}

\begin{proof}
	We have 
	\begin{eqnarray*}
		Q_{j0}(-t)Q_{j0}(t)&=&\left(\sum\limits_{\lambda _{i}\in \Lambda _j }\mathrm{%
			e}^{-\mathrm{i}\omega _{i}t}P_{i0}\right)\left(\sum\limits_{\lambda _{k}\in
			\Lambda _j }\mathrm{e}^{\mathrm{i}\omega _{k}t}P_{k0}\right) \\
		&=&\sum\limits_{\lambda _{i},\lambda_k\in \Lambda _j }\mathrm{%
			e}^{\mathrm{i}(-\omega _{i}+\omega_{k})t}P_{i0}P_{k0}=\sum\limits_{\lambda _{i}\in \Lambda _j }P_{i0},
	\end{eqnarray*}
	since $P_{i0}P_{k0}=0$ for $\lambda_i\neq \lambda_k$ and $P_{i0}P_{k0}=P_{i0}$ for $\lambda_i=\lambda_k$. Thus,
	$$
	Q_{j0}(-t)|_{U_j}Q_{j0}(t)|_{U_j}=I_{U_j}\text{\ \ and\ \ }Q_{j0}(-t)|_{U_j^e}Q_{j0}(t)|_{U_j^e}=I_{U_j^e}.
	$$
\end{proof}
\begin{remark}
	\label{remQ00} 
As a consequence of the previous Proposition \ref{Qj0}, we have
		\begin{equation*}
			\min\limits_{\substack{\widehat{u}\in U_{j} \\ \left\Vert \widehat{u}\right\Vert=1}}\left\Vert Q_{j0}(t)%
			\widehat{u}\right\Vert=\frac{1}{\left\Vert Q_{j0}(-t)|_{U_{j}} \right\Vert}\text{\ \ and\ \ }\min\limits_{\substack{\widehat{u}\in U^e_{j} \\ \left\Vert \widehat{u}\right\Vert=1}}\left\Vert Q_{j0}(t)%
			\widehat{u}\right\Vert=\frac{1}{\left\Vert Q_{j0}(-t)|_{U^e_{j}} \right\Vert}. 
		\end{equation*}
\end{remark}

\begin{remark}\label{remQ0} 
	By (\ref{Qj01}), we see that
	\begin{equation}
	\sup\limits_{t\in\mathbb{R}}\left\Vert Q_{j0}(t)|_{U_j}\right\Vert<+\infty \text{\ \ and\ \ }\sup\limits_{t\in\mathbb{R}}\left\Vert Q_{j0}(t)|_{U^e_j}\right\Vert<+\infty. \label{supsup}
	\end{equation}
By the previous Remark \ref{remQ00}  and (\ref{supsup}) (as applied with time $-t$), we also see that
	$$
	\inf\limits_{t\in\mathbb{R}}\left\Vert Q_{j0}(t)|_{U_j}\right\Vert>0.\text{\ \ and\ \ }\inf\limits_{t\in\mathbb{R}}\left\Vert Q_{j0}(t)|_{U^e_j}\right\Vert>0.
	$$
\end{remark}

\begin{remark}
	For $j\in\{1,\ldots,q\}$ and $u\in\mathbb{C}^n$, the fact
	$$
	\inf\limits_{t\in\mathbb{R}}\Vert Q_{j0}(t)u\Vert>0\text{\ if\ }\Lambda_j(u)\neq \emptyset,
	$$
	included in 4) of Proposition \ref{lemmaQjl} of Appendix \ref{lastsection}, follows by the previous Remarks \ref{remQ00}  and \ref{remQ0}. In fact, if $\Lambda_j(u)\neq \emptyset$, we have
	$$
	v=\sum\limits_{\lambda_i\in\Lambda_j(u)}P_{i0}u\neq 0
	$$
	and then
	\begin{eqnarray*}
		\inf\limits_{t\in\mathbb{R}}\Vert Q_{j0}(t)u\Vert&=&\inf\limits_{t\in\mathbb{R}}\Vert Q_{j0}(t)v\Vert=\inf\limits_{t\in\mathbb{R}}\frac{\Vert Q_{j0}(t)v\Vert}{\Vert v\Vert}\Vert v\Vert\\
		&\geq&\inf\limits_{t\in\mathbb{R}}\min\limits_{\substack{\widehat{u}\in U_j\\ \Vert \widehat{u}\Vert =1}}\Vert Q_{j0}(t)\widehat{u}\Vert \Vert v\Vert=\inf\limits_{t\in\mathbb{R}}\frac{\Vert v\Vert}{\Vert Q_{j0}(-t)|_{U_j}\Vert}\\
		&=&\frac{\Vert v\Vert}{\sup\limits_{t\in\mathbb{R}}\Vert Q_{j0}(-t)|_{U_j}\Vert}>0.
	\end{eqnarray*}
\end{remark}

\end{document}